\documentclass[authoryear,12pt]{elsarticle}

\usepackage[latin1]{inputenc}
\usepackage{amsmath,amssymb,amsthm}
\usepackage{graphicx,graphics,color,epsfig}
\usepackage{hyperref}

\journal{Transportation Research C - 2011}

\graphicspath{{./figures/}}

\newcommand{\un}{\textbf{1}\normalsize}

\newcommand{\CU}{\mathcal U}

\newcommand{\BN}{\mathbb N}

\newcommand{\BR}{\mathbb R}

\theoremstyle{plain}
\newtheorem{theorem}{Theorem}
\newtheorem{corollary}{Corollary}

\theoremstyle{remark}
\newtheorem{remark}{Remark}

\begin{document}

\begin{frontmatter}

\title{The Traffic Phases of Road Networks}

\author[adr1]{Nadir Farhi\corref{cor1}}
\ead{nadir.farhi@ifsttar.fr}

\author[adr2]{Maurice Goursat}
\ead{maurice.goursat@inria.fr}

\author[adr2]{Jean-Pierre Quadrat}
\ead{jean-pierre.quadrat@inria.fr}

\address[adr1]{INRIA - Rh\^one Aples - Grenoble, 655 Avenue de l'Europe, 38330, \\ Montbonnot-Saint-Martin, France.}
\address[adr2]{INRIA - Paris - Rocqencourt, Domaine de Voluceau - B.P. 105,\\ 78153, Rocquencourt, France.}

\cortext[cor1]{Corresponding author.}

\begin{abstract}
We study the relation between the average traffic flow and the
vehicle density on road networks that we call 2D-traffic
fundamental diagram. We show that this diagram presents mainly
four phases. We analyze different cases. First, the case of a
junction managed with a priority rule is presented, four traffic
phases are identified and described, and a good analytic
approximation of the fundamental diagram is obtained by computing
a generalized eigenvalue of the dynamics of the system. Then, the
model is extended to the case of two junctions, and finally to a
regular city. The system still presents mainly four phases. The
role of a critical circuit of non-priority roads appears clearly
in the two junctions case. In Section 4, we use traffic light
controls to improve the traffic diagram. We present the
improvements obtained by open-loop, local feedback, and global
feedback strategies. A comparison based on the response times to
reach the stationary regime is also given. Finally, we show the
importance of the design of the junction. It appears that if the
junction is enough large, the traffic is almost not slowed down by
the junction.
\end{abstract}


\end{frontmatter}

\newpage

%-------------------------
\section{Introduction}
%-------------------------

The purpose of this paper is the discussion of the traffic phases
appearing in road networks. These phases appear clearly on what we
call the 2D traffic fundamental diagram.  This diagram gives the
relation between the average flow and the car density in the
network. The 2D traffic
fundamental diagram is an extension to road networks of the
well-known fundamental diagram studied mainly for a unique road in
traffic literature such as
in~\citep{DER,FUKUI,NS,BLA,WW,CSS,HEL,LMQ}.
This relation between the car density and the car flow has been
observed on a highway since 1935 by \cite{Green35}.

Microscopic traffic on a road network has been studied in
statistical physics with a cellular automata point of view
in~\citep{CSS,BML,Cues,Cues2,Cues3} similar to the modeling done
here. These works are mainly about numerical simulations showing a
threshold density where blocking appears for somewhat different
modeling, often stochastic. The case of two roads with one
junction without turning possibilities with a stochastic modeling
has been studied in detail in~\citep{FUKUI,FUKUI1,FUKUI2}.
In~\citep{BBSS}, the traffic light optimization of a city is done
with an automatic control point of view.

Recently, other attempts to extend this fundamental diagram
to the 2D cases have been done in~\citep{GD07,DG08,GD08,BL09,Hel09} in
traffic literature.

\citet*{DG08}, using a variational
traffic theory~\citep{Dag05,GD07,GD08}, have shown the existence of
a concave macroscopic fundamental diagram first on a ring, then on
a network using an aggregation method. Experimental
studies~\citep[e.g.][]{BL09} have shown that heterogeneity in traffic
measurements may have a strong impact on the shape of the
fundamental diagram built using the aggregation approach
of~\cite{DG08}.

More recently, \cite{Hel09}, using queuing theory, has derived
fundamental diagrams for a traffic model with junctions in the
saturated and unsaturated cases.

Here we propose a new 2D traffic model obtained from Petri
net modeling and maxplus algebra. This paper belongs to a sequence
of works presented in conferences \citep{CDC05,FAR2}, available
online \citep{Russ,ARX}, a thesis \citep{PHD}, and a companion paper
\citep{IEEE}. We give the first full synthesis of what
we have understood on the traffic phases. In
\citep{CDC05,FAR2,ARX,PHD,IEEE} the models are presented using
Petri nets and maxplus algebra formalisms. Here we avoid these
formalisms and use only dynamical systems in the standard algebra
framework. Our point of view is close to the cellular automata one
where interactions of a large number of simplified vehicles are
studied in asymptotic regimes. Our point of view is different from
the two approaches (Daganzo-Geroliminis and Helbing) discussed
previously. Our models are simpler but the systems can be analyzed
more deeply. In particular the fundamental diagram is obtained
(numerically or analytically) for the complete set of possible
densities. The different traffic phases can be deduced from the
obtained diagram instead of making a specific analysis for each
traffic situation (under or over saturated traffic), as for example,
in~\citep{Hel09}.

We consider various closed regular road networks cut in cells
which can contain at most one car. At a junction, a car leaving
can turn. To enter the junction, we consider the priority to the
right policy (a vehicle gives way to vehicles approaching from the
right at intersections) or the use of traffic lights with various
strategies. The roads are one-way streets without possibility of
overtaking. A car in a cell enters the next cell as soon as it is
freed. We observe through simulation or by mathematical analysis
(in the simplest case) four phases~:
\begin{itemize}
  \item \emph{Free phase}:  This phase corresponds to low densities. After a
    finite transitional time, the vehicles separate enough to be able to move freely. The average flow
    is equal to the density.
  \item \emph{Saturation phase}: This phase corresponds to average densities.
    The junctions serve with their maximal flow. There is no jam downstream of the junctions.
    The average flow is constant (independent of the density).
  \item \emph{Recession phase}: This phase corresponds to densities large enough
   to jam places downstream of the junctions. The average flow decreases with the density.
  \item \emph{Freeze phase}: This phase corresponds to high densities, where the
    number of vehicles is large enough to fill a circuit of
    roads. Once such a full circuit is constituted, the traffic is frozen and the average flow vanishes.
\end{itemize}

The simplest case studied is a system of two circular roads with
only one junction
managed with the priority to the right rule. In
this case, the dynamics is given. It is additively homogeneous of
degree one. We can prove the existence of an average flow which is
close to a generalized eigenvalue obtained analytically. We are
particularly interested in the dependence of these quantities with
respect to the density, but the ratio ($r$) between the number of
cells of the non-priority road and the number of cells of the
priority road also plays an important role in the determination of
the four phases.

Then, we study the case of four roads with two junctions.  We show
numerically that the fundamental diagram has the same shape as in
the previous case, but now $r$ has to be interpreted as the ratio
between the length of a non-priority road circuit and the length
of a priority road circuit. When we extend the analysis to a
larger network, the same kind of results persists.

The four phases exist in the presence of traffic light control,
but the length of the saturation phase can be enlarged and the
freeze phase reduced, almost without any impact on the free phase.
Moreover, we show that the number of cells of the junction plays
an important role. In our simple model, if we consider a junction
of two cells, instead of one cell, we suppress the saturation
phase and multiply by two the maximal flow.

Throughout this paper, we consider closed networks with a fixed
number of vehicles and wait until a stationary or
periodic regime is established\footnote{We have not seen clearly chaotic regime though
its presence is not excluded by the fact that the dynamics are
always additively homogeneous of degree 1 but not monotone.}
to compute the average flow as a function of
a given density. Nevertheless, the fundamental diagrams have a
sense for open systems when the density evolves slowly with
respect to the car speed. The diagram can be useful to describe
the behavior of zones of real towns with homogeneous structure of
roads and junctions.

Numerical simulations are done with the \emph{Scicoslab} software. They use the maxplus arithmetic of this software
very useful to describe large dynamical systems in a matrix may.
These simulations have played an important role to understand the different
phases of the system before being able to derive them sometimes analytically.

%------------------------------------------------------------------------------------
\section{Traffic modeling and analysis of two circular roads with one junction}
%------------------------------------------------------------------------------------
\label{sec2}

\begin{figure}[htbp]
  \begin{center}
    \includegraphics[width=60mm]{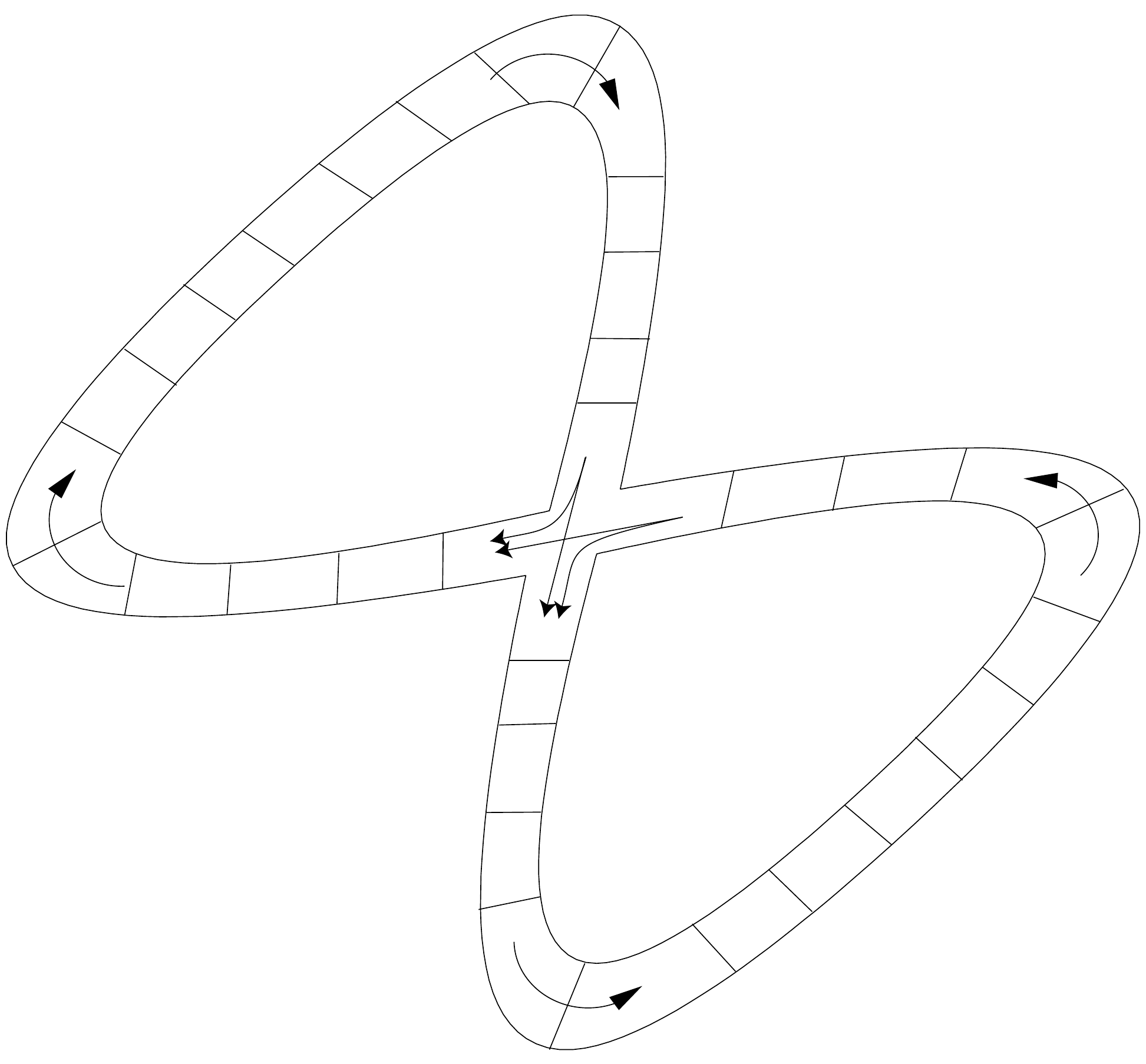}\hspace{7mm}
    \includegraphics[width=60mm]{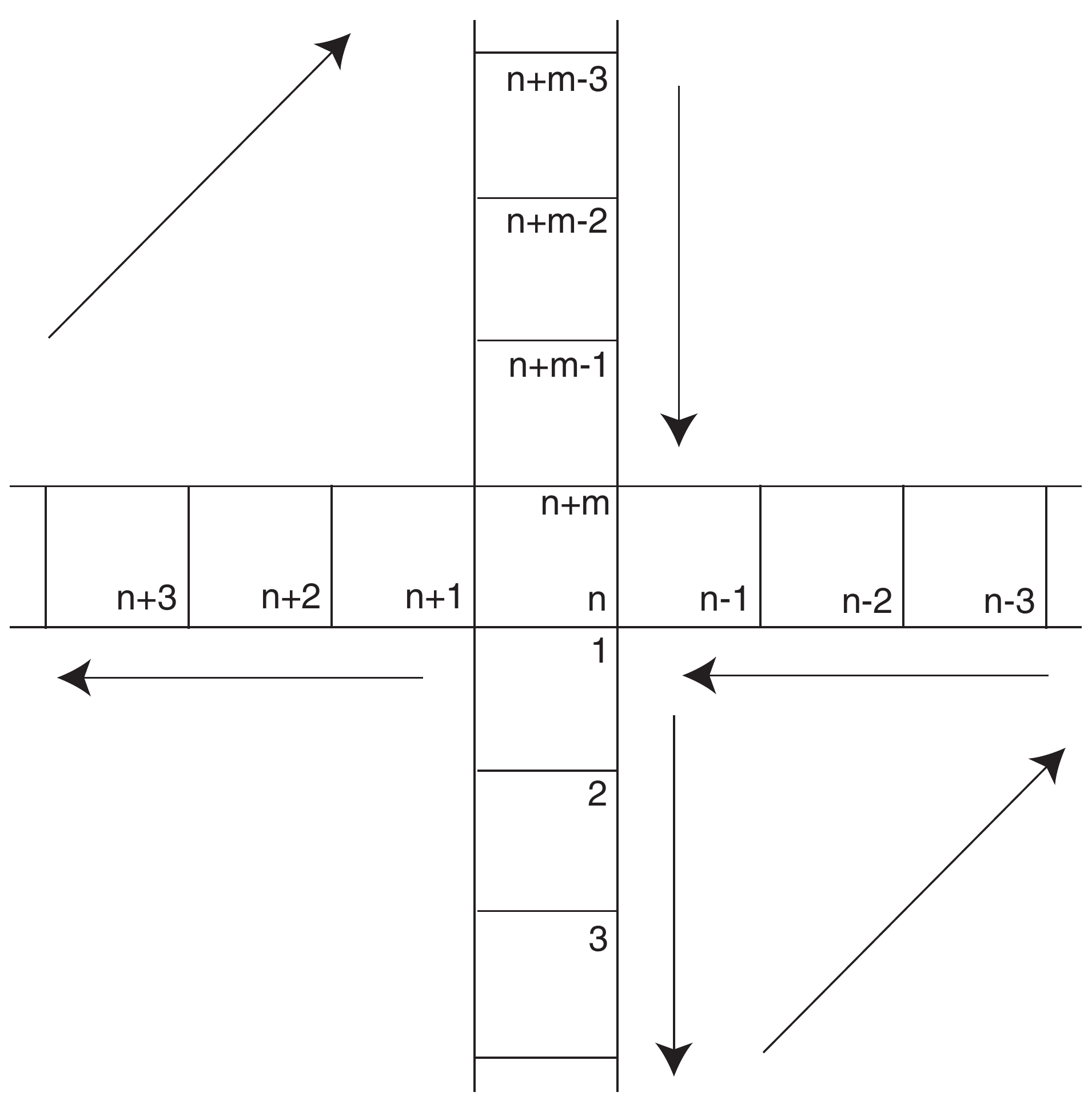}
    \caption{Two circular roads crossing on a junction with possibility of turning.}
    \label{crois}
  \end{center}
\end{figure}

We consider a circular traffic system of two roads with one
junction, in the shape of the numeral 8~\footnote{The roads being circular, the number of cars
in the system stays unchanged over time. This assumption allows us to maintain 
constant the car density.}, as shown on
Figure~\ref{crois}.
At the junction, the priority to the right
rule is applied to enter, and the vehicles have the possibility to
go straight or to turn. Thus, the northwest road has priority over
the southeast road.

Each road is cut in cells. We denote by $m$ and
$n$ the sizes (in number of cells) of the
priority and of the non-priority roads, respectively. The
cells are numbered. The junction is considered as a special
cell having two inputs and two outputs with two
sub-cells inside ($n$ and $n+m$ containing the vehicles going
respectively West and South); see Figure~\ref{crois}. Each
cell, including the junction, can contain at most
one vehicle. In a unit of time, on the roads outside the junction, a vehicle moves to the
cell ahead if the latter is free, and stays unmoved otherwise. At
the junction, the proportion of vehicles leaving and going west is
equal to the proportion of vehicles going south, and thus
equal to $1/2$. To enter the junction, the priority to the right
rule is applied. That is, a vehicle coming from the north has
priority over a vehicle coming from the east.

%-----------------------------
\subsection{The Dynamics}
%-----------------------------
\label{subsec2-1}

The car positions at the initial time (time $0$) are given by
numbers $a_i, 1\leq i\leq n+m$. For $i\neq n, n+m$, $a_i$ gives
the number of cars in the cell $i$ at time $0$. $a_n$ (resp.
$a_{n+m}$) gives the number of cars being in the junction at time
$0$, and that are going west (resp. south). Therefore the total
number of cars in the junction at time $0$ is $a_n+a_{n+m}$. To
describe the car dynamics on this system, we use  the variable
$x_i^k$ which denotes~:\\
-- for $i\neq n, n+m$, the cumulative car inflow into cell $i$ up
to time $k$, \\
-- for $i=n$, the cumulative inflow of vehicles coming from
the non-priority road (from cell $n-1$), into the junction
up to time $k$, \\
-- for $i=n+m$, the cumulative inflow of
vehicles coming from the priority road (from cell $n+m-1$), into
the junction  up to time $k$.
\begin{enumerate}
\item For each cell $i$, $i\neq 1,n, n+1, n+m$, the cumulative car
inflow into cell $i$ up to time $k+1$, denoted by $x_i^{k+1}$
satisfies:
\begin{itemize}
  \item $x_i^{k+1}$ is less than or equal to $a_{i-1}+x_{i-1}^k$, which
    is the number of vehicles present at cell $i-1$ at time~$0$, given
    by $a_{i-1}$, plus the cumulative car inflow into cell $i-1$ up to time~$k$,
    given by $x_{i-1}^k$.
  \item $x_i^{k+1}$ is less than or equal to $1-a_i+x_{i+1}^k$,
    which is the number of free places at cell $i$ at time $0$, given by
    $1-a_i$, plus the cumulative car outflow of cell $i$ (that is the cumulative car inflow into cell $i+1$ up to time $k$),
    given by $x_{i+1}^k$.
  \item $x_i^{k+1}$ is given by one of these two bounds since
    a car enters a cell as soon as possible.
\end{itemize}
Thus, we have:
\begin{equation}\label{eqxi}
  x_i^{k+1}=\min\{a_{i-1}+x_{i-1}^k, 1-a_i+x_{i+1}^k\}.
\end{equation}

\item For cell 1:
\begin{itemize}
  \item $x_1^{k+1}$ is less than or equal to $a_{n+m}+(x_n^k+x_{n+m}^k)/2$,
    which is the number of vehicles being in the junction and going south at time $0$, given
    by $a_{n+m}$, plus one half of the total cumulative car
    inflow into the junction up to time $k$, given by
    $(x_n^k+x_{n+m}^k)/2$.
  \item $x_1^{k+1}$ is less than or equal to $1-a_1+x_2^k$
    which is the number of free places in cell $1$ at time $0$, given by $1-a_1$,
    plus the cumulative car outflow of cell $1$ (that is the cumulative car inflow into cell $2$ up to time $k$),
    given by $x_2^k$.
  \item $x_1^{k+1}$ is given by one of these two bounds since
    a car enters a cell as soon as possible.
\end{itemize}
Hence:
\begin{equation}\label{c1}
  x_1^{k+1}=\min\{a_{n+m}+\frac{x_n^k+x_{n+m}^k}{2}\;, 1-a_1+x_2^k\}.
\end{equation}
If we want to maintain the integer property of counting events
and obtain a discrete version of equation~(\ref{c1}), we
could use the dynamics:
\begin{equation}\label{ceil}
  x_1^{k+1}=\min\{a_{n+m}+\left\lceil\frac{x_n^k+x_{n+m}^k}{2}\right\rceil\;,
    1-a_1+x_2^k\},
\end{equation}
where $\lceil\cdot\rceil$ is the rounding up operator to the
nearest integer. Equation
(\ref{ceil}), together with equation (\ref{floor}) given below,
tell that, with respect to the total counting of vehicles entered
the junction, the odd vehicles go south and the even vehicles go
west.
\item For cell~$n+1$, similarly with
cell~1, we obtain:
\begin{equation}\label{cn+1}
  x_{n+1}^{k+1}=\min\{a_n+\frac{x_n^k+x_{n+m}^k}{2}\;,
    1-a_{n+1}+x_{n+2}^k\}.
\end{equation}
\begin{equation}\label{floor}
  x_{n+1}^{k+1}=\min\{a_n+\left\lfloor\frac{x_n^k+x_{n+m}^k}{2}\right\rfloor\;,
    1-a_{n+1}+x_{n+2}^k\}.
\end{equation}
where $\lfloor\cdot\rfloor$ is the rounding down operator to the
nearest integer. 

\item For cell $n+m$:
\begin{itemize}
  \item $x^{k+1}_{n+m}$ is less than or equal to $a_{n+m-1}+x_{n+m-1}^k$,
    which is the number of vehicles present at cell $n+m-1$ at time $0$, given
    by $a_{n+m-1}$, plus the cumulative car inflow into cell $n+m-1$ up to time $k$,
    given by $x_{n+m-1}^k$.
  \item $x^{k+1}_{n+m}$ is less than or equal to
    $1-a_n-a_{n+m}+x_1^k+x_{n+1}^k-x_n^k$, which is the number of free places
    in the junction at time $0$, given by $1-a_n-a_{n+m}$,
    plus the cumulative car inflow into cell $1$ up to time $k$, given by $x_1^k$,
    plus the cumulative car inflow into cell $n+1$ up to time $k$, given by $x_{n+1}^k$,
    minus the cumulative car inflow into cell $n$ up to time $k$, given by
    $x_n^k$. Indeed, $1-a_n-a_{n+m}+x_1^k+x_{n+1}^k-x_n^k$ gives
    the number of entry authorizations by the
    North into the junction up to time $k$ (there is an authorization when the junction
    is free). The number $x_1^k+x_{n+1}^k$ gives the cumulative inflow of cars having left
    the junction up to time $k$. The number $x_n^k$ is subtracted here because entry authorizations to the
    junction are also used by vehicles entering from the East, and
    here, vehicles coming from the North have priority to use these authorizations.
  \item $x^{k+1}_{n+m}$ is given by one of these two bounds since
    a car enters a cell as soon as possible.
\end{itemize}
Thus:
\begin{equation}%\label{eqxn+m}
  x_{n+m}^{k+1}=\min\{a_{n+m-1}+x_{n+m-1}^k\;,
    1-a_n-a_{n+m}+x_1^k+x_{n+1}^k-x_n^k\}.
\end{equation}

\item For cell $n$:
\newline Similarly with cell $n+m$, we obtain~:
\begin{equation}%\label{eqxn}
  x_n^{k+1}=\min\{a_{n-1}+x_{n-1}^k\;, 1-a_n-a_{n+m}+x_1^k+x_{n+1}^k-x_{n+m}^{k+1}\}.
\end{equation}
where we have to subtract the
number of authorizations used by entering cars from the North
\emph{up to time $k+1$}, since these cars have priority.
\end{enumerate}
The dynamics of the whole system, in the continuous case
(equations (\ref{c1}) and (\ref{cn+1}) are considered rather than
equations (\ref{ceil}) and (\ref{floor})), is then given by:
\begin{align}
  & x_i^{k+1}=\min\{a_{i-1}+x_{i-1}^k, 1-a_i+x_{i+1}^k\},\;i\neq 1,n,n+1,n+m\;,\label{eqxi}\\
  & x_n^{k+1}=\min\{a_{n-1}+x_{n-1}^k\;, 1-a_n-a_{n+m}+x_1^k+x_{n+1}^k-x_{n+m}^{k+1}\}\; ,\label{eqxn}\\
  & x_{n+m}^{k+1}=\min\{a_{n+m-1}+x_{n+m-1}^k\;, 1-a_n-a_{n+m}+x_1^k+x_{n+1}^k-x_n^k\}\; ,\label{eqxn+m}\\
  & x_1^{k+1}=\min\{a_{n+m}+(x_n^k+x_{n+m}^k)/2\;, 1-a_1+x_2^k\}\;, \label{eqx1}\\
  & x_{n+1}^{k+1}=\min\{a_n+(x_n^k+x_{n+m}^k)/2\;, 1-a_{n+1}+x_{n+2}^k\}\;.\label{eqxn+1}
\end{align}
The parameters $a_i$ satisfy the constraints
\begin{equation}\label{const}
  \begin{cases}
    0\leq a_i\leq 1,\;  1\leq i\leq n+m\;,\\
    0\leq a_n+a_{n+m}\leq 1.
  \end{cases}
\end{equation}
The car density denoted by $d$ is given by
\begin{equation}\label{nota}
  d=\frac{1}{n+m-1}\sum_{i=1}^{n+m} a_i.
\end{equation}

In Table~\ref{simcon}, we show a short simulation of the
system (\ref{eqxi})--(\ref{eqxn+1}), with the parameters $n=m=5$,
$a=[0,1,0,1,0,1,0,0,1,0]$, and the initial condition $x^0=[0,0,0,0,0,0,0,0,0,0]$.

\begin{table}[htbp]
  \begin{center}
  \begin{tabular}{|c||c|c|c|c|c|c|c|c|c|c|}
    \hline
    Time & $x_1$ & $x_2$ & $x_3$ & $x_4$ & $x_5$ & $x_6$ & $x_7$ & $x_8$ & $x_9$ & $x_{10}$\\
    \hline
    \hline
    0 & 0 & 0 & 0 & 0 & 0 & 0 & 0 & 0 & 0 & 0\\
    \hline
    1 & 0 & 0 & 1 & 0 & 0 & 0 & 1 & 0 & 0 & 1\\
    \hline
    2 & 1/2 & 0 & 1 & 0 & 0 & 1/2 & 1 & 1 & 0 & 1\\
    \hline
    3 & 1/2 & 1/2 & 1 & 0 & 1 & 1/2 & 3/2 & 1 & 1 & 1\\
    \hline
    4 & 1 & 1/2 & 1 & 1 & 1 & 1 & 3/2 & 3/2 & 1 & 1\\
    \hline
    5 & 1 & 1 & 3/2 & 1 & 1 & 1 & 2 & 3/2 & 1 & 2\\
    \hline
  \end{tabular}
  \caption{A simulation of the system (\ref{eqxi})--(\ref{eqxn+1}) (continuous dynamics).}
  \label{simcon}
  \end{center}
\end{table}

In table~\ref{simdis}, we show a simulation of the
discrete version of this system; that is a simulation of
((\ref{eqxi}), (\ref{eqxn}), (\ref{eqxn+m}), (\ref{ceil}),
(\ref{floor})).

\begin{table}[htbp]
  \begin{center}
  \begin{tabular}{|c||c|c|c|c|c|c|c|c|c|c|}
    \hline
    Time & $x_1$ & $x_2$ & $x_3$ & $x_4$ & $x_5$ & $x_6$ & $x_7$ & $x_8$ & $x_9$ & $x_{10}$\\
    \hline
    \hline
    0 & 0 & 0 & 0 & 0 & 0 & 0 & 0 & 0 & 0 & 0\\
    \hline
    1 & 0 & 0 & 1 & 0 & 0 & 0 & 1 & 0 & 0 & 1\\
    \hline
    2 & 1 & 0 & 1 & 0 & 0 & 0 & 1 & 1 & 0 & 1\\
    \hline
    3 & 1 & 1 & 1 & 0 & 1 & 0 & 1 & 1 & 1 & 1\\
    \hline
    4 & 1 & 1 & 1 & 1 & 1 & 1 & 1 & 1 & 1 & 1\\
    \hline
    5 & 1 & 1 & 2 & 1 & 1 & 1 & 2 & 1 & 1 & 2\\
    \hline
  \end{tabular}
  \caption{A simulation of the system ((\ref{eqxi}), (\ref{eqxn}), (\ref{eqxn+m}), (\ref{ceil}),
    (\ref{floor})) (discrete dynamics).}
  \label{simdis}
  \end{center}
\end{table}

Using the discrete dynamics, we can compute the car
positions on the roads. Let us denote by $y_i^k$ the boolean
variable telling whether if there is a car in cell $i$ at time $k$
 ($y_i^k=1$) or not ($y_i^k=0$). The value $y_n^k$ (resp.
$y_{n+m}^k$) gives the presence of a car in the junction
going west (resp. south). The presence of a car
in the junction, no matter if it is going south or west, is given
by $y_n^k+y_{n+m}^k$. The values of $y_i^k, 1\leq i\leq n+m,
k\in\BN$ are deduced from the values of $x_i^k, 1\leq i\leq n+m,
k\in\BN$ as follows:
\begin{equation}\label{yi}
  \begin{cases}
    y_i^k=a_i+x_i^k-x_{i+1}^k, \quad i\neq n, n+m, &\\
    y_{n+m}^k=a_{n+m}+\lceil (x_n^k+x_{n+m}^k)/2\rceil -x_1^k, & \\
    y_n^k=a_n+\lfloor (x_n^k+x_{n+m}^k)/2\rfloor -x_{n+1}^k. &
  \end{cases}
\end{equation}
Note that, with the initial condition $x^0=0$ (which is natural
here since $x_i^k$ is a cumulative flow)
we have $y_i^0=a_i, \forall 1\leq i\leq n+m$. This
is also true for any initial condition $x^0$ satisfying
$x_p^0=x_q^0, \forall 1\leq p,q\leq n+m$. Indeed, this property
comes from the additive homogeneity of degree $1$ of the dynamical system, as it
is mentioned below. In Table~\ref{simveh}, we give the car
positions computed using~(\ref{yi}) where $x_i^k, 1\leq i\leq
n+m, \;\; 0\leq k\leq 5$ are the values given in
Table~\ref{simdis}.

\begin{table}[htbp]
  \begin{center}
  \begin{tabular}{|c||c|c|c|c||c|c||c|c|c|c|}
    \hline
     & \multicolumn{4}{c||}{Non-priority road} & \multicolumn{2}{c||}{Junction} & \multicolumn{4}{c|}{Priority road}\\
     \cline{6-7}
     & \multicolumn{4}{c||}{$\longrightarrow$} & \multicolumn{1}{c|}{West} & \multicolumn{1}{c||}{South} & \multicolumn{4}{c|}{$\longleftarrow$}\\
    \hline
    Time & $y_1$ & $y_2$ & $y_3$ & $y_4$ & $y_5$ & $y_{10}$ & $y_9$ & $y_8$ & $y_7$ & $y_6$\\
    \hline
    \hline
    0 & 0 & 1 & 0 & 1 & 0 & 0 & 1 & 0 & 0 & 1\\
    \hline
    1 & 0 & 0 & 1 & 1 & 0 & 1 & 0 & 0 & 1 & 0\\
    \hline
    2 & 1 & 0 & 1 & 1 & 0 & 0 & 0 & 1 & 0 & 0\\
    \hline
    3 & 0 & 1 & 1 & 0 & 1 & 0 & 1 & 0 & 0 & 0\\
    \hline
    4 & 0 & 1 & 0 & 1 & 0 & 0 & 1 & 0 & 0 & 1\\
    \hline
    5 & 0 & 0 & 1 & 1 & 0 & 1 & 0 & 0 & 1 & 0\\
    \hline
  \end{tabular}
  \caption{Car positions deduced from the simulation of the discrete dynamics.}
  \label{simveh}
  \end{center}
\end{table}

Let us remark that some car flow bounds can be obtained.
Outside the junction, at most one vehicle can pass through each
cell in two units of time (one time unit to enter the cell, one
time unit to leave the cell). Thus the flow is bounded by $1/2$ on
the roads. The junction having to serve the two roads it is
natural to think that the junction induces a bound of 1/4 on the
flow. This is justified in Theorem~\ref{thLAA} below.

To discuss the properties of the system
(\ref{eqxi})--(\ref{eqxn+1}), we recall some definitions. Let
$g:~\BR^p\to\BR^p$. The dynamical system $x^{k+1}=g(x^k)$ is
\emph{additively homogeneous of degree~$1$} if $g$ satisfies
$\forall x\in\BR^n, \forall\alpha\in\BR,
g(\alpha\un+x)=\alpha\un+g(x),$ where we denote
$\un={}^t(1,1,...,1)$. The dynamical system $x^{k+1}=g(x^k)$ is
\emph{monotone} if $g$ satisfies $\forall x,y\in\BR^p, x\leq
y\Rightarrow g(x)\leq g(y).$ It admits a \emph{growth rate} if the
vector $\chi\in\BR^p$, defined by $\chi=\lim_{k\to +\infty}x^k/k$,
exists. It admits an \emph{additive eigenvalue} if there exists
$x\in\BR^p$ such that $\lambda + x=g(x)$. When $\chi_i=\chi_j,
\forall 1\leq i,j\leq p$, we may also say, by abuse of language,
that the growth rate $\chi$ is any of these components.

The dynamical system (\ref{eqxi})--(\ref{eqxn+1}) is additively
homogeneous of degree $1$ but not monotone (because of
equations~(\ref{eqxn}) and~(\ref{eqxn+m})). It is known
\citep{GG04} that if an additively homogeneous system is monotone
and connected\footnote{A dynamical system $x^{k+1}=g(x^k)$, where
$g: \BR^p\to\BR^p$, is connected if: $\forall 1\leq i,j\leq p,
\lim_{\alpha\to\infty}g_i(\alpha e_j)=\infty$, where $e_j$ is the
$j^{th}$ vector of the canonical basis of $\BR^p$.}, then the
average growth rate of the dynamical system is unique (independent
of the initial state of the system), and its components are equal
and coincide with its unique additive eigenvalue. In our case, we
have the following result (see \citep{IEEE,PHD} for the proof).
\begin{theorem}\label{thLAA}
   The system~(\ref{eqxi})-(\ref{eqxn+1}) with the initial condition $x^0=0$ have the following properties~:
   \begin{enumerate}
     \item The trajectories of the states are non-negative and non-decreasing.
     \item The distances between any pair of states stay bounded~: \\
       $\exists c\in\BR, \sup_{k} |x_i^k-x_j^k| \leq c, \forall i,j$.
     \item If a growth rate $\chi$ exists, then $\chi\leq 1/4$.
   \end{enumerate}
\end{theorem}

From Theorem~\ref{thLAA} and from the property of homogeneity, it
is easy to apply the ergodic theorem and prove the existence of a
growth rate; see \citep{IEEE}. Moreover, the eigenvalue problem can
be solved explicitly as a function of the vehicle density;
see \citep{PHD, ARX}. In our case this eigenvalue is not
equal to the growth rate, but we can derive a good approximation
(justified by numerical simulations) of the growth rate from its
analytic expression. 

%----------------------------------------------------------
\subsection{The eigenvalue as a function of the density}
%----------------------------------------------------------
\label{subsec2-2}
We are interested in this subsection by solving the additive eigenvalue
problem associated to the dynamical system~(\ref{eqxi})--(\ref{eqxn+1}).
The objective of determining the eigenvalue is to compare it to the
growth rate of the dynamical system, which is interpreted in terms of traffic
as the average car flow. We will see that, although the eigenvalue
is different from the average car flow, it gives a good approximation of it.
This approximation helps a lot to identify and understand all the phases
of traffic.

The eigenvalue problem associated to the dynamical
system~(\ref{eqxi})--(\ref{eqxn+1}) is the computation of the
couple $(\lambda,x)$ solution of:
\begin{align}
  & \lambda+x_i=\min\{a_{i-1}+x_{i-1}, 1-a_i+x_{i+1}\},\quad i\neq 1,n,n+1,n+m\;,\label{veqxi}\\
  & \lambda+x_n=\min\{a_{n-1}+x_{n-1}\;, 1-a_n-a_{n+m}+x_1+x_{n+1}-(\lambda+x_{n+m})\}.\label{veqxn}\\
  & \lambda+x_{n+m}=\min\{a_{n+m-1}+x_{n+m-1}\;, 1-a_n-a_{n+m}+x_1+x_{n+1}-x_n\}\;,\label{veqxn+m}\\
  & \lambda+x_1=\min\{a_{n+m}+(x_n+x_{n+m})/2\;, 1-a_1+x_2\}\;,\label{veqx1}\\
  & \lambda+x_{n+1}=\min\{a_n+(x_n+x_{n+m})/2\;, 1-a_{n+1}+x_{n+2}\}\;.\label{veqxn+1}
\end{align}

\begin{figure}[htbp]
  \begin{center}
    \includegraphics[width=4cm,height=3.5cm]{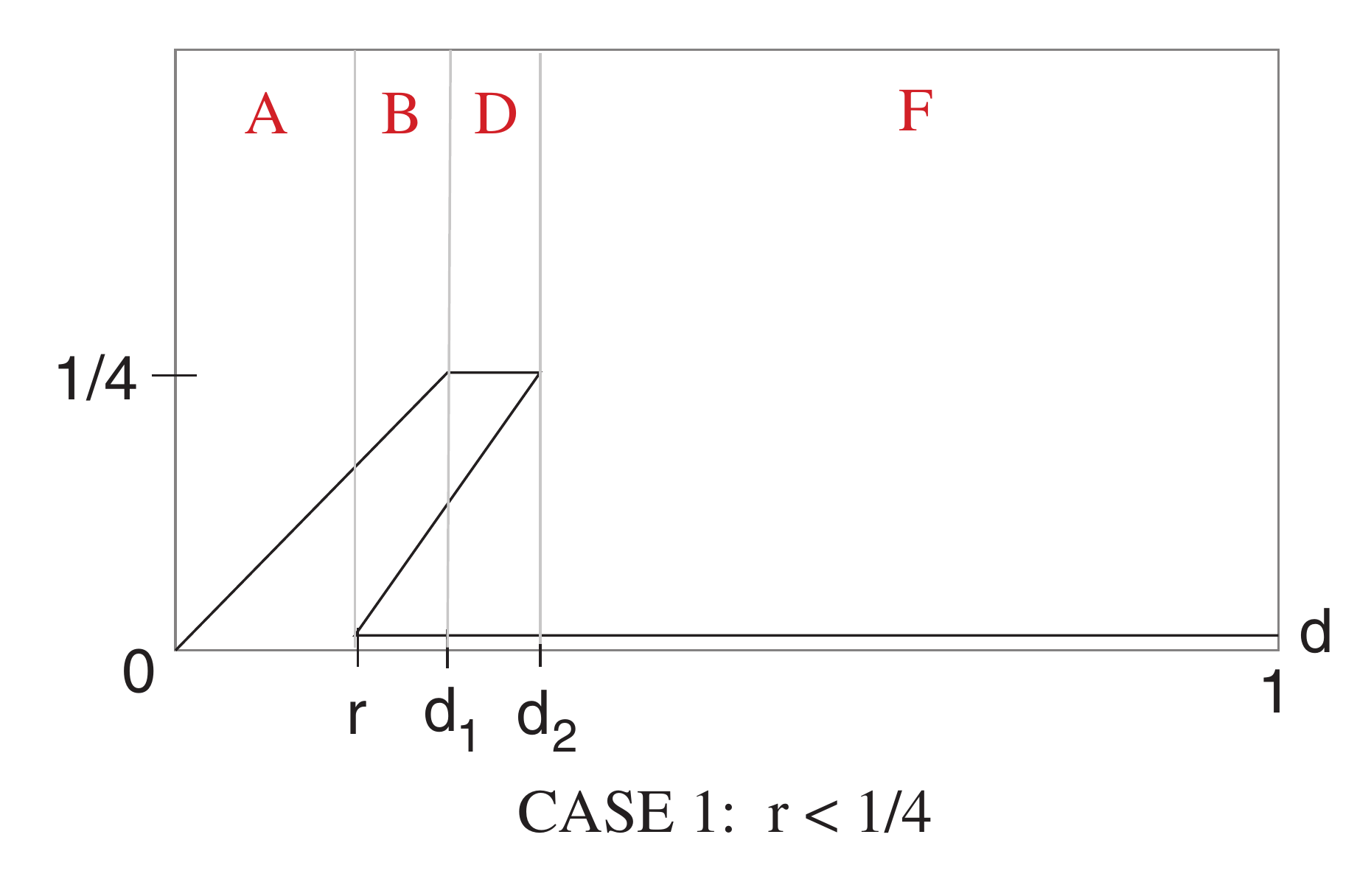}\hspace{5mm}
    \includegraphics[width=4cm,height=3.5cm]{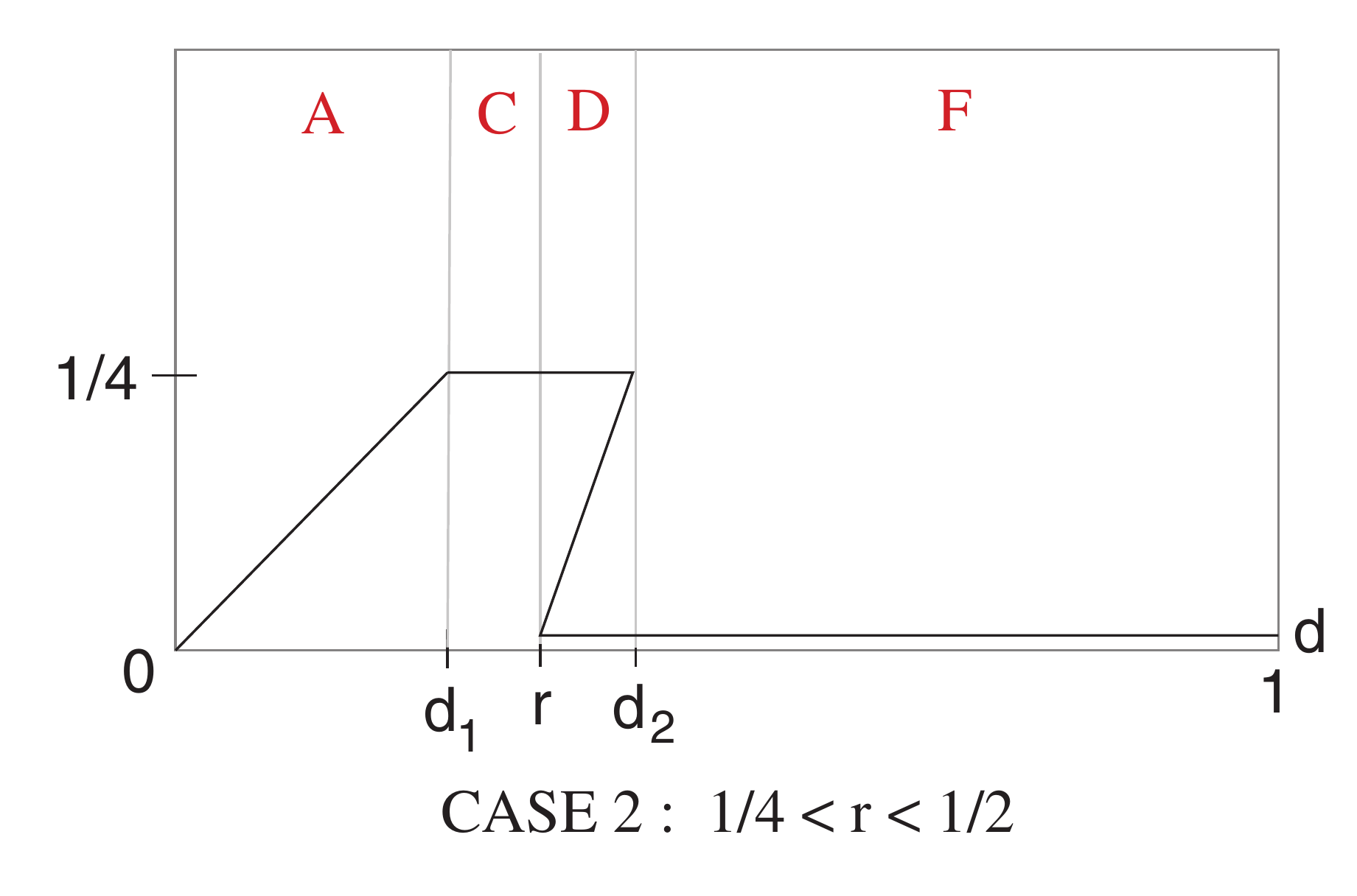}\hspace{5mm}
    \includegraphics[width=4cm,height=3.5cm]{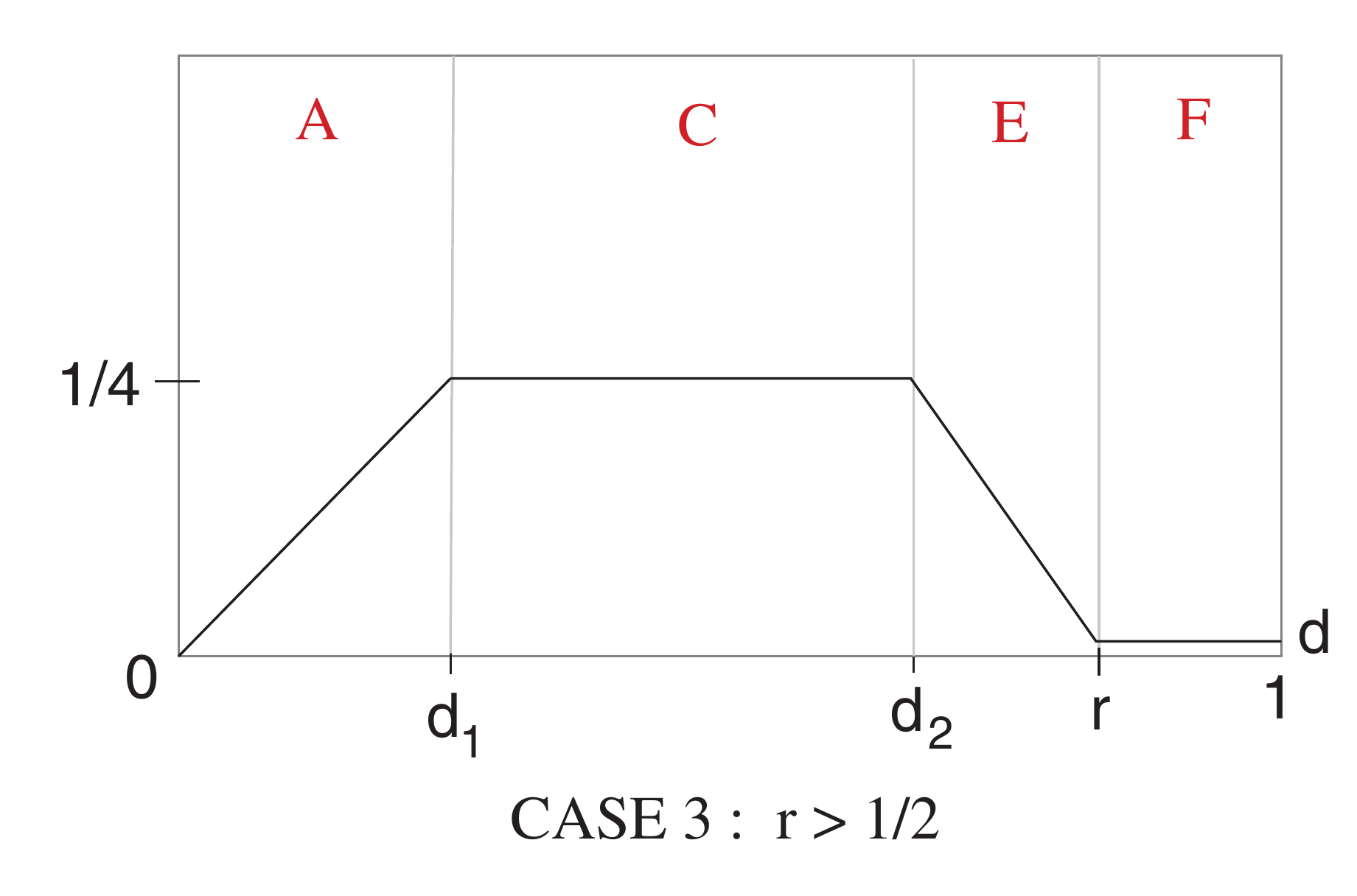}
    \caption{The curve of $\lambda$ given in Theorem~\ref{eigen} depending on $d$.}
    \label{diag}
  \end{center}
\end{figure}

Let us use the notations $r=n/(n+m-1)$ and $\rho=1/(n+m-1)=r/n$.
The following results are proved in \citep{ARX}.
\begin{theorem}\label{eigen}
  There exists a non-negative additive eigenvalue $\lambda$, solution
  of the system (\ref{veqxi})--(\ref{veqxn+1}), satisfying:
  $$0=\max\left\{\min \left\{d-(1+\rho)\lambda,\;\frac{1}{4}-\lambda,\;
      r-d-\left(2r-1+\rho\right)\lambda\right\},\;-\lambda \right\}.$$
\end{theorem}
\begin{remark} Using the notations:
\begin{itemize}
  \item[] $d_1=(n+m)/[4(n+m-1)]=(1+\rho)(1/4)\;,$
  \item[] $d_2=(3n+m-2)/[4(n+m-1)]=(2r+1-\rho)/4\;,$
\end{itemize}
we can explain the eigenvalues given by Theorem~\ref{eigen} as
follows (see Figure~\ref{diag}):
  \begin{itemize}
    \item If $0\leq d\leq d_1$ then
       $\lambda=d/(1+\rho)\;$,
    \item If $d_1\leq d\leq d_2$ then
      $\lambda=1/4\;$,
    \item If $d_2< d\leq r$ or $r\leq d< d_2$ which cases
      correspond respectively to $r>1/2$ or $r<1/2$ then
     $\lambda=(r-d)/(2r-1+\rho)\;$,
    \item If $r\leq d\leq 1$  then $\lambda=0\;$.
  \end{itemize}
\end{remark}
\begin{corollary} \label{fond-cor1}
  In the case $r\geq 1/2$, the non negative eigenvalue $\lambda$ solution of the system (\ref{veqxi})--(\ref{veqxn+1}) is
  given by:
  $$\lambda=\max\left\{\min\left\{\frac{1}{1+\rho}\;d\;,\frac{1}{4}\;,\frac{r-d}{2r-1+\rho}\right\},0\right\}.$$
\end{corollary}
%\proof follows directly from Theorem~\ref{eigen}\finprv
\begin{corollary} \label{fond-cor2}
  For large values of $n$ and $m$ such that $n>m-2$ (which is the case $r\geq 1/2$),
  a non negative eigenvalue $\lambda$ solution of the system (\ref{veqxi})--(\ref{veqxn+1}) is given by:
  $$\lambda=\max\left\{\min\left\{d,\;\frac{1}{4},\;\frac{r-d}{2r-1}\right\},0\right\}.$$
\end{corollary}
%\proof follows directly from Corollary~\ref{fond-cor2}\finprv
\begin{remark}\label{rem1}
  We can check that as soon as we assume $m>1$, we get $d_1<d_2$.
  The position of $r$ with respect to $d_1$
  and $d_2$ gives three cases shown in Figure~\ref{diag}:
  \begin{itemize}
    \item[A.] $d\in[0,\min(d_1,r)) \quad \Rightarrow \quad \lambda=\frac{1}{1+\rho}\;d$,
    \item[B.] $d\in[\min(d_1,r),d_1) \quad \Rightarrow \quad
      \begin{cases}
        & \lambda=\frac{1}{1+\rho}\;d,\\
        \text{or} &
        \lambda=\frac{r-d}{2r-1+\rho},\\
        \text{or} & \lambda=0.
      \end{cases}$,
    \item[C.] $d\in[d_1,\min(d_2,r)) \quad \Rightarrow \quad \lambda=\frac{1}{4}$,
    \item[D.] $d\in[\max(d_1,r),d_2) \quad \Rightarrow \quad
      \begin{cases}
        & \lambda=\frac{1}{4},\\
        \text{or} &
        \lambda=\frac{r-d}{2r-1+\rho},\\
        \text{or} & \lambda=0.
      \end{cases}$,
    \item[E.] $d\in[d_2,\max(d_2,r)) \quad \Rightarrow \quad \lambda=\frac{r-d}{2r-1+\rho}$,
    \item[F.] $d\in[\max(d_2,r),1] \quad \Rightarrow \quad \lambda=0$.
  \end{itemize}
\end{remark}

%----------------------------------------------
\subsection{The traffic fundamental diagram}
%----------------------------------------------
\label{subsec2-3}

\begin{figure}[htbp]
  \begin{center}
    \includegraphics[width=62mm]{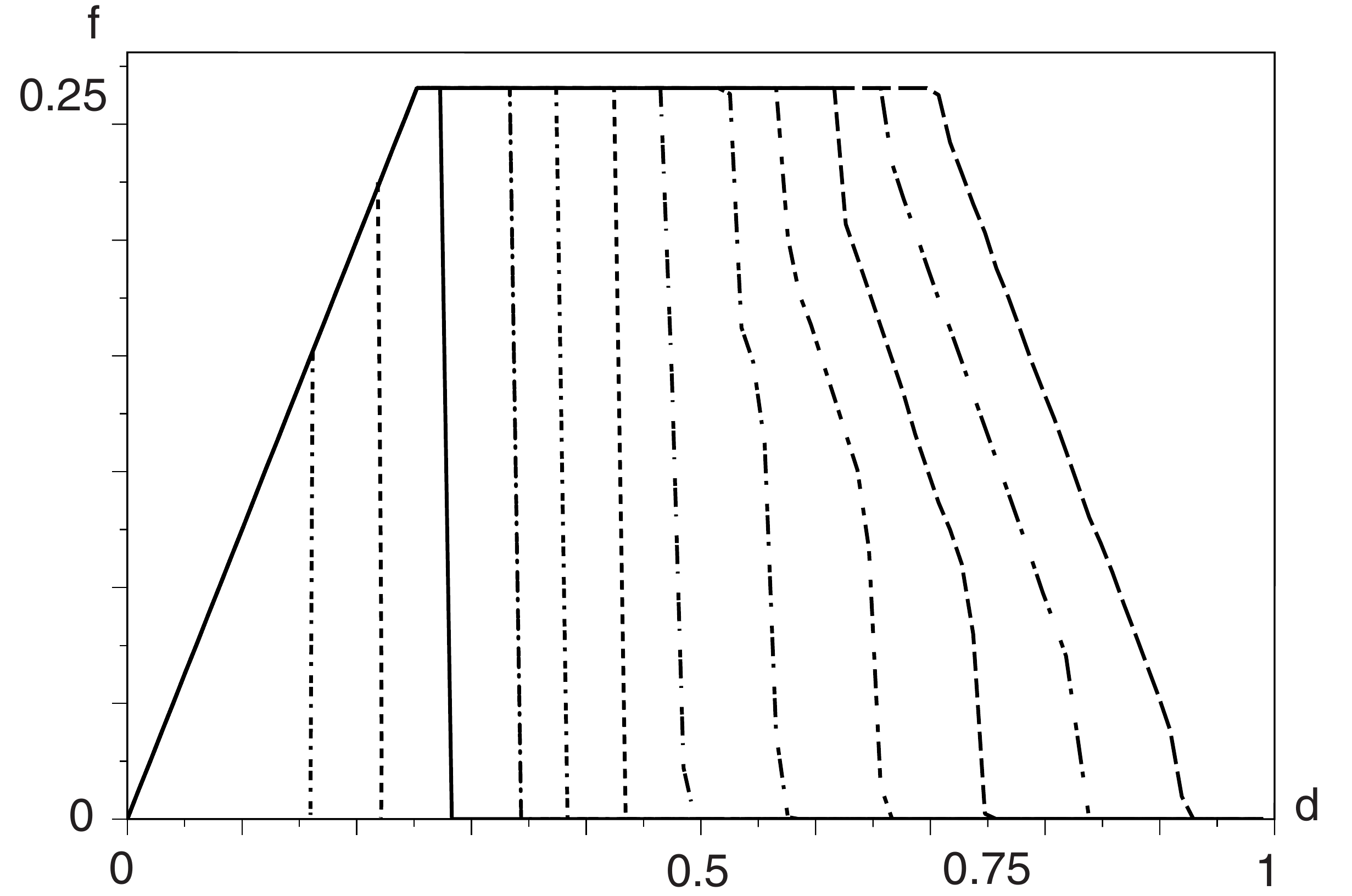}\hspace{7mm}
    \includegraphics[width=57mm]{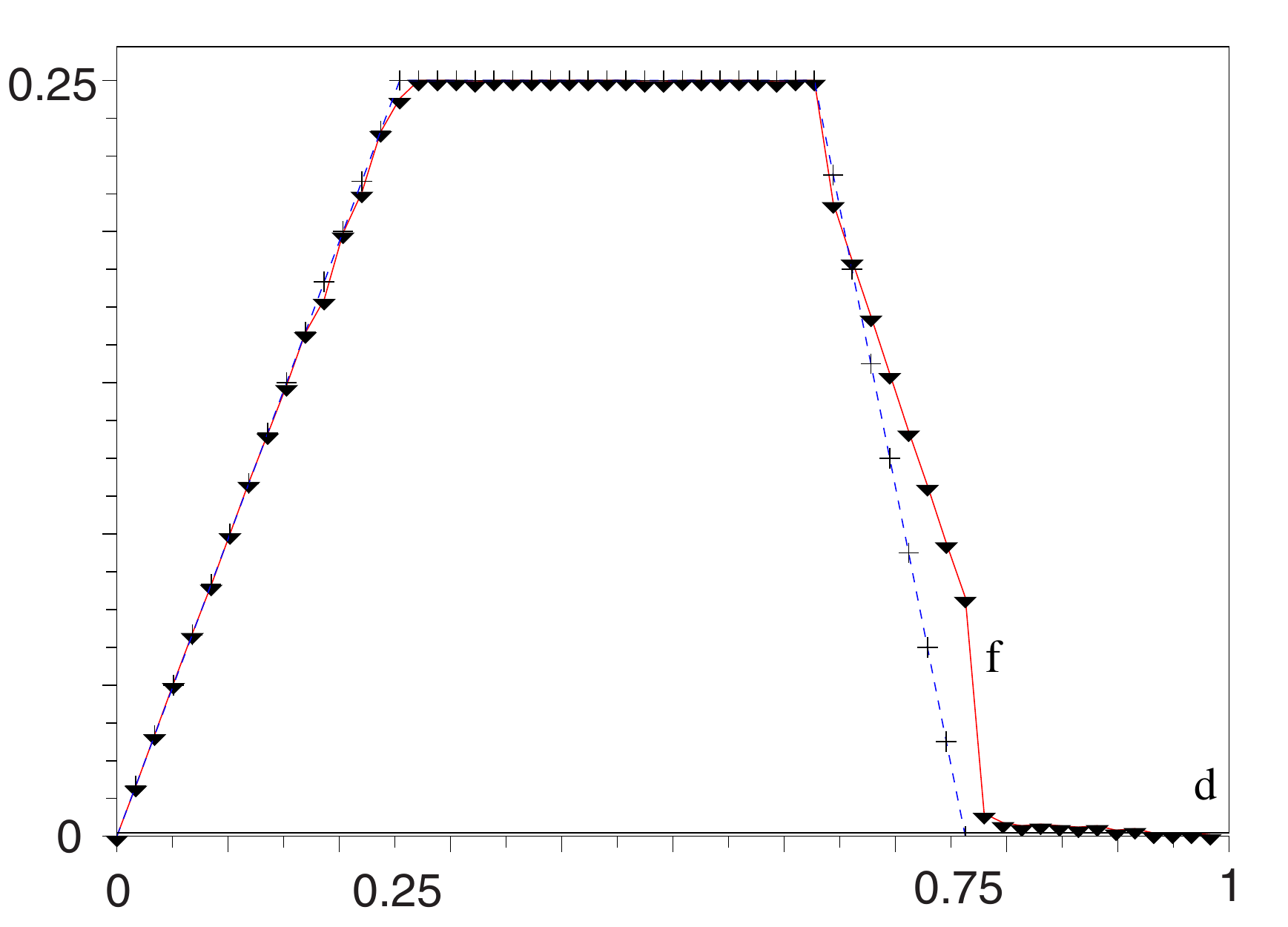}
    \caption{On the left side: The fundamental diagram of 2D-traffic
      depending on the ratio $r$, where $r$ takes the values: 0.15, 0.20,
      0.25, 0.35, 0.40, 0.45, 0.50, 0.60, 0.70, 0.75, 0.85 and 0.95,
      respectively from the left to the right. On the right side: Comparison of $f$
      and $\lambda$ in the case $n=45, m=15$, i.e. $r=3/4>1/2$.}
    \label{turn-45-15}
  \end{center}
\end{figure}

The system (\ref{eqxi})--(\ref{eqxn+1}) admits a growth rate
$f=\lim_{k\to+\infty}x^k/k$ for trajectories starting from $0$.
This quantity has the interpretation of the average car flow (for
this reason, we denote it by $f$ rather than $\chi$). It is
difficult to obtain an analytical expression for $f$.
Nevertheless, it is easy to compute an approximation by numerical
simulation. Choosing a $K$ large enough, we have $f\simeq x^K/K$.
Clearly this quantity depends on the car density and on the car
positions at initial time. However, we have remarked
by simulation that the dependence with the
initial car positions is not very important and disappears when
the size of the system (size of the roads and number of cars)
increases. Asymptotically, the average car flow depends only on
the density $d$ and on the relative sizes of the roads $r$.
Moreover, since the system is not monotone, the eigenvalue is not
equal to the growth rate (see \citep{IEEE}). However, we see
by simulation that the average car flow
$f$ is always close to one of the
eigenvalues $\lambda$. The left side of
Figure~\ref{turn-45-15} shows numerical simulations for different
values of $r$. This figure tells that when the eigenvalue is not
unique, the average car flow chooses the eigenvalue zero (see also
Figure~\ref{diag2-num1} below). In the right side of
Figure~\ref{turn-45-15}, we see that $\lambda$ and $f$ are very
close to each other when $\lambda$ is unique (i.e. when $r>1/2$).
This figure tells that $f$ and $\lambda$ differ mainly in
one phase among from the four traffic phases (see the traffic
phases analysis in section~\ref{subsec2-4} below).

In the left side of Figure~\ref{turn-45-15}, we give the
dependence of $f$ with $d$ and $r$. For a given $r$, the function
$f(d)$ is a generalization to 2D-traffic systems (presence of
junctions) of what is called \emph{fundamental diagram} in traffic
literature.

The eigenvalue, when it is unique, gives a good analytical
approximation of the average flow.  When it is not unique, we have
to choose the good one (in this case the eigenvalue 0). Finally, we obtain a good
approximation of the flow by the formula~:
%which is in the case of big sizes of the roads (to compare to the
%formula given in Corollary~\ref{fond-cor2}):
\begin{equation}\label{for2}
  f=\max\left\{\min\left\{d,\;\frac{1}{4},\;\frac{r-d}{\max\{2r-1,\;0\}}\right\},0\right\},
\end{equation}
where $a/0=\text{sign}(a)\infty$. Therefore when $r\leq 1/2$
we have $\max\{2r-1,0\}=0$, then, for $d<r$ we obtain $f=\min\{d,1/4\}$ and for $d\geq r$
we obtain $f=0$.

In Figure~\ref{diag2-num1} (to be compared with
Figure~\ref{diag}), we give a graphic summary of
the different possible cases.

\begin{figure}[htbp]
  \begin{center}
    \includegraphics[width=4cm,height=3.5cm]{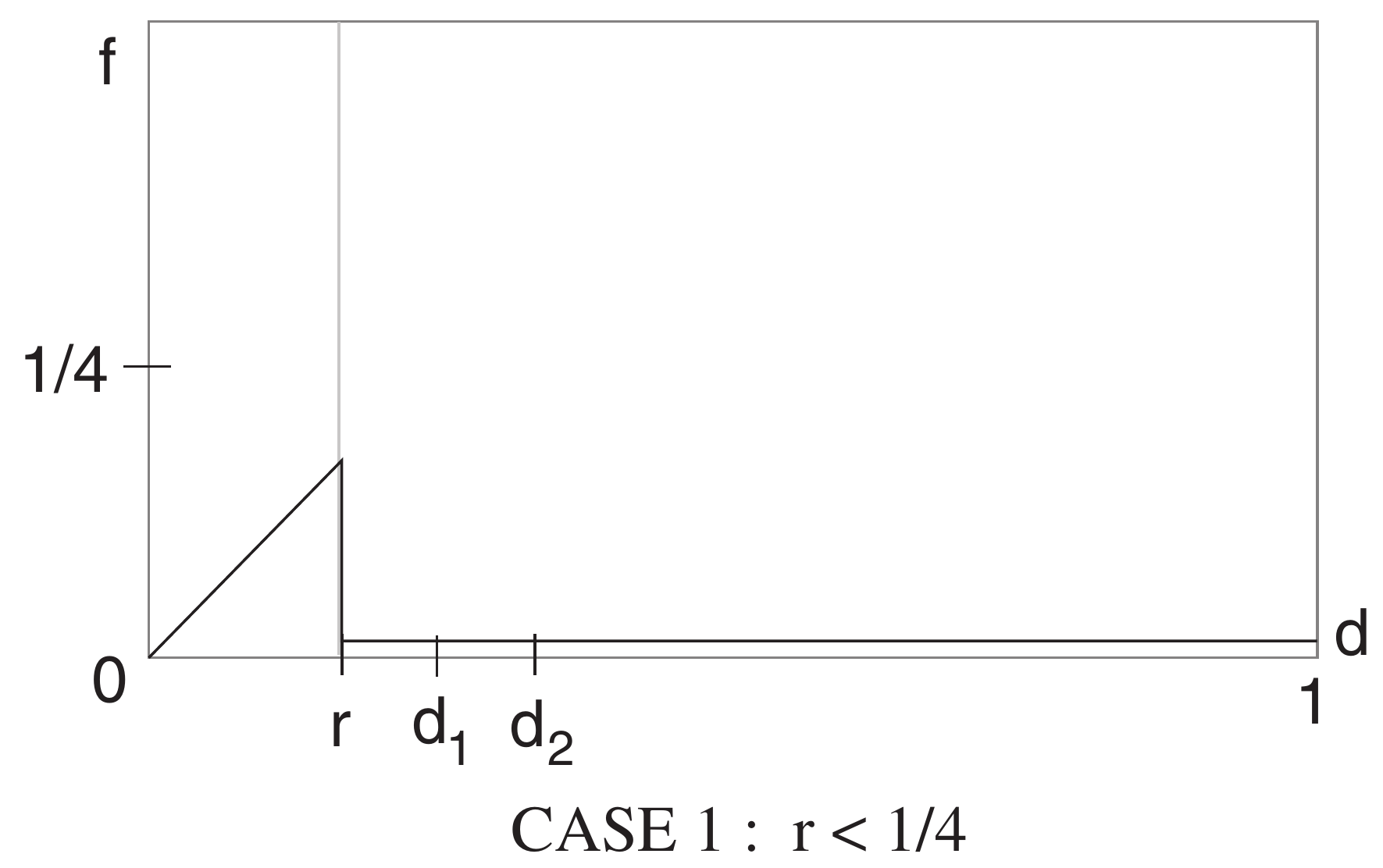}\hspace{5mm}
    \includegraphics[width=4cm,height=3.5cm]{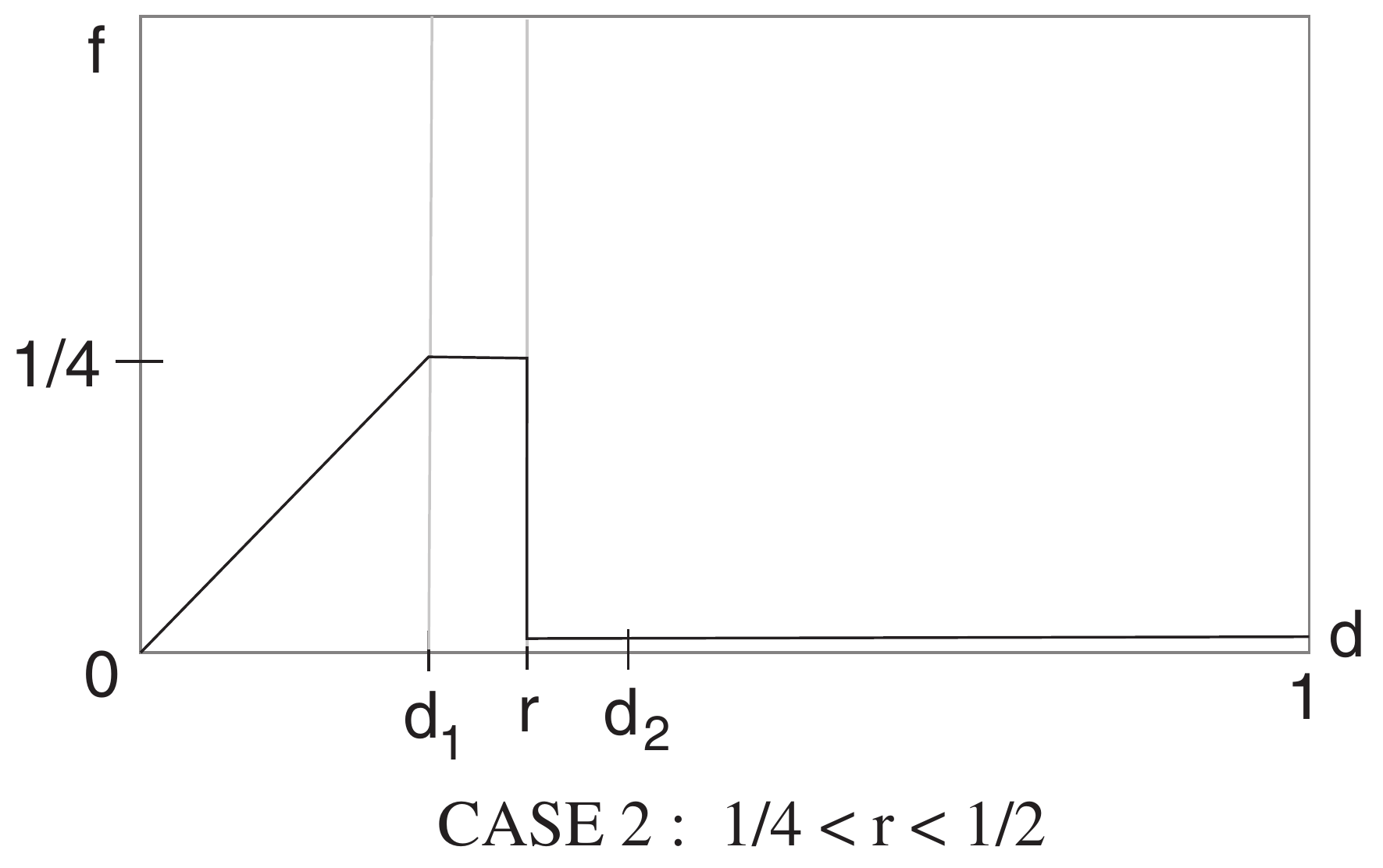}\hspace{5mm}
    \includegraphics[width=4cm,height=3.5cm]{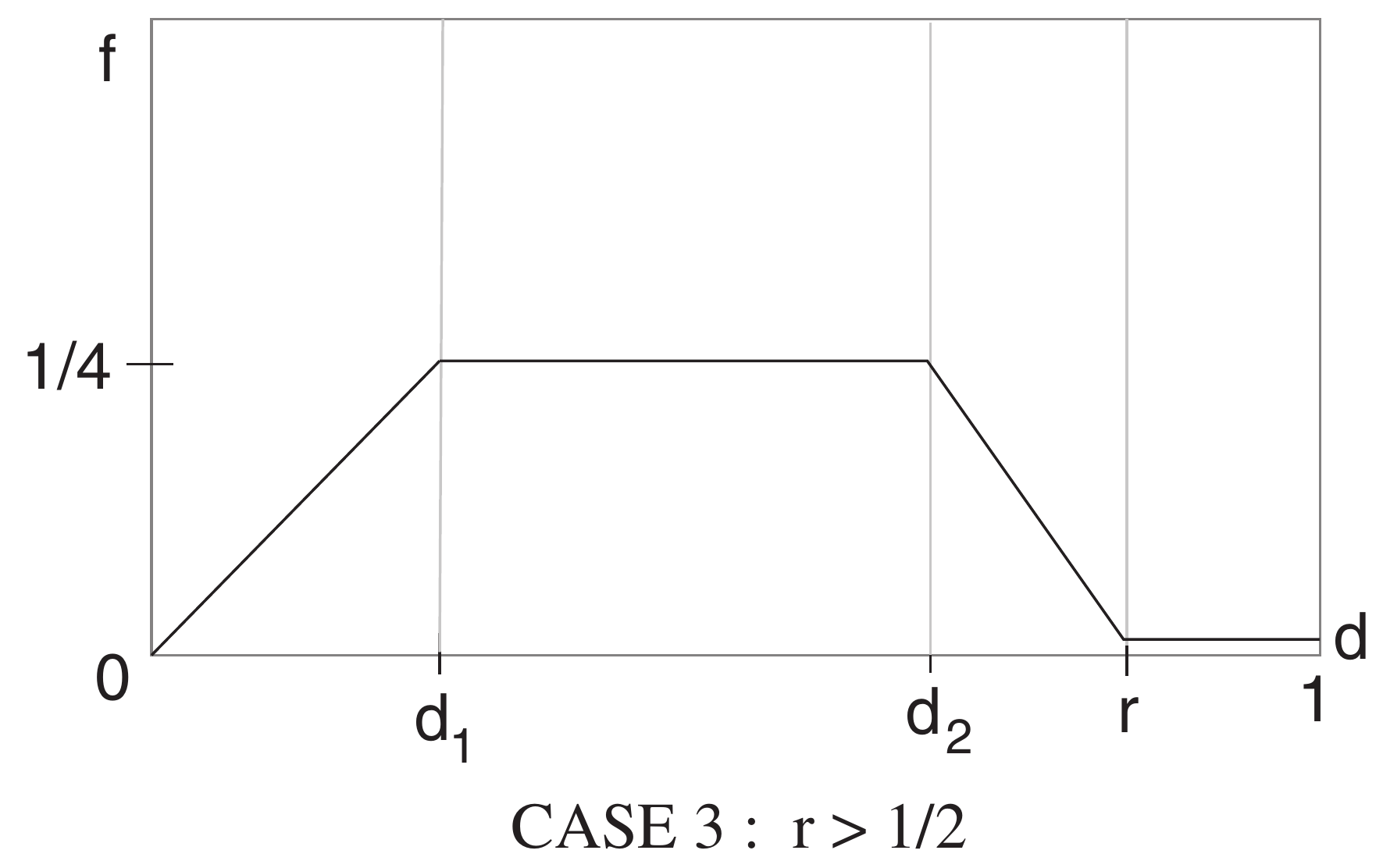}
    \caption{Summary of the fundamental diagrams obtained numerically.}
    \label{diag2-num1}
  \end{center}
\end{figure}

%------------------------------------------------------
\subsection{The traffic phases of the global fundamental diagram}
%------------------------------------------------------
\label{subsec2-4}

When $r>1/2$ there are four phases. When $r<1/2$ the third phase
vanishes, and when $r<1/4$ even the second phase vanishes. Let us
discuss the physical interpretation of these four phases.

  \begin{enumerate}
    \item \emph{Free phase:} $0\leq d\leq d_1$.

      In this case, at the stationary regime, the cars are separated by free cells and
      move freely. There is no effect of the intersection, and the
      behavior of the system is similar to a simple circular road.
      The flow is $f=d$.

\begin{figure}[htbp]
  \begin{center}
     \includegraphics[width=4cm]{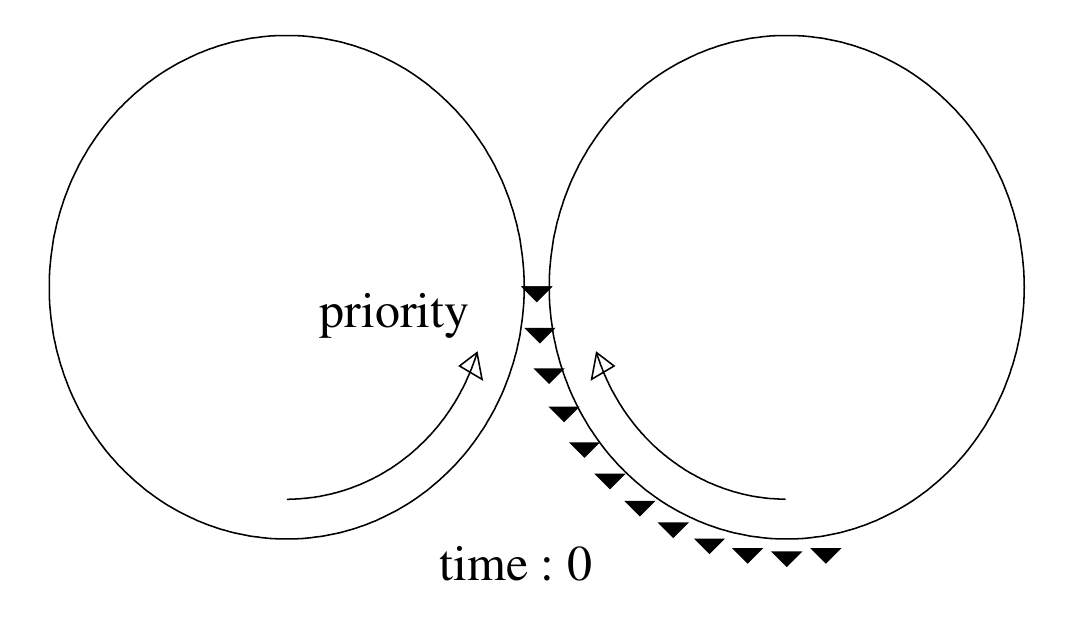}\hspace{0.5cm}
     \includegraphics[width=4cm]{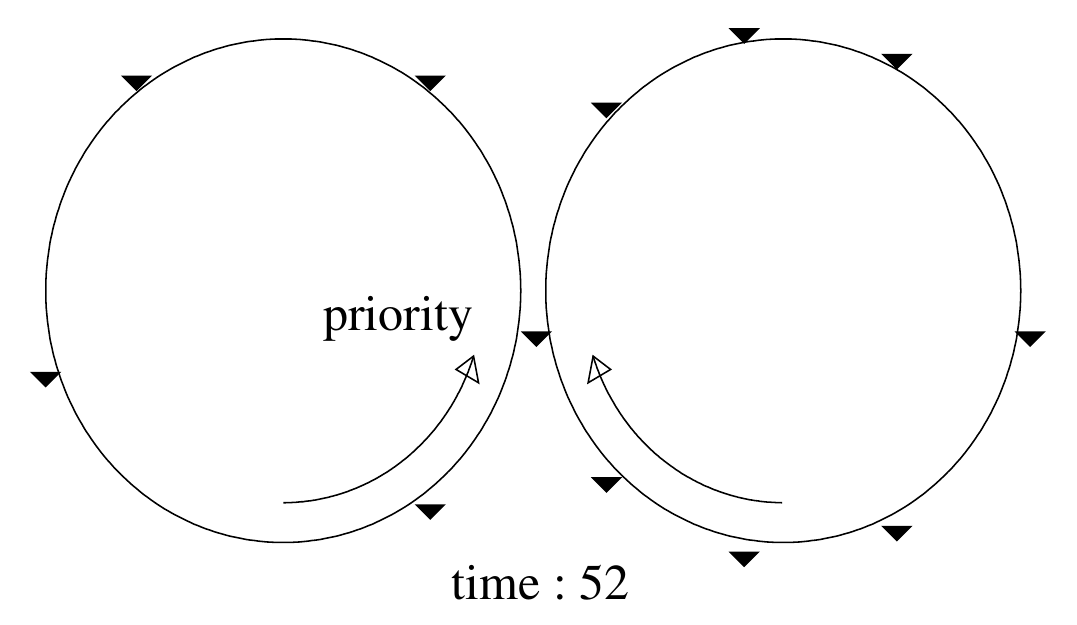}
     \caption{Positions of cars in the free phase: initial
        and periodic asymptotic. Sizes: $n=40, m=20$. The number of vehicles is $12$.}
     \label{phase1}
   \end{center}
\end{figure}     
      
    \item \emph{Saturation phase:} $d_1\leq d\leq d_2$.

      This phase appears as soon as $r>1/4$.
      In this case, the intersection is never free, and its output flow is
      maximum and equal to $1/4$. Thanks to the priority rule, after an initial transition period,
      the flow on the priority road is regular, and the density of cars
      on this road is given by the maximal output flow of the intersection
      (which is $1/4$).
      All the other cars are on the non-priority road. $N$ being the
      total number of cars, $m/4$ cars stay on the priority road and
      $N-m/4$ stay on the non-priority road. The flow on the non-priority
      road must be equal to $1/4$ in such a way that the junction be always
      fed. The non-priority road is seen as a circular road with a retarder (a
      slow cell), see \citep{CDC05}. The role of the retarder is played here by
      the intersection. Thanks to the works presented in \citep{CDC05}, it is understood that the
      flow $1/4$ is reached when the density belongs to $[1/4,3/4]$ on circular roads with retarder.
      This last constraint gives $n/4\leq N-m/4\leq 3n/4$, which is obtained when $d_1\leq d\leq
      d_2$.

\begin{figure}[htbp]
  \begin{center}
    \includegraphics[width=4cm]{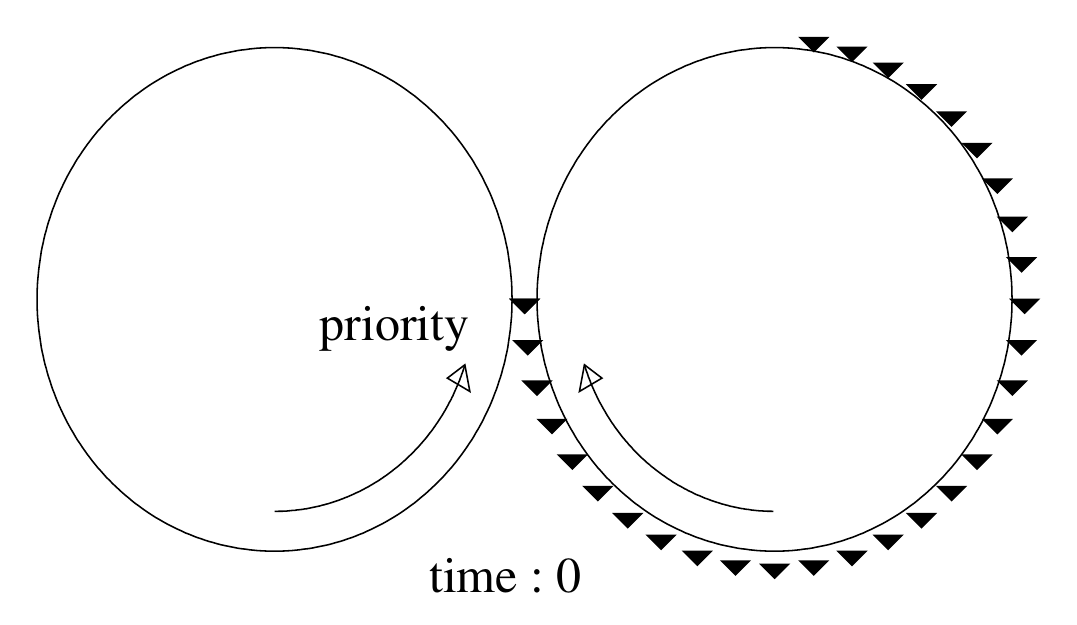}\hspace{0.5cm}
    \includegraphics[width=4cm]{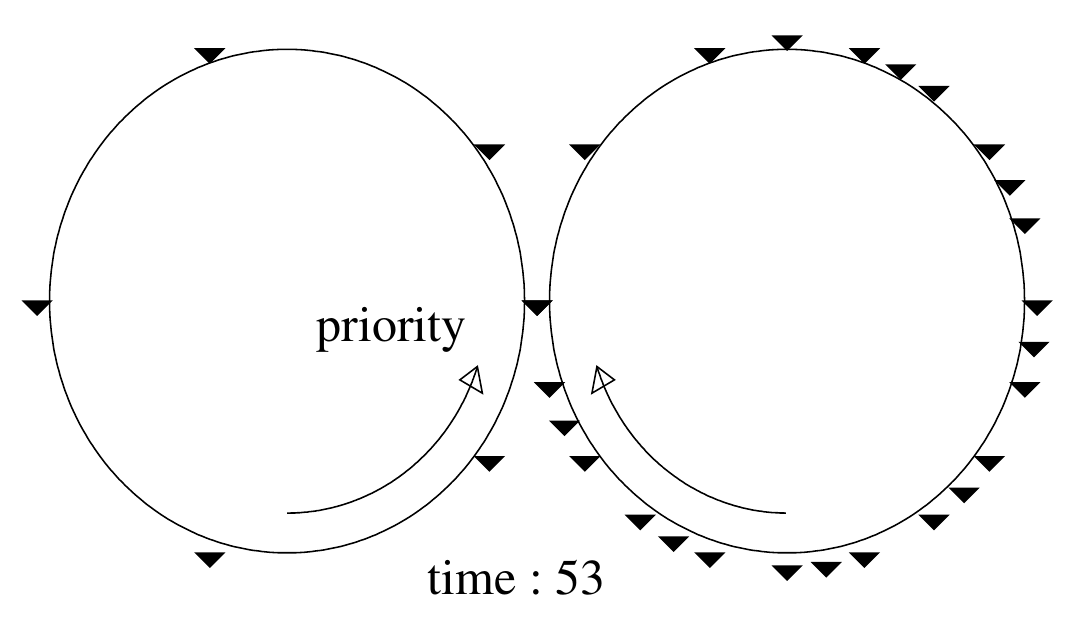}
    \caption{Positions of cars in the saturation phase: initial
       and periodic asymptotic. Sizes: $n=40, m=20$. The number of vehicles is $30$.}
    \label{phase2-1}
  \end{center}
\end{figure}     

    \item \emph{The recession phase:} $d_2\leq d\leq r$.

      This phase appears only when the non-priority road is longer
      than the priority road, that is when $r>1/2$. In this case,
      this phase appears when the global density exceeds $d_2$, that is when the
      number of vehicles in the system exceeds $m/4+3n/4$. The
      flow on the non-priority road is less than $1/4$ because the
      density of vehicles on this road exceeds $3n/4$. The
      cells immediately downstream from the junction in the non-priority road are
      crowded, and thus the junction cannot serve with its maximal
      flow, so the global flow is less than $1/4$. Since the
      density of vehicles in the priority road is given by the
      output flow of the junction, this density is less than $1/4$
      and it decreases by increasing the global density. So when
      we increase the global density slowly in this phase, the
      whole number of added vehicles is added to the non-priority
      road, and moreover, vehicles on the priority road are pumped
      onto the non-priority road.

      This phase starts when the density is such that $m/4$ vehicles are in the
      priority road, and $3n/4$ vehicles
      are on the non-priority road. Thus, we
      have $N=m/4+3n/4$, which gives $d=d_2$.

      Since this phase is characterized by the pumping of the
      priority road vehicles, the end of this phase appears when the density
      is such that no vehicle stays on the priority road.
      When  we have $n$ vehicles on the non priority road
      its density is $1$. Thus we get $N=n$, which gives $d=r$,
      thus $d_2<d<r.$

\begin{figure}[htbp]
  \begin{center}
    \includegraphics[width=4cm]{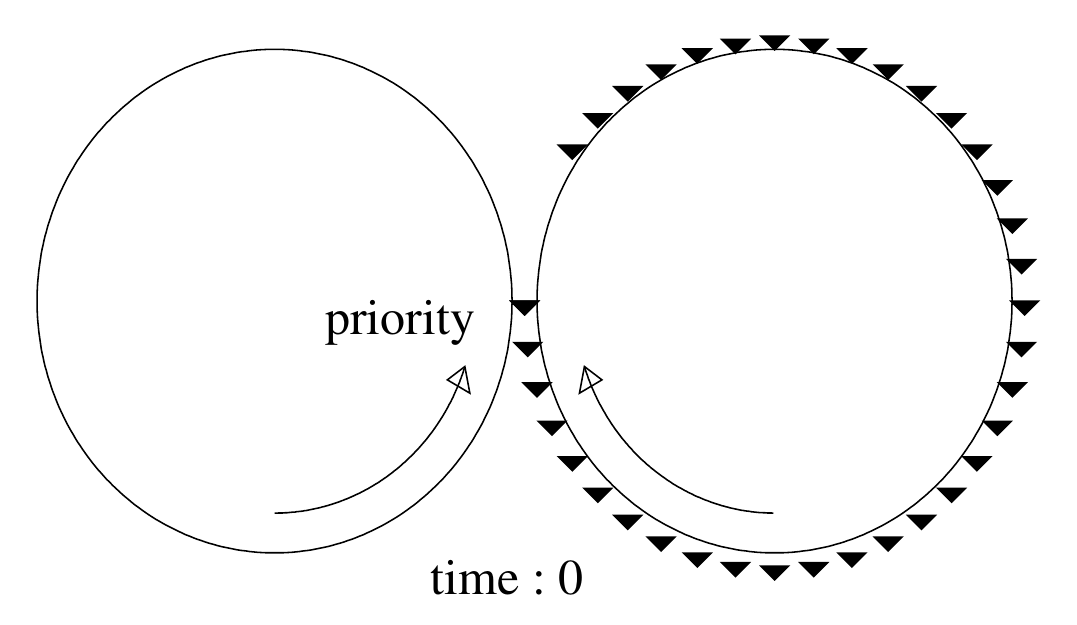}\hspace{0.5cm}
    \includegraphics[width=4cm]{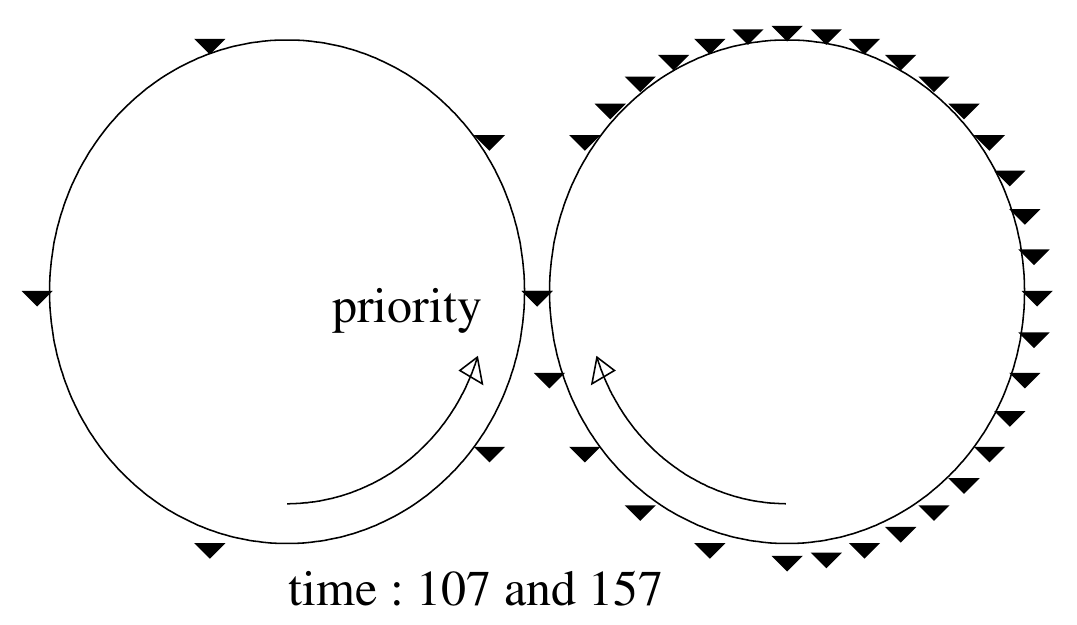}\hspace{0.5cm}
    \includegraphics[width=4cm]{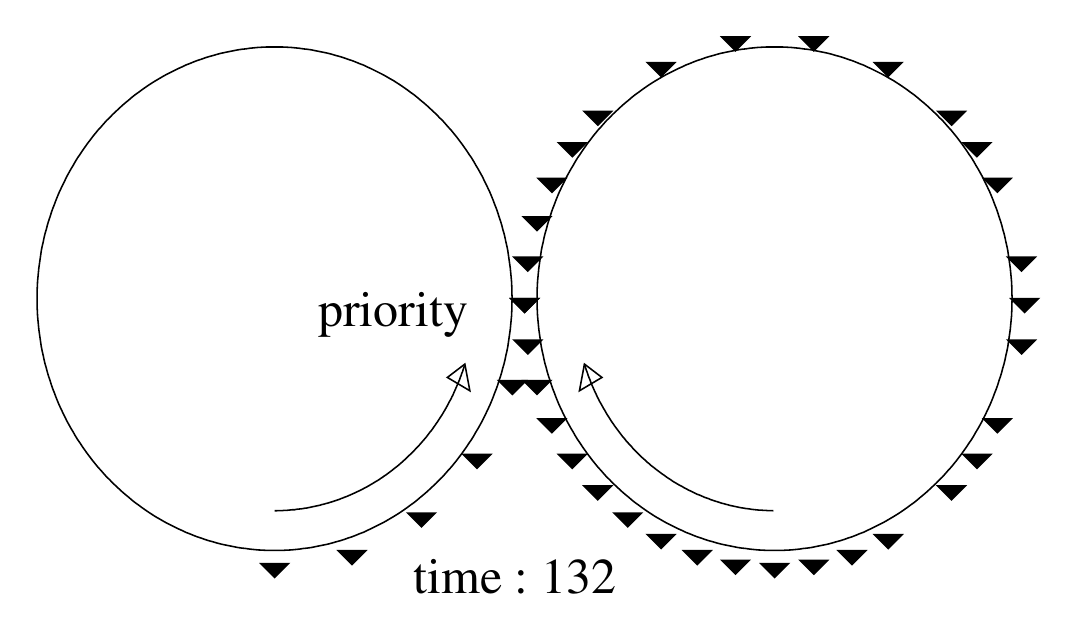}\hspace{0.5cm}
    \caption{Positions of cars in the recession phase: initial
       and periodic asymptotic. Sizes: $n=40, m=20$. The number of
        vehicles is $37$.}
    \label{phase2-2}
  \end{center}
\end{figure}    

    \item \emph{The freeze phase:} $r\leq d\leq 1$.

      When the total number of cars exceeds the size of the non-
      priority road, a jam appears. In this case, the non-priority road fills
      up. When a car in the intersection wants to enter this
      road, the system blocks.

\begin{figure}[htbp]
  \begin{center}
    \includegraphics[width=4cm]{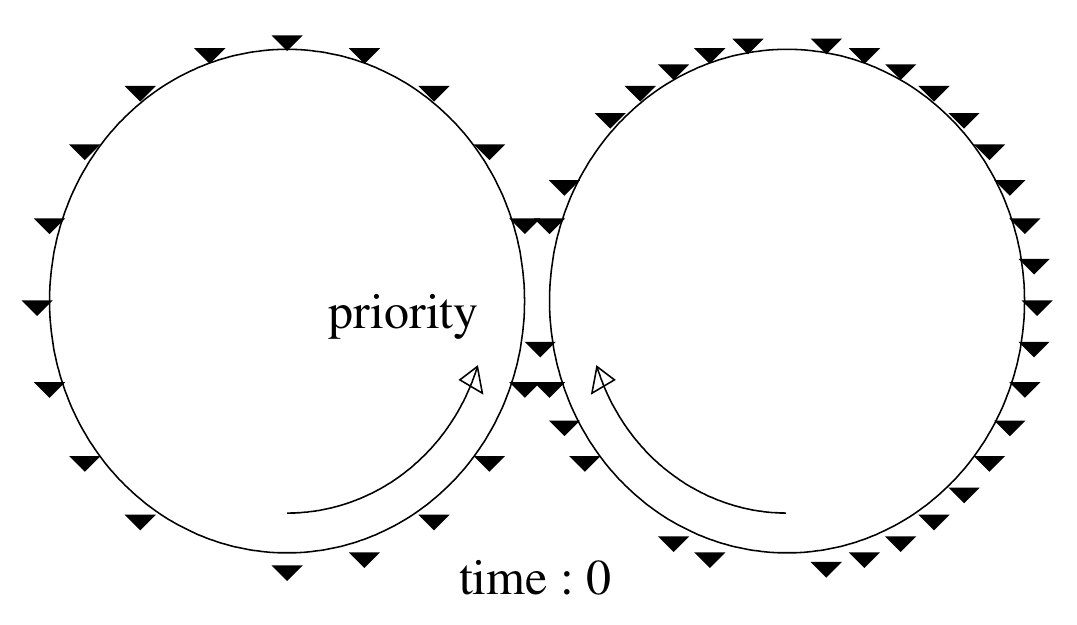}\hspace{0.5cm}
    \includegraphics[width=4cm]{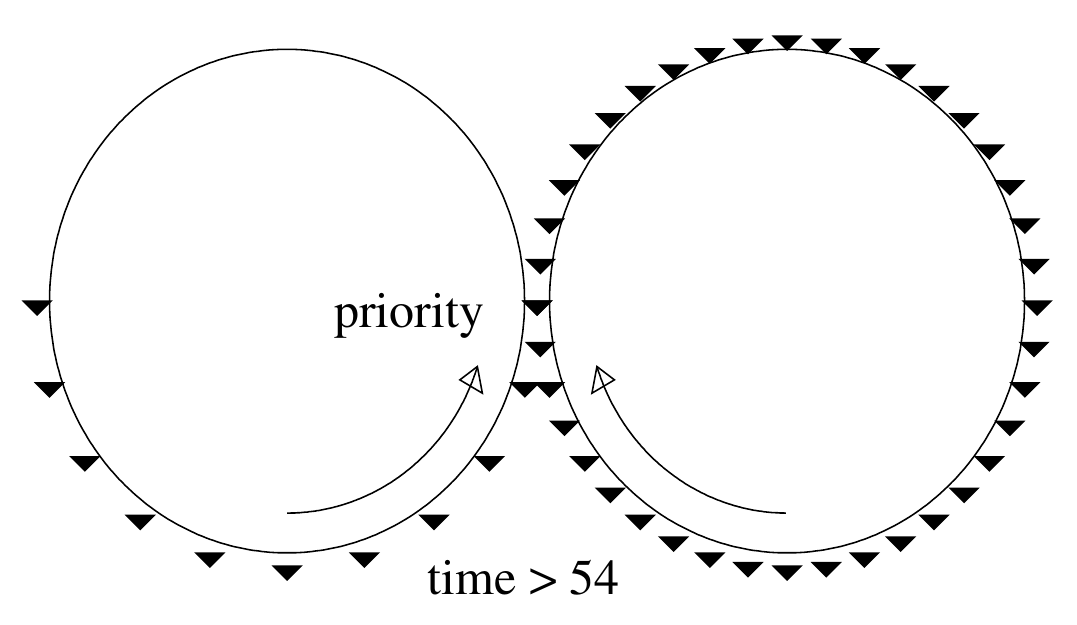}\hspace{0.5cm}
    \caption{Positions of cars in the freeze phase: initial
       and periodic asymptotic. Sizes: $n=40, m=20$. The number of vehicles is
       $50$.}
    \label{phase3}
  \end{center}
\end{figure}     
      
  \end{enumerate}

%----------------------------------------------------------------
\subsection{The fundamental diagram of the individual roads}
%----------------------------------------------------------------

\begin{figure}[htbp]
  \begin{center}
    \includegraphics[width=6cm]{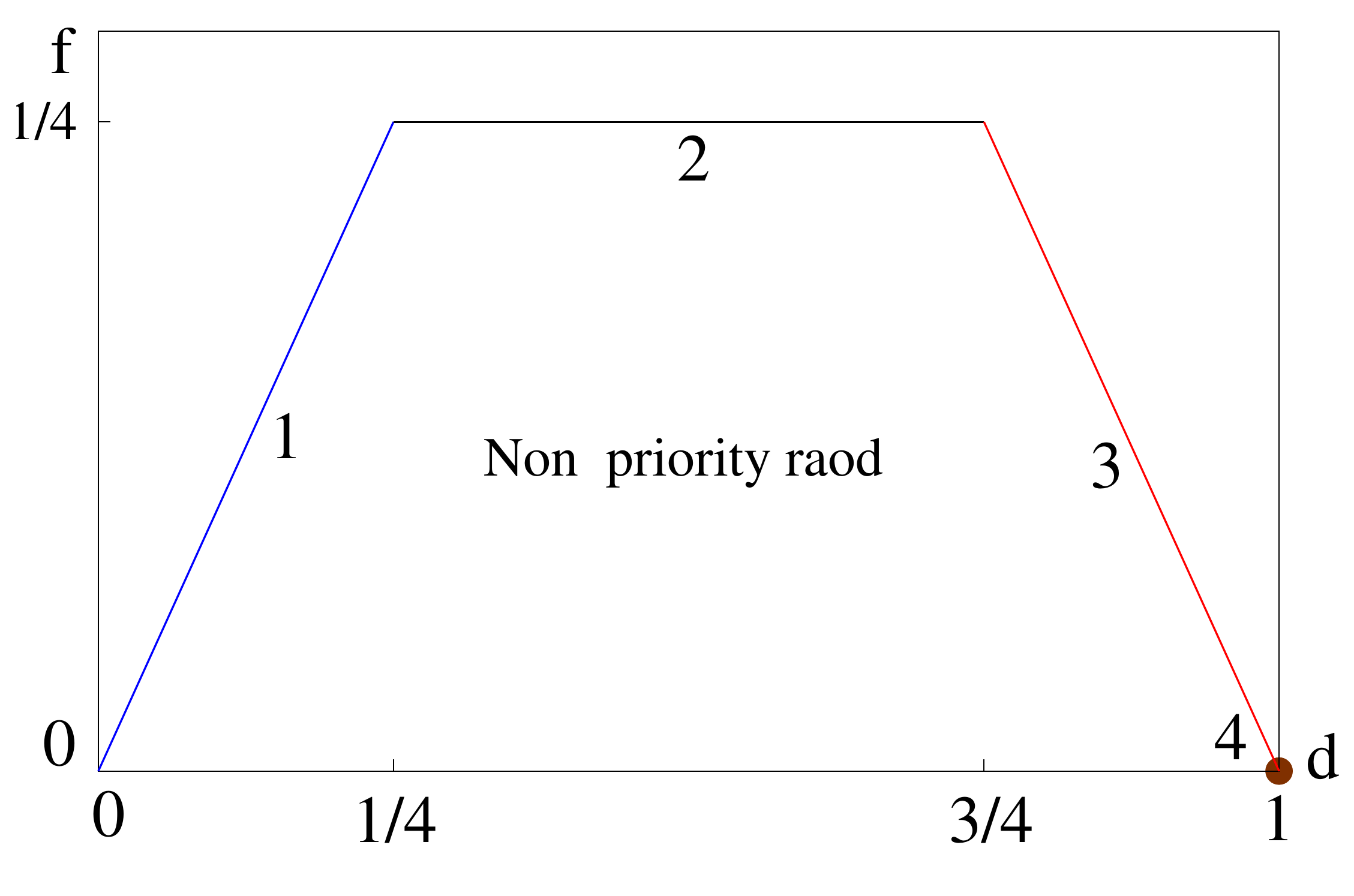}\hspace{1cm}
    \includegraphics[width=6cm]{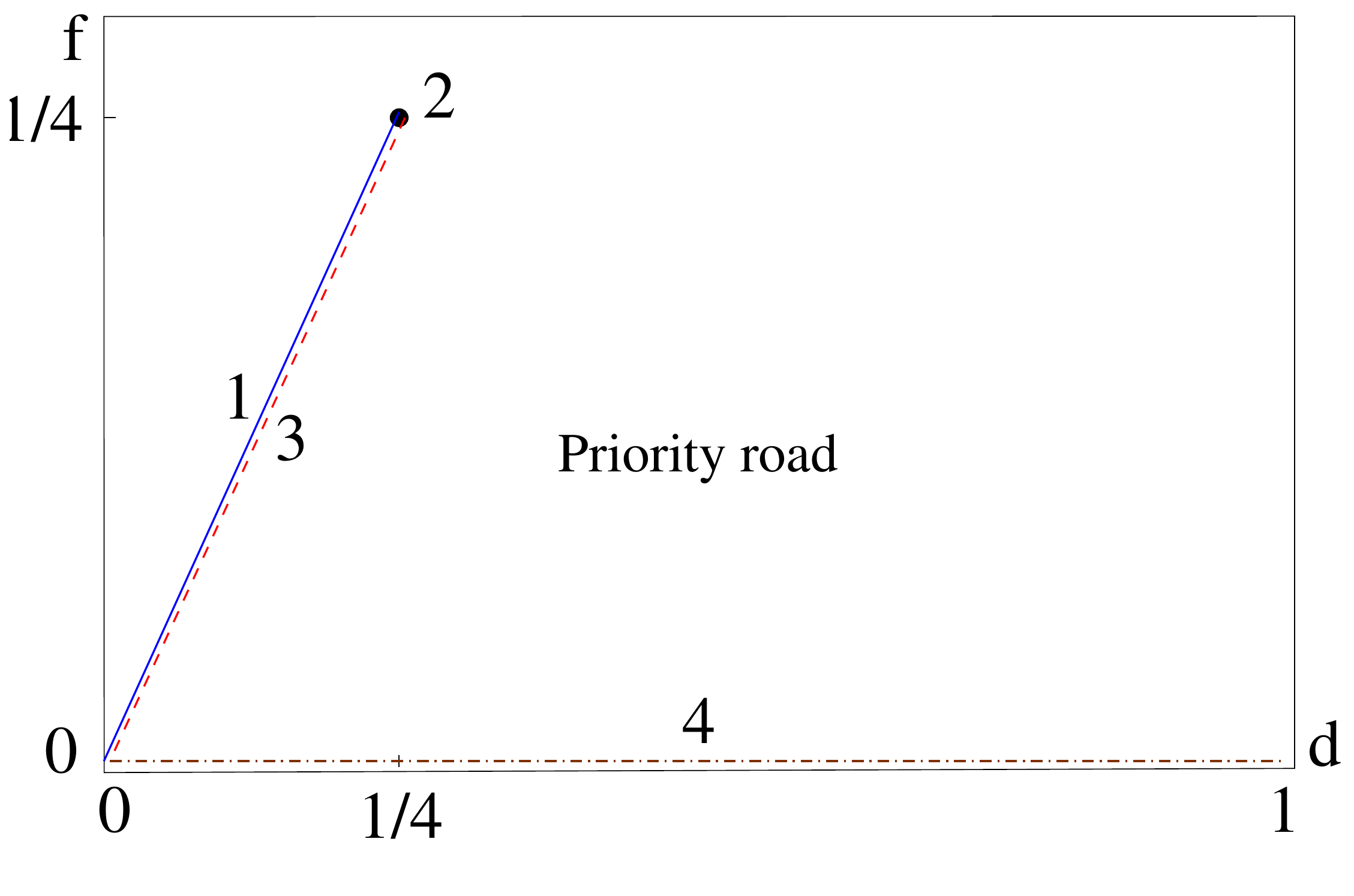}
    \caption{Approximate fundamental diagram on each road. 1: Free phase, 2: Saturation phase,
         3: Recession phase, 4: Freeze phase.}
    \label{diag-prior}
  \end{center}
\end{figure}

Let us discuss the fundamental diagrams for the two roads. For this
we need the densities $d_m$ on the priority road and $d_n$ on the
non-priority road. It is easy to check that:
\begin{equation}\label{relation}
  d=r d_n+(1-r) d_m.
\end{equation}
Using (\ref{relation}) we get the diagrams plotted on Figure~\ref{diag-prior}.

\begin{enumerate}
  \item \emph{Free phase:}
    In this phase, we have $d=d_m=d_n=f$.  Thus, the diagrams are $f=d_m$ and
    $f=d_n$ with $0\leq d_m\leq 1/4$ and $0\leq d_n\leq 1/4$.
  \item \emph{Saturation phase:}
    In this case, the global flow $f$ and the density $d_m$ are constant
    and are equal to $1/4$. Then the diagram of the priority road is restricted to
    the point $(1/4,1/4)$.

    From (\ref{relation}), we get
    $d_n=d/r-(1-r)/4r,$
    then $1/4 \leq d_n\leq 3/4$ from
    $d_1\leq d\leq d_2$.
    Hence the diagram of the non-priority road is given by the segment
    $[(1/4,1/4)\;,\;(3/4,1/4)]$.
  \item \emph{Recession phase:}
    Let us give a simplified description of what happens in this phase.
    In fact, we see on numerical simulations a more complicated periodic regime
    where the densities of the car do not stay approximately constant
    as in the other phases.

    The pumping phenomena makes the density of the
    priority road low. Thus the vehicles are moving freely on this
    road. We have $d_m=f=(r-d)/(2r-1)$. Then $d_2\leq d\leq r$ gives $0\leq d_m\leq
    1/4$. Therefore, the diagram for the priority road is the segment
    $[(1/4,1/4)\;,\;(0,0)]$.

    From (\ref{relation}), we get
    $d_n=d/(2r-1)d-(1-r)/(2r-1),$
    which must satisfy $3/4<d_n<1$ to be crowded enough. That is implied by $d_2<d<r$.
    Moreover, we have
    $f=1-d_n=(r-d)/(2r-1)$ for a road without intersection (see \citep{CDC05}).
    Therefore, the diagram for the non-priority road is the
    segment $[(3/4,1/4)\;,\;(1,0)]$.
  \item \emph{Freeze phase:}
    In this case, we have $d_n=1$ and $f=0$, and the diagram of the
    non-priority road is restricted to the point $(1,0)$.
    From (\ref{relation}), we get
    $d_m=d/(1-r)-r/(1-r)$.
    Therefore, $r\leq d\leq 1$ gives $0\leq d_m\leq 1$.
    The diagram for the priority road in this phase is the segment
    $[(0,0)\;,\;(1,0)]$.
\end{enumerate}

In this study, the turning movement percentage is taken
equal to $1/2$, but the same approach can be used to study the
general case, where this percentage is a parameter $\alpha$.
Preliminary works, show that a function of $r$ and $\alpha$ could
play the same role as $r$ here. 

%-----------------------------------------------
\section{Extension to more than one junction}
%-----------------------------------------------

In this section, we extend the 
one junction model to the case
with two junctions first, then to a regular city on a torus.
We use the same vehicle dynamics and junction management
(right priority, turning policy) as in the case of one junction.
To be able to fix
the car densities, we consider, as in section~\ref{sec2}, closed
systems. In the case of a city, the simplest way to obtain a
closed system is to consider a regular city on a torus; see
Figure~\ref{regular}. We show that the shape of the fundamental
diagram is mainly the same as in the case of one junction.

%---------------------------------------
\subsection{The case of two junctions}
%---------------------------------------

\begin{figure}[htbp]
  \begin{center}
    \includegraphics[width=5cm]{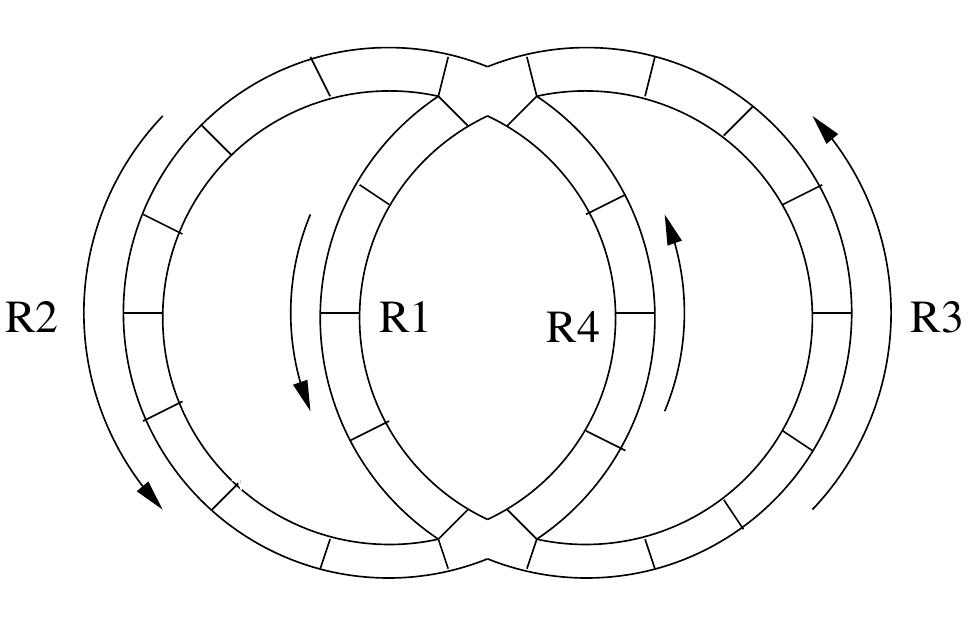}
    \caption{Two circular roads crossing on two junctions.}
    \label{fouroads}
  \end{center}
\end{figure}

We consider a system of two circular roads crossing on two
junctions, as shown on Figure~\ref{fouroads}. This configuration
gives four roads without intersections (R1, R2, R3 and R4). The road
R2 (resp. R3) has priority over the road R1 (resp. R4).

Using the same approach as in the case of one junction, we can
derive easily the dynamics of this system. A systematic way to
determine such a kind of model is described in \citep{PHD} . Here
we will discuss only the fundamental diagram obtained and its
phases.

\begin{figure}[htbp]
  \begin{center}
    \includegraphics[width=5cm]{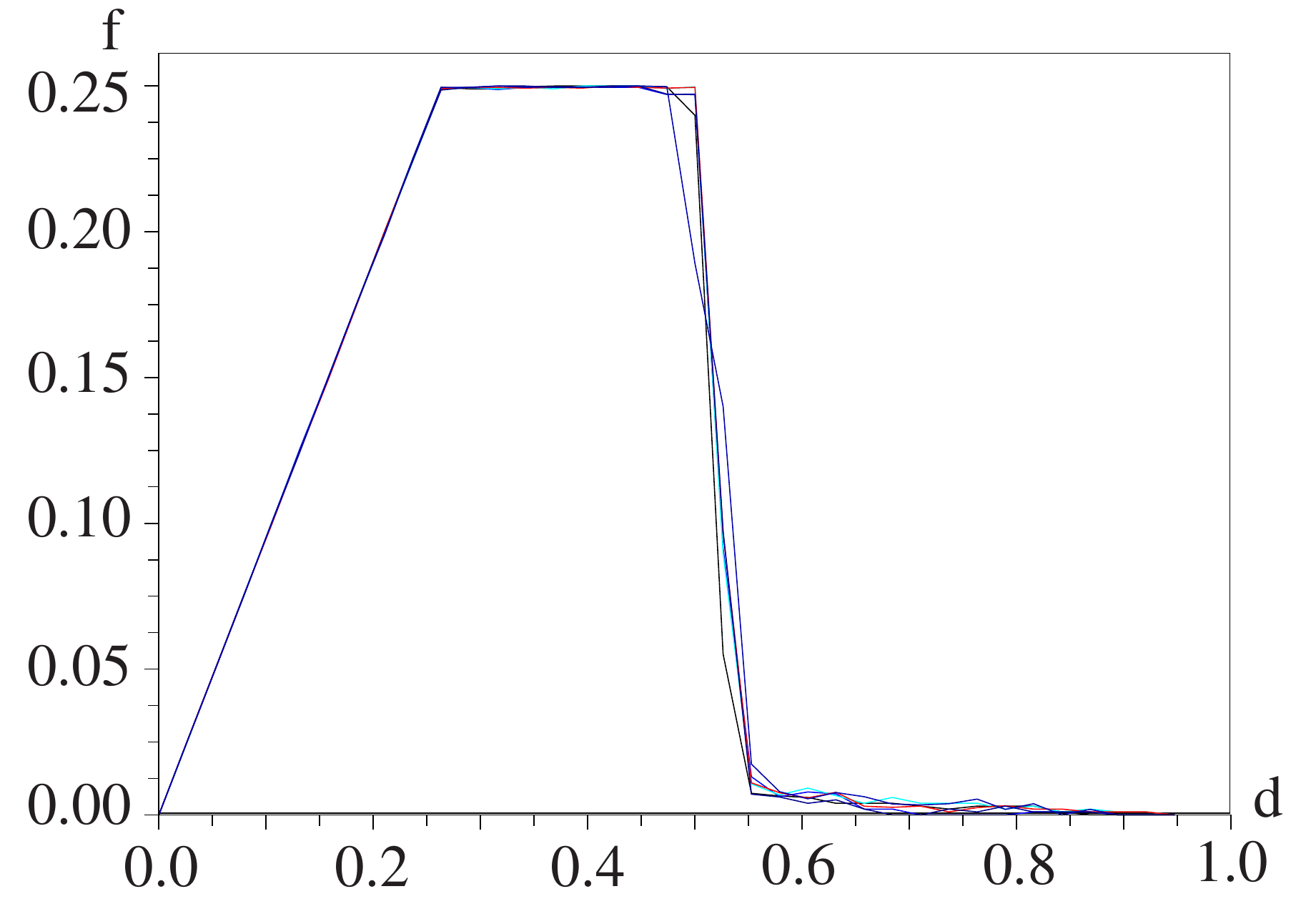}\hspace{15mm}
    \includegraphics[width=5cm]{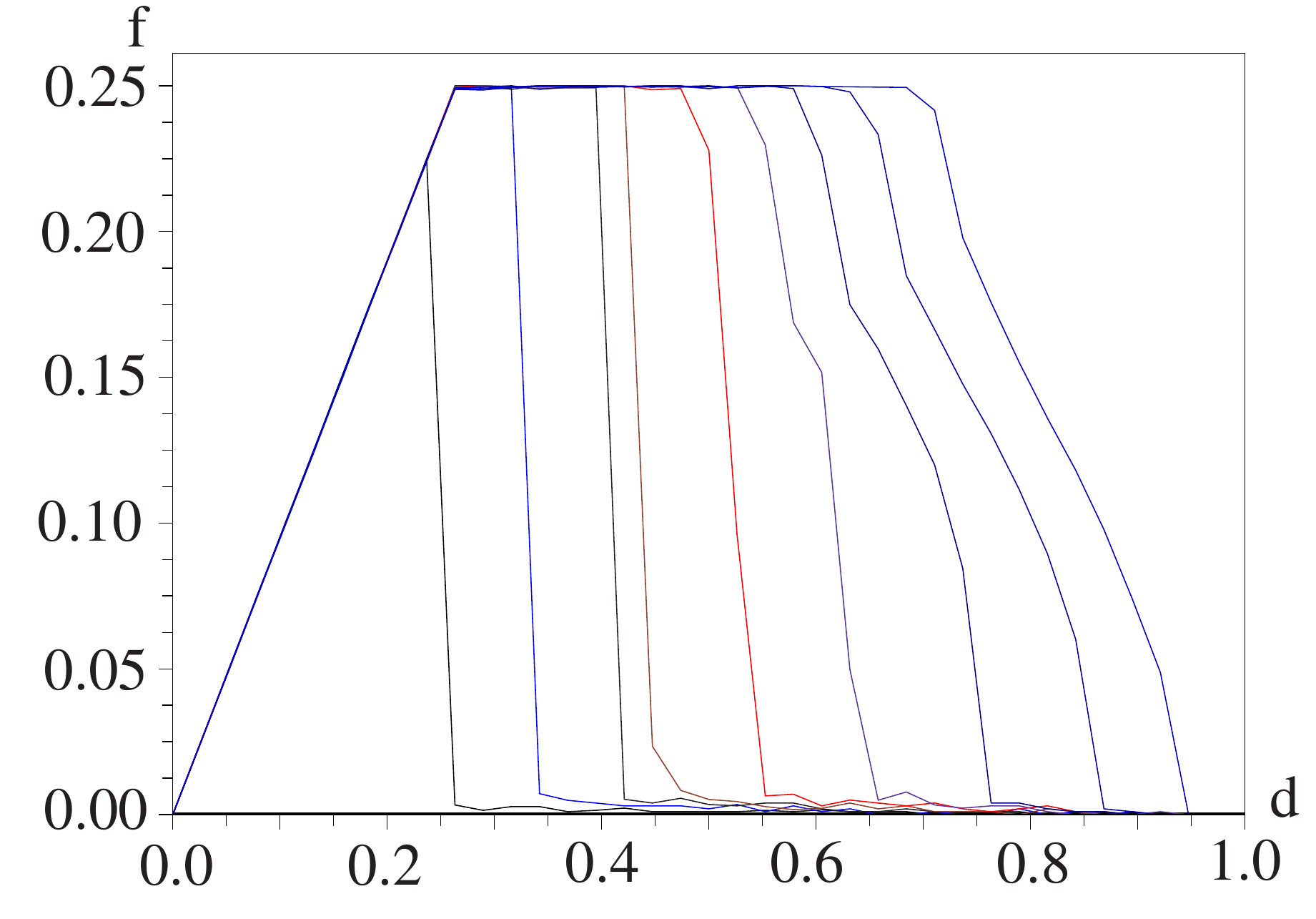}
    \caption{On the left side: the fundamental diagram on the whole system for various
     sizes of the four roads with $r$ maintained constant: $r=1/2$. On the right side:
     the fundamental diagram on the whole system for different values of $r$,
     respectively from the left to the right: $2/10, 3/10, 4/10, 5/10, 6/10, 7/10, 8/10,
     9/10, 9.5/10$.}
    \label{dep}
  \end{center}
\end{figure}

$\bullet$ The first important observation that we can make is that
the fundamental diagram depends on the sizes of the four roads
only through the ratio between the sum of the priority road sizes
and the sum of the non-priority road sizes called $r$; see
Figure~\ref{dep}.

\begin{figure}[htbp]
  \begin{center}
    \includegraphics[width=5cm]{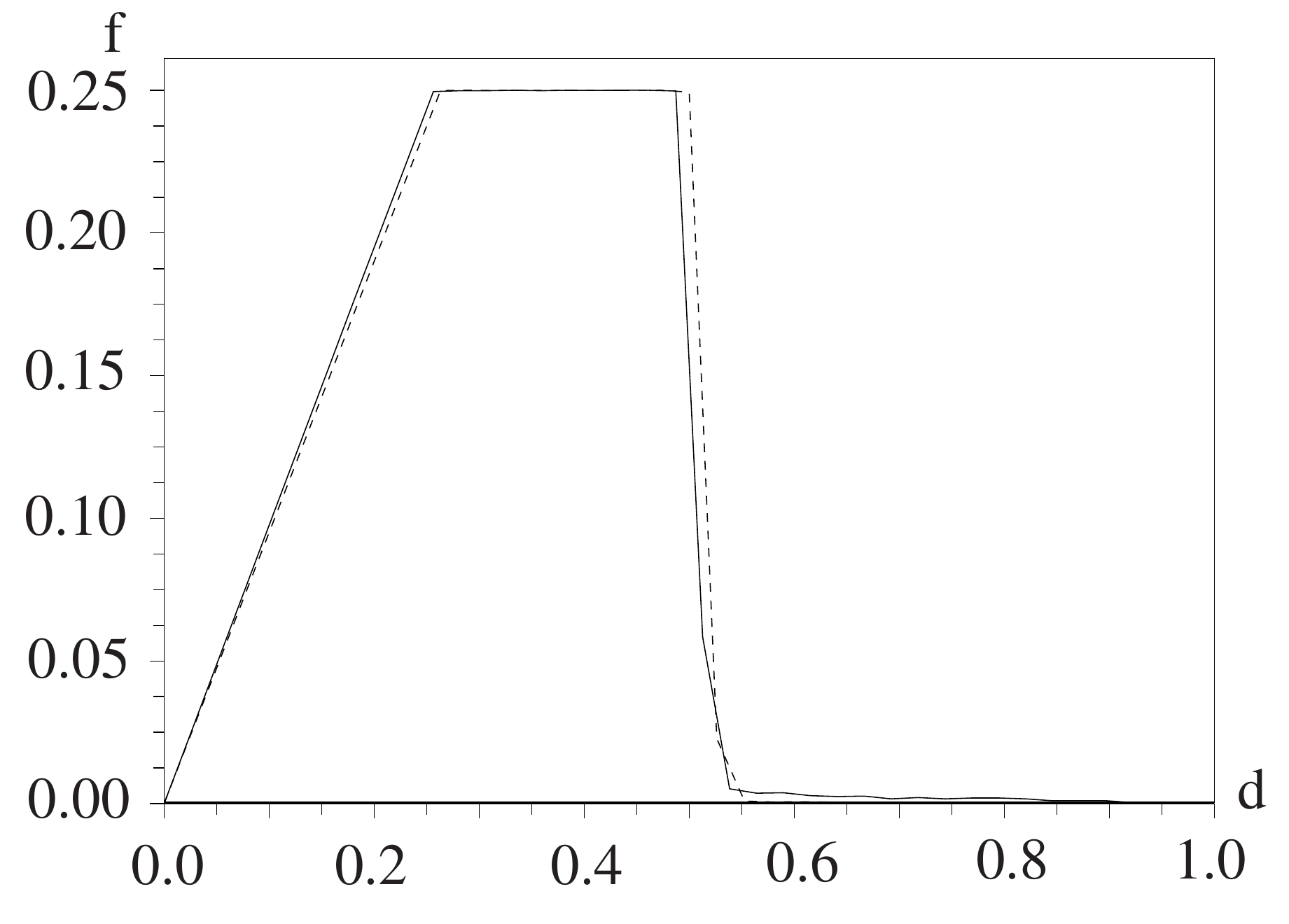}
    \caption{Comparison of the fundamental diagrams for one and two junctions with the same $r=1/2$.}
    \label{un=deux}
  \end{center}
\end{figure}

$\bullet$ The second important observation is that in terms of the
fundamental diagram of traffic, a system of four roads with two
junctions with a ratio $r$ between the sum of the priority road
sizes and the sum of the non-priority road sizes behaves like a
system of two roads with one junction with the same ratio $r$
between the size of the priority road and the size of the
non-priority road; see Figure~\ref{un=deux}.

\begin{figure}[htbp]
  \begin{center}
    \includegraphics[width=3cm]{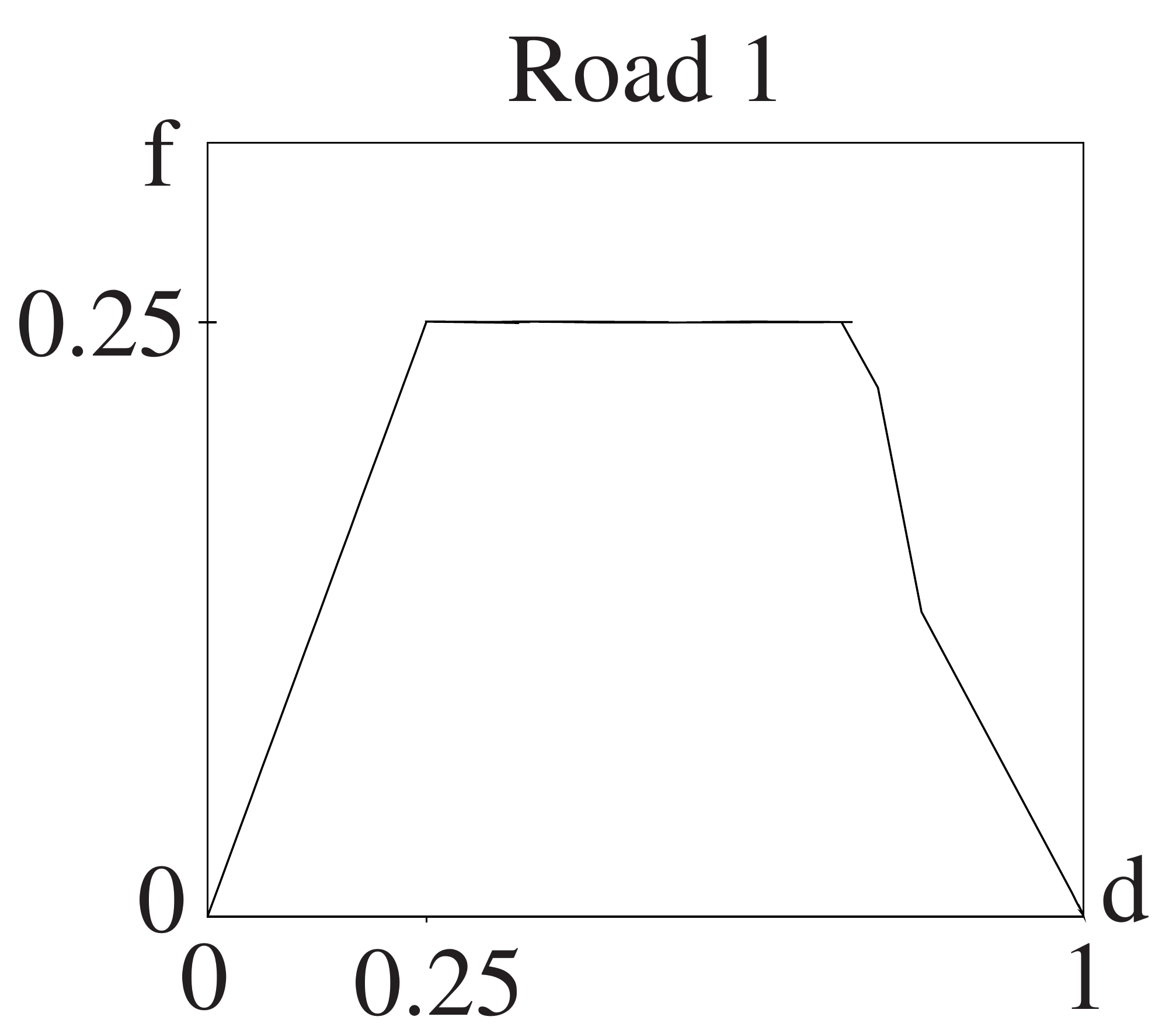} \hspace{2mm}
    \includegraphics[width=3cm]{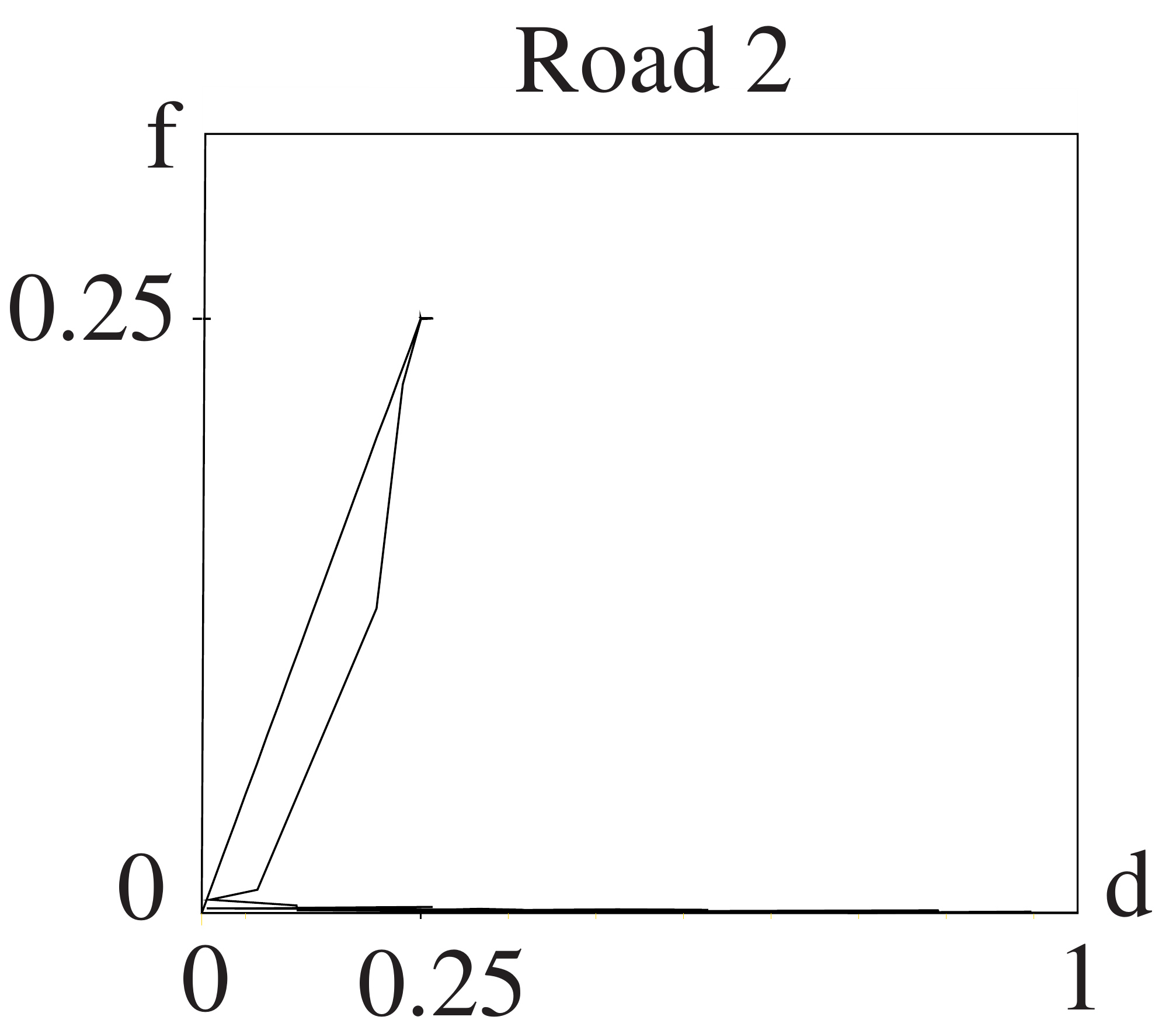} \hspace{2mm}
    \includegraphics[width=3cm]{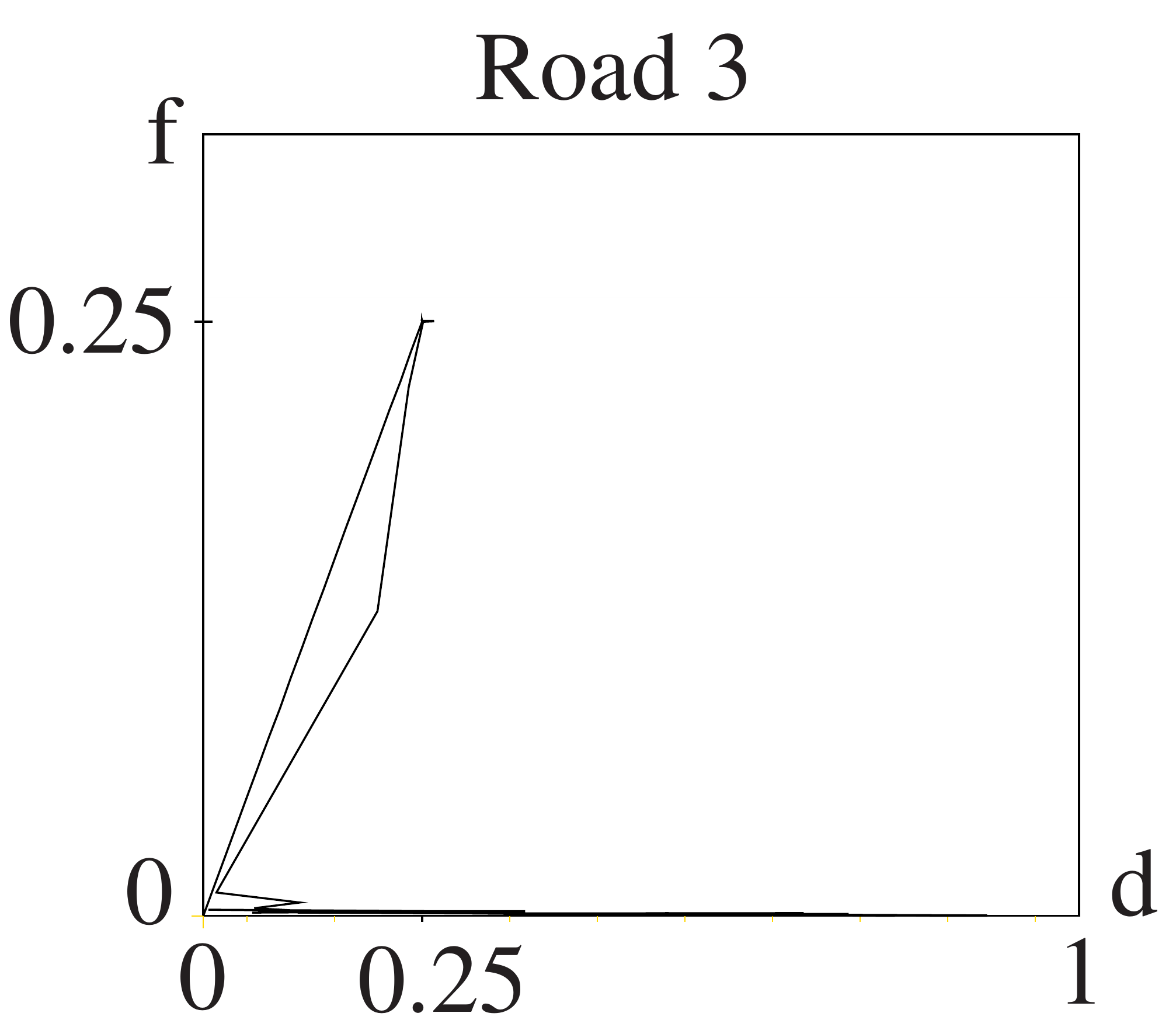} \hspace{2mm}
    \includegraphics[width=3cm]{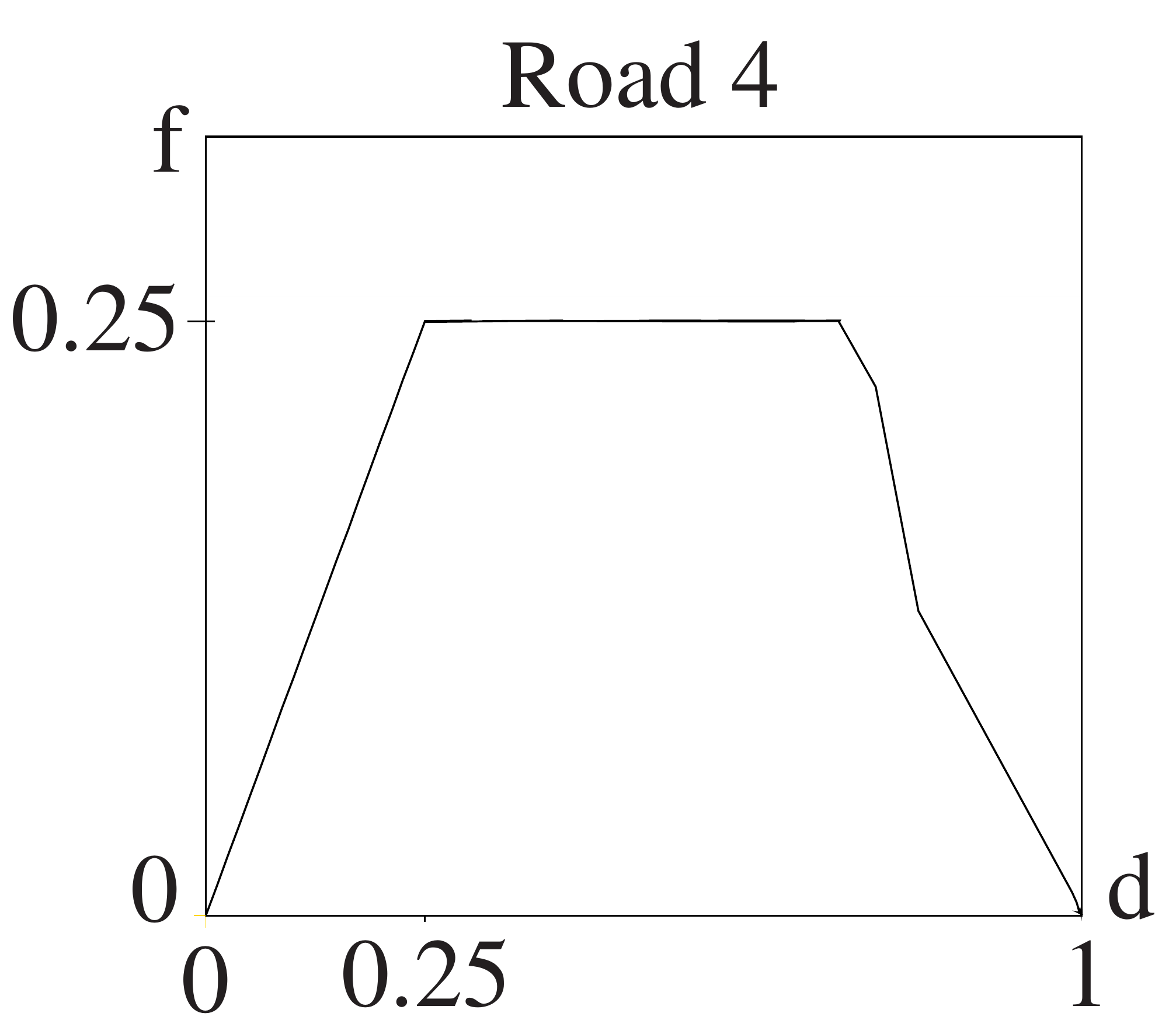}
    \caption{Fundamental diagram on the roads 1, 2, 3 and 4.}
    \label{loi1234}
  \end{center}
\end{figure}

\begin{figure}[htbp]
  \begin{center}
    \includegraphics[width=3.3cm]{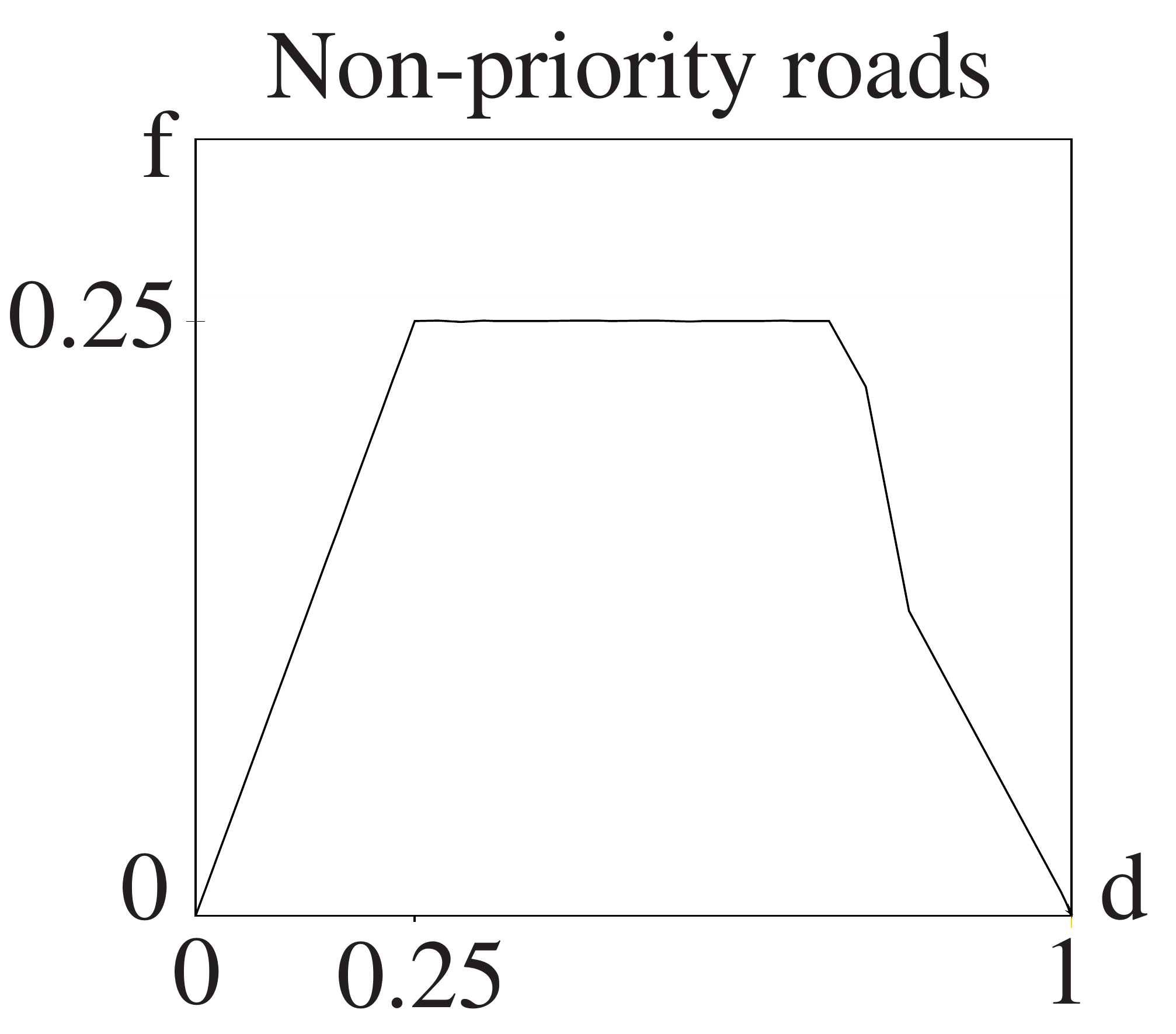} \hspace{1cm}
    \includegraphics[width=3.3cm]{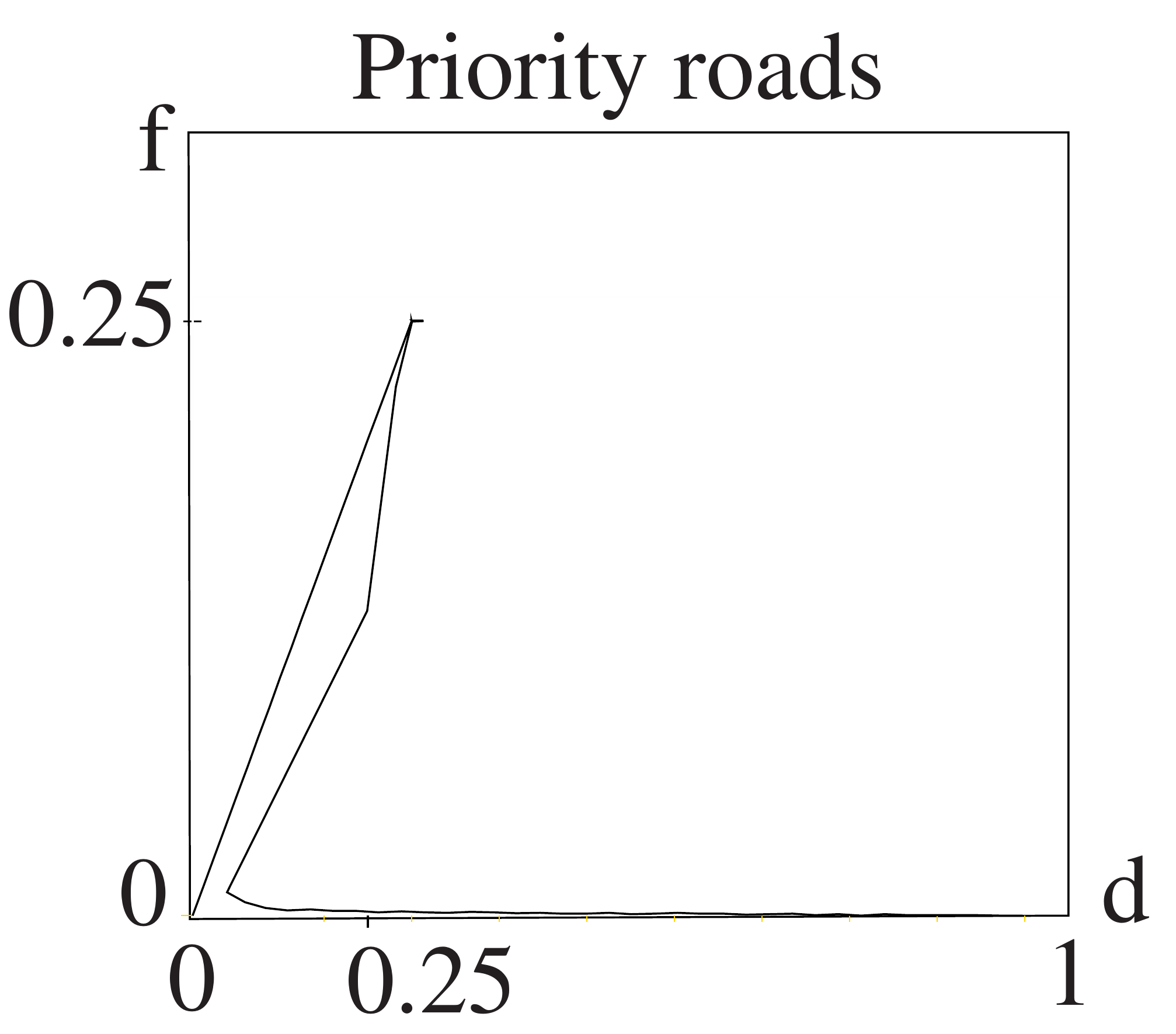}
    \caption{The fundamental diagram of the priority roads (2 and 3) and that of the non priority roads (1 and 4).}
    \label{loi14-23}
  \end{center}
\end{figure}

$\bullet$ The third observation is that the traffic phases are
similar to those obtained in the case of one junction. Here also,
the asymptotic regimes are well understood.
Figures~(\ref{ph1-deux}), (\ref{ph2-deux}), (\ref{ph3-deux}) and
(\ref{ph4-deux}) show the stationary regimes corresponding to
the four phases.
  \begin{enumerate}
    \item \emph{Free phase:} After a mixing regime where the cars split
      on the different roads in equal populations, a periodic regime appears where the
      vehicles move freely and where the priority rule is never
      applied; see Figure~(\ref{ph1-deux}). During the periodic
      regime, the vehicles are sufficiently separated so as to not disturb
      each other on the roads or in the intersections.
      Thus, all the vehicles move every time.
      That gives the same flow on all the roads, which is equal to
      the density of vehicles on the network. This phase
      corresponds to densities from $0$ to~$1/4$.

\begin{figure}[htbp]
  \begin{center}
    \includegraphics[width=4cm]{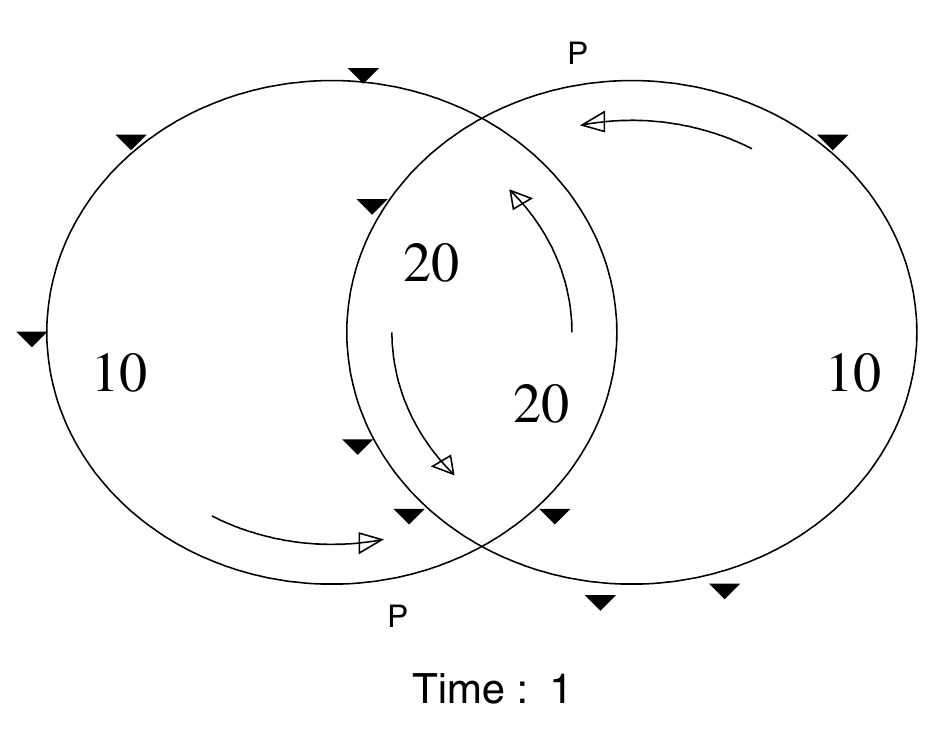} \hspace{1cm}
    \includegraphics[width=4cm]{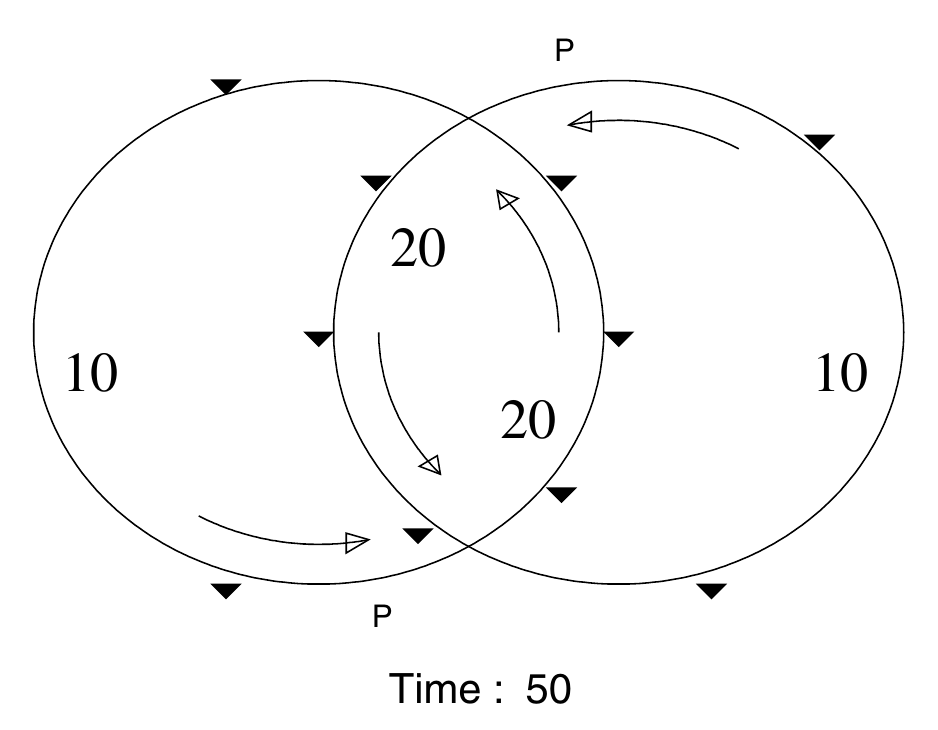}
    \caption{Initial and asymptotic periodic configurations
          during the free phase. The size of each priority road is 10.
          The size of each non-priority road is 20. The
          total number of vehicles is 10.}
    \label{ph1-deux}
  \end{center}
\end{figure}     
            
    \item \emph{Saturation phase:} A periodic regime is reached
      where the junctions serve with their maximal flow applying
      the priority rule. Thus, the vehicles are accumulated in the
      non-priority roads; see Figure~\ref{ph2-deux}. During the
      periodic regime, the average number of vehicles moving on
      the priority roads remains the same as in the case of the
      density $d=1/4$. This means that all the added vehicles,
      with respect to the density $1/4$, are absorbed by the non-priority roads.
      The junctions serve with the
      maximal flow. Hence, the average flow is equal to $1/4$. The saturation phase
      is given by a horizontal segment on the global diagram (Figure~\ref{dep}) and on the
      diagrams corresponding to the non-priority roads
      (Figure~\ref{loi1234}, Roads 1 and 4. ~Figure~\ref{loi14-23}, left side).
      The saturation phase is given by a dot on the diagrams corresponding to the priority
      roads, because the number of vehicles in the priority roads is unchanged (Figure~\ref{loi1234}, Roads 2 and 3.
      Figure~\ref{loi14-23}, right side).

\begin{figure}[htbp]
  \begin{center}
    \includegraphics[width=4cm]{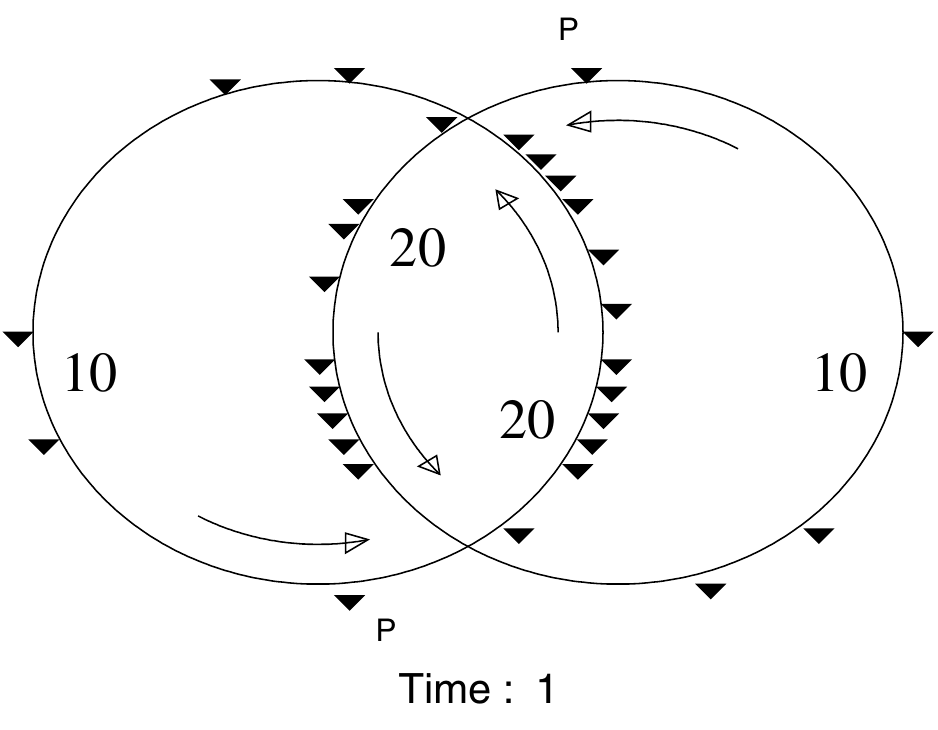} \hspace{1cm}
    \includegraphics[width=4cm]{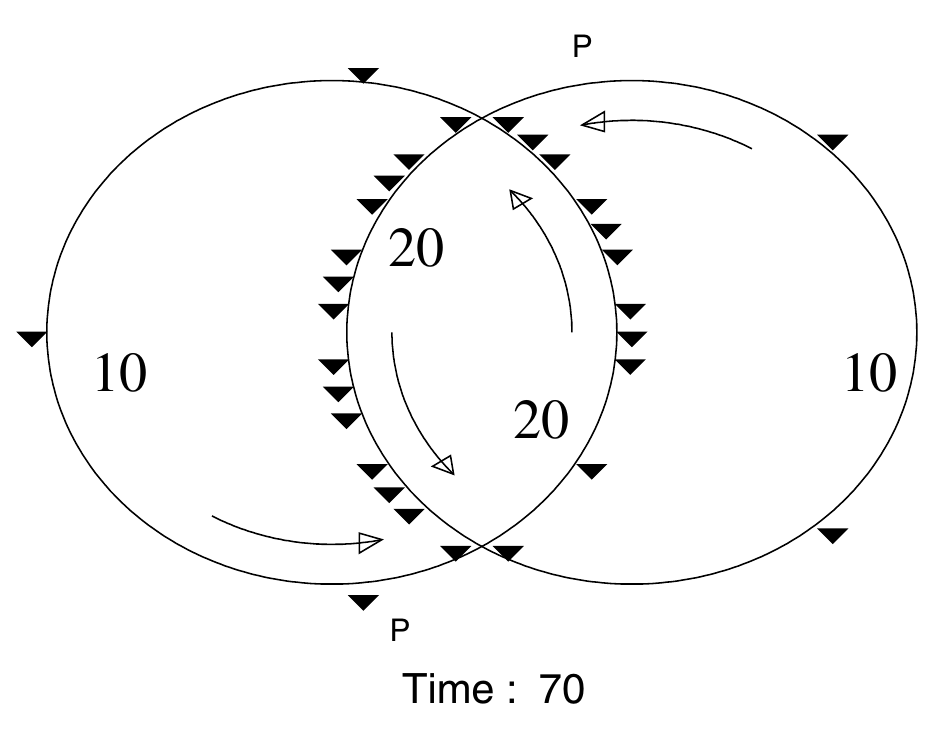}
    \caption{Initial and asymptotic periodic configurations
        during the saturation phase. The size of each priority road is 10.
        The size of each non-priority road is 20. The
        total number of vehicles is 30.}
    \label{ph2-deux}
  \end{center}
\end{figure}     
            
    \item \emph{Recession phase:} This phase appears only when the
      sum of the priority road sizes is less than the sum of the
      non-priority road sizes. In this case, the vehicles continue
      to be accumulated in the non-priority roads by increasing the
      density (with respect to the \emph{saturation phase}). At
      the periodic regime, the number of vehicles on the non-priority
      roads is so big that the flow is less than $1/4$,
      but this number is insufficient to fill up the non-priority roads
      and freeze the traffic (which corresponds to the \emph{freeze phase}, see below).

\begin{figure}[htbp]
  \begin{center}
    \includegraphics[width=4cm]{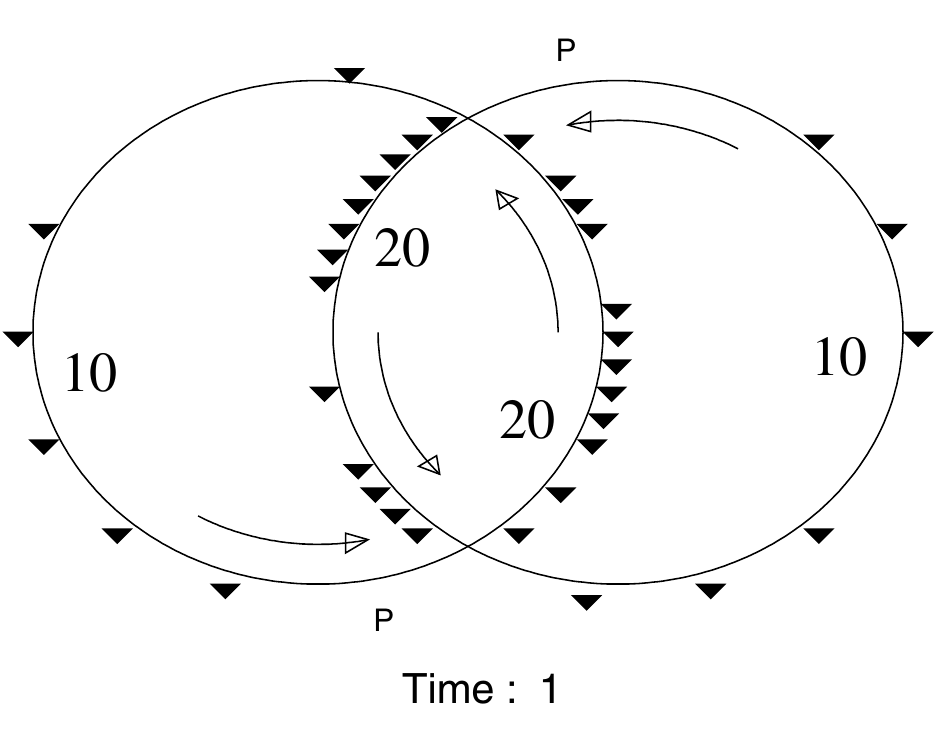} \hspace{0.5cm}
    \includegraphics[width=4cm]{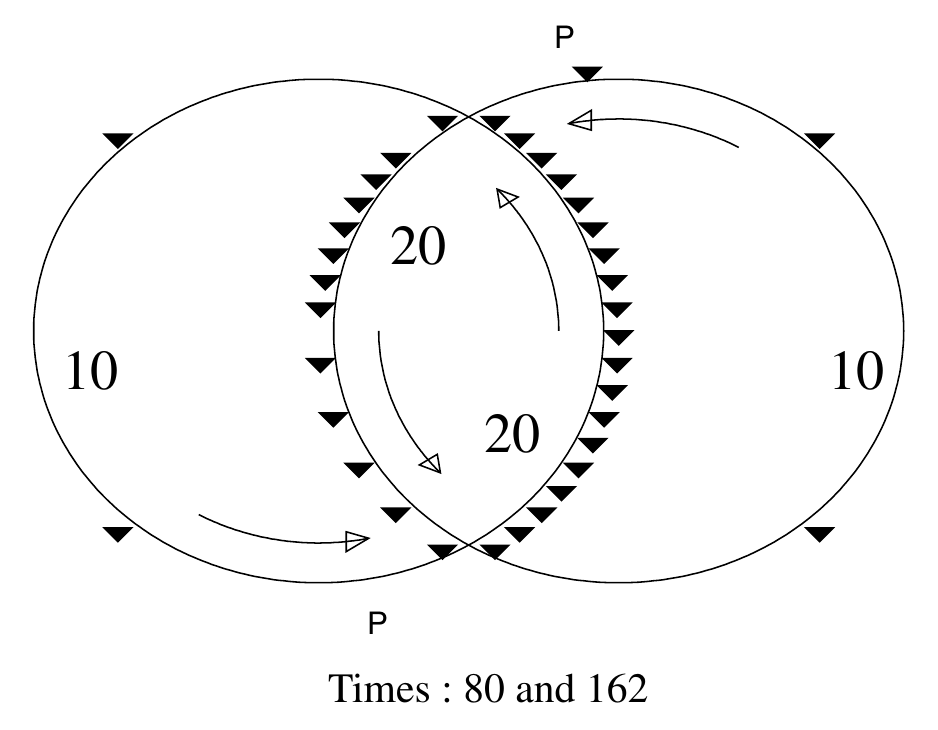} \hspace{0.5cm}
    \includegraphics[width=4cm]{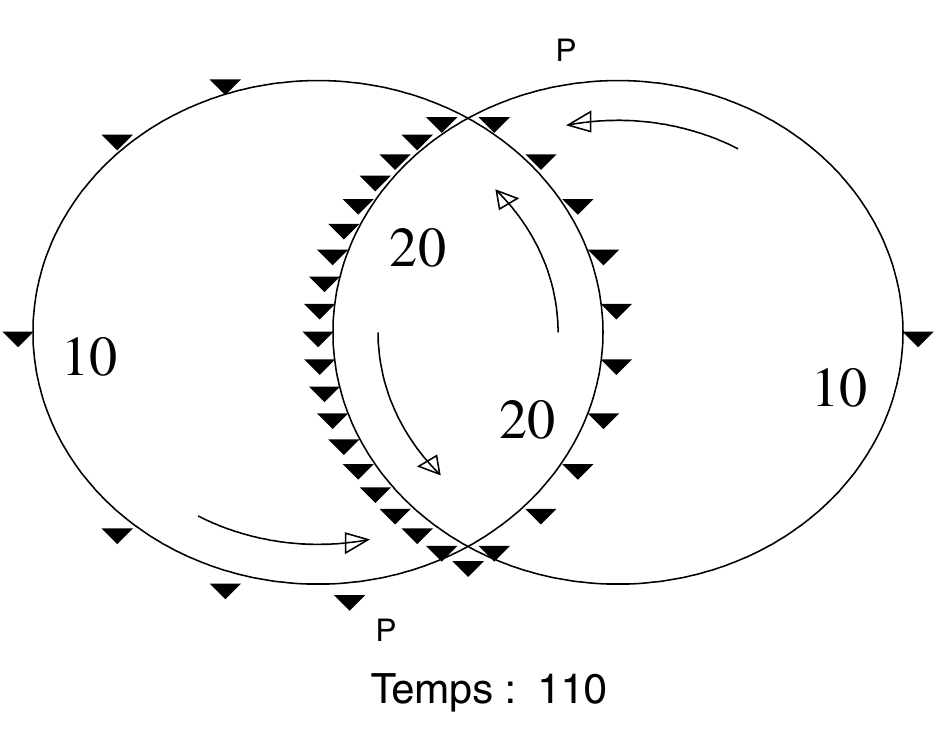}
    \caption{Initial and asymptotic periodic configurations
           during the recession phase. The size of each priority road is 10.
           The size of each non-priority road is 20. The
           total number of vehicles is 37.}
    \label{ph3-deux}
  \end{center}
\end{figure}     

      During this phase, a phenomenon of  pumping vehicles from
      priority roads to non-priority roads is observed. This can
      be seen on the fundamental diagram of the priority roads
      where the density of vehicles on these roads decreases when
      the global density in the whole system is increased
      (Figure~\ref{loi1234}, Roads 2 and 3, Figure~\ref{loi14-23}, right side.).

      Since the density of vehicles in the non-priority roads
      exceeds $3/4$, the global flow is less than $1/4$ (see 1D-traffic phases in
      \citep{PHD,CDC05}). Since the priority roads are served by a
      flow less than $1/4$ and the output flow of these
      roads is bounded by the maximum flow of the
      junctions (which is $1/4$), the density on these roads is
      equal to the global flow.  The global flow is less than $1/4$ and
      decreases by increasing the global density.
      Hence, the density on the priority roads decreases when we increase the global density.
    \item \emph{Freeze phase:} A freeze appears as soon as the number of vehicles
      on the whole system reaches the sum of the sizes of the non-priority roads.
      In this case, the traffic freeze results from the non-priority roads filling up and
      blocking the two junctions; see Figure~\ref{ph4-deux}.

\begin{figure}[htbp]
  \begin{center}
    \includegraphics[width=4cm]{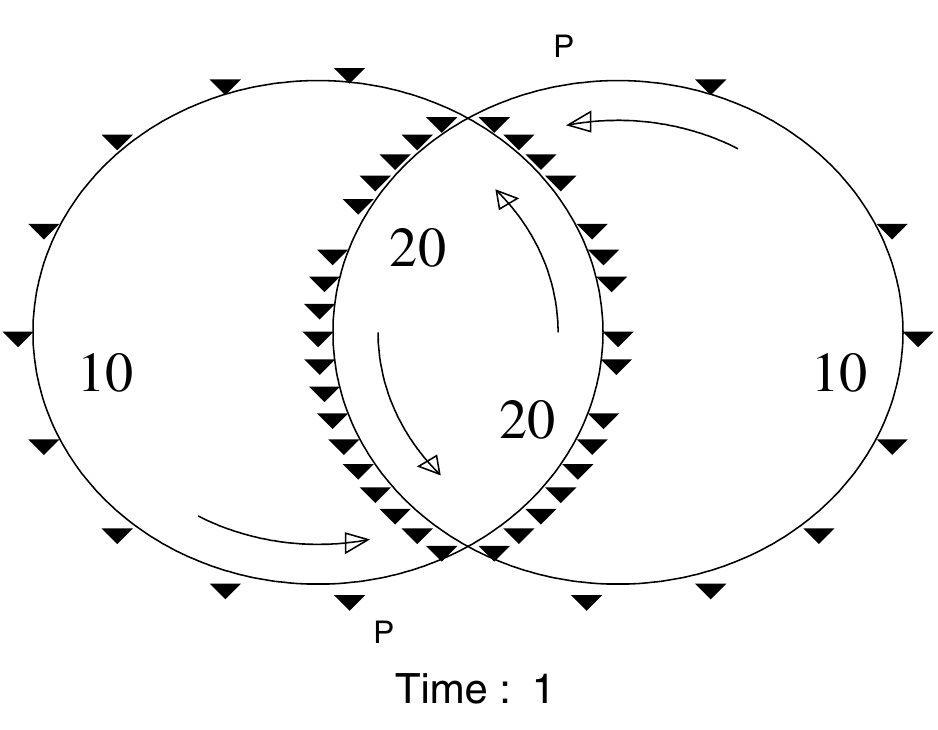} \hspace{1cm}
    \includegraphics[width=4cm]{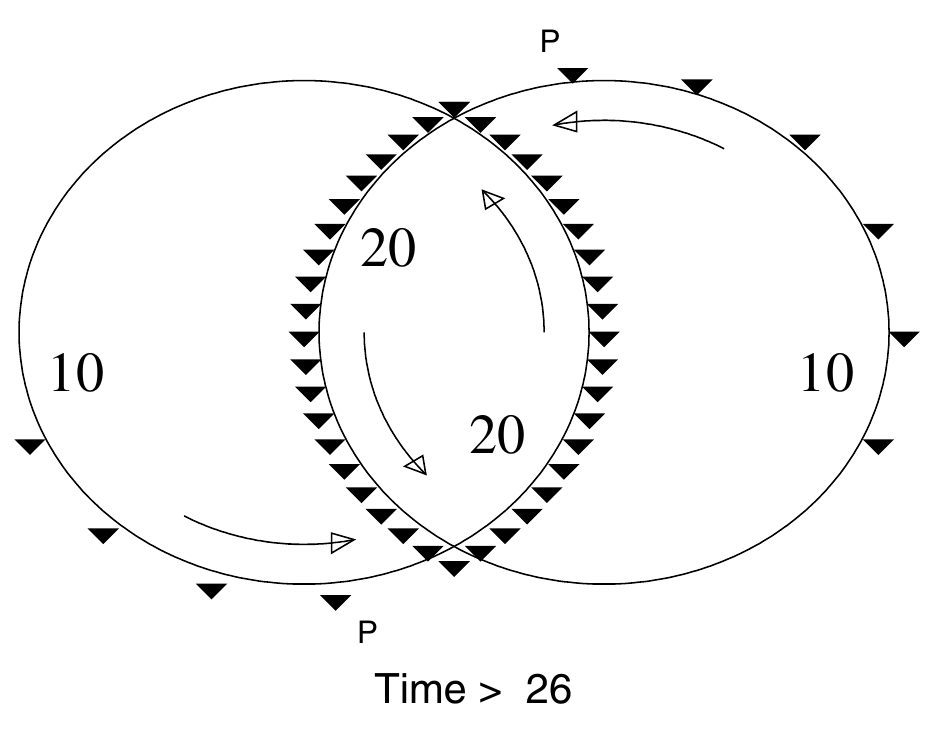}
    \caption{Initial and asymptotic periodic configurations
           during the freeze phase. The size of each priority road is~10.
           The size of each non-priority road is~20. The
           total number of vehicles is~50.}
    \label{ph4-deux}
  \end{center}
\end{figure}     
           
  \end{enumerate}

In Figure~\ref{loi1234}, we give the fundamental diagram for each
road.

%-----------------------------------------
\subsection{A Regular city on a torus}
%-----------------------------------------

\begin{figure}[htbp]
  \begin{center}
    \includegraphics[width=5cm,height=4cm]{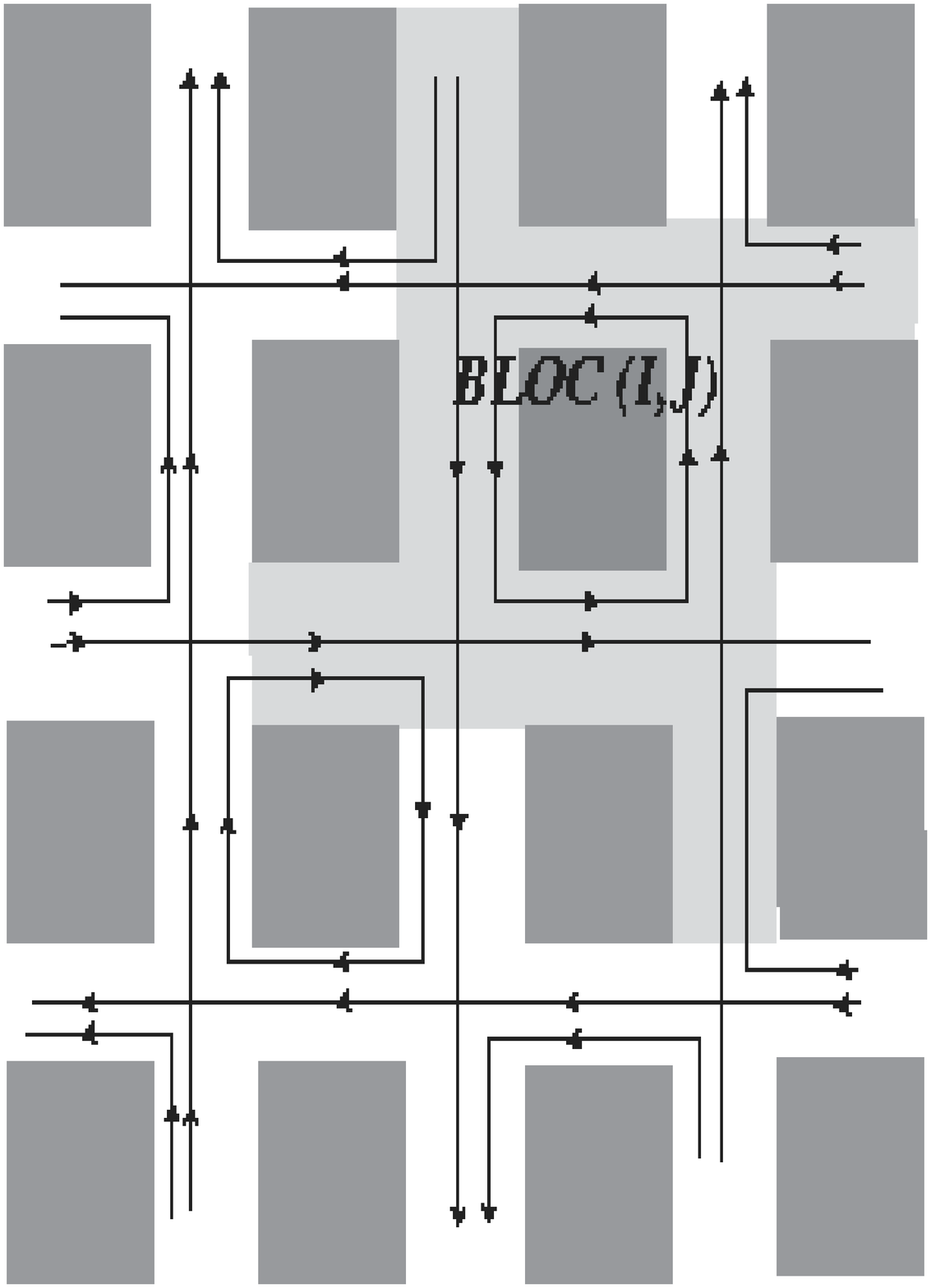}\hspace{1cm}
    \includegraphics[width=4.3cm]{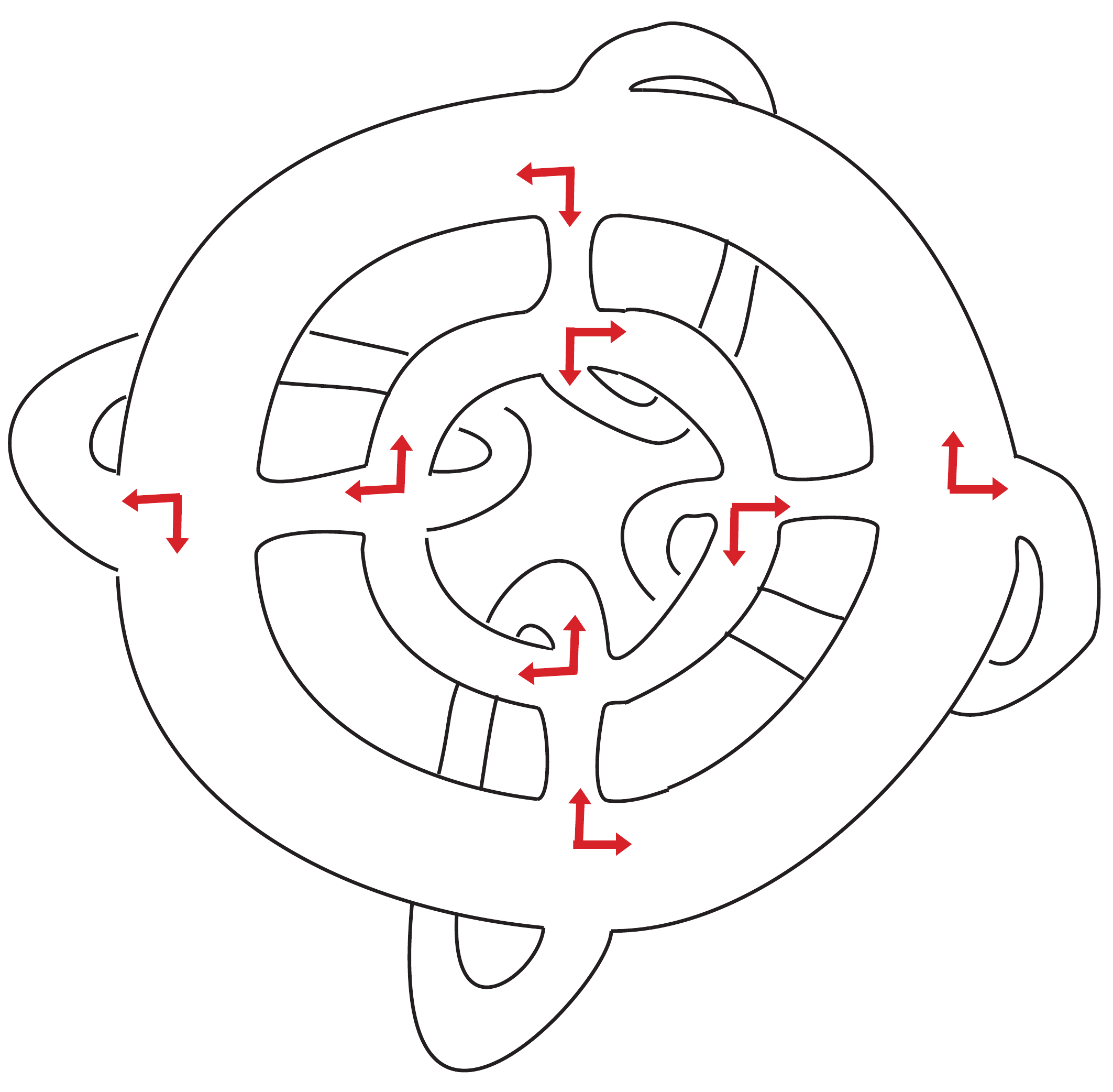}
    \caption{A 4 $\times$ 4 regular city, and a 2 $\times$ 4 regular city on a torus.}
    \label{regular}
  \end{center}
\end{figure}

In this section, we extend the model described in
section~\ref{sec2} to the case of a regular city. To fix the
density, the city is set on a torus; see Figure~\ref{regular}.
We will consider here only the case where all the streets have the
same sizes which catch the main qualitative results (that is the
presence of mainly four phases). It is not easy to obtain the
dynamics of the city but this has been done in a modular way.
Details on the construction are available in \citep{PHD}.
From this dynamics, by numerical
simulation, we see that the average flow mainly does not depend on
the car initial positions, but only on the car density. Therefore
the fundamental diagram associated to the city exists
empirically. We discuss here these numerical results.

As above, we simulate the system for different global densities of
vehicles and we derive the traffic fundamental diagram.
%and we observe the asymptotic periodic regimes. We show
Numerically we observe the existence of the average flow becoming
independent of the initial car distribution when the size of the
system grows. As above, numerical diagrams give two, three, or
four traffic phases depending on the sizes of the roads. We
confirm the observations made in the case of two junctions on the
dependence of this diagram on the ratio between the sum of the
sizes of the non-priority roads and the sum of the sizes of the
priority roads. However, we observed that the beginnings of the
recession and of the freeze phases depend also on the lengths of
road circuits. Indeed, besides the free density phase where the
vehicles move freely, more vehicles circulate on the non-priority
roads during the other phases. As soon as the average density of
vehicles on a circuit made by non-priority roads exceeds 3/4, the
recession phase appears\footnote{Note that a non-priority road
circuit is more susceptible to be filled than a priority-road
circuit. The filling of the latter is non stable, and is
considered as a singular case.}. Similarly, for much higher
densities, as soon as a circuit of non-priority roads fills up,
the traffic is frozen. However, we think that this dependence of
the traffic phases on the lengths of the non-priority road
circuits may disappear when the size of the system grows. To
confirm or reverse this assumption, we shall extend the analytical
results obtained on small systems to larger systems.

\begin{figure}
    \begin{center}
      \begin{tabular}{ccc}
        \includegraphics[width=4cm]{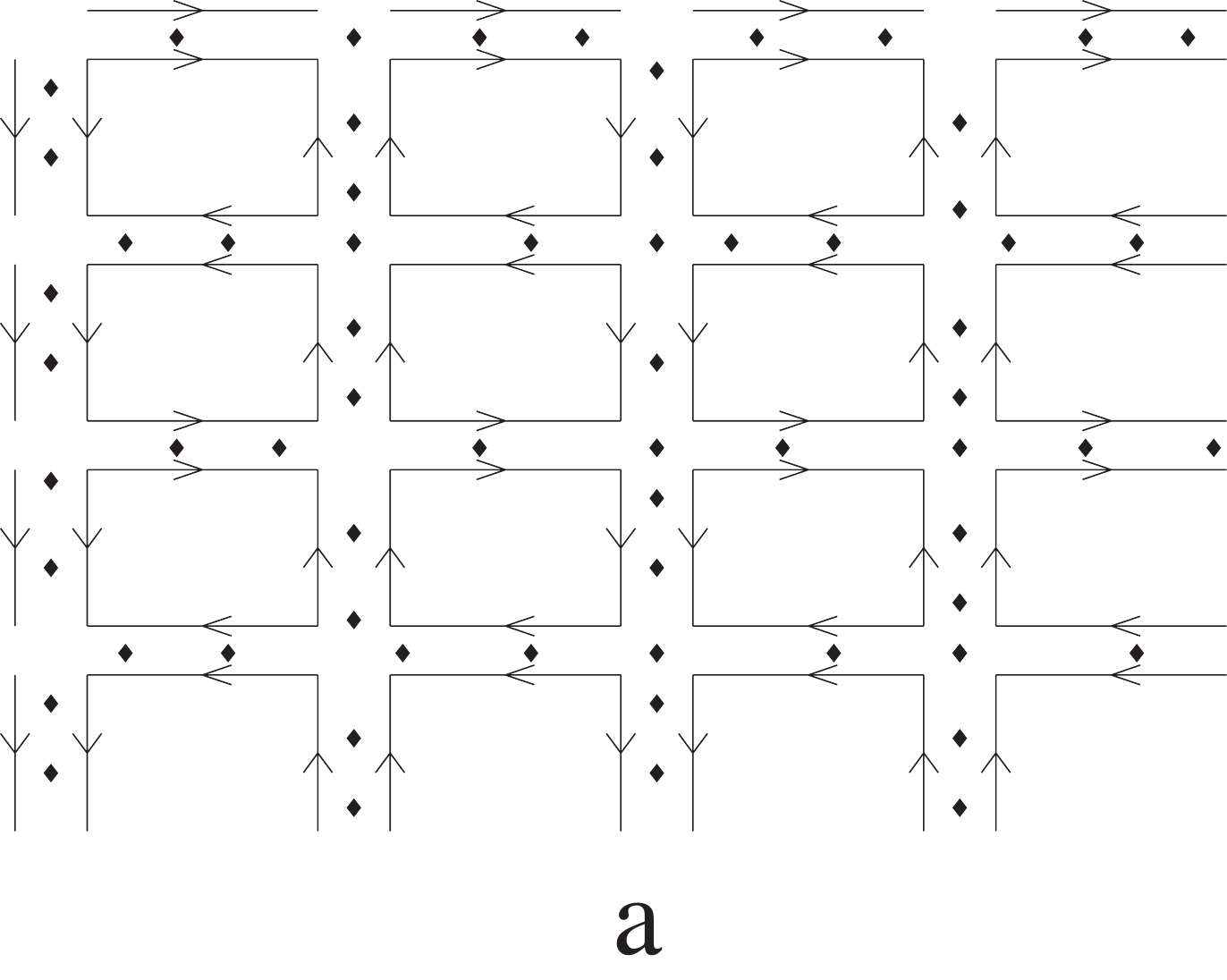} &
        \includegraphics[width=4cm]{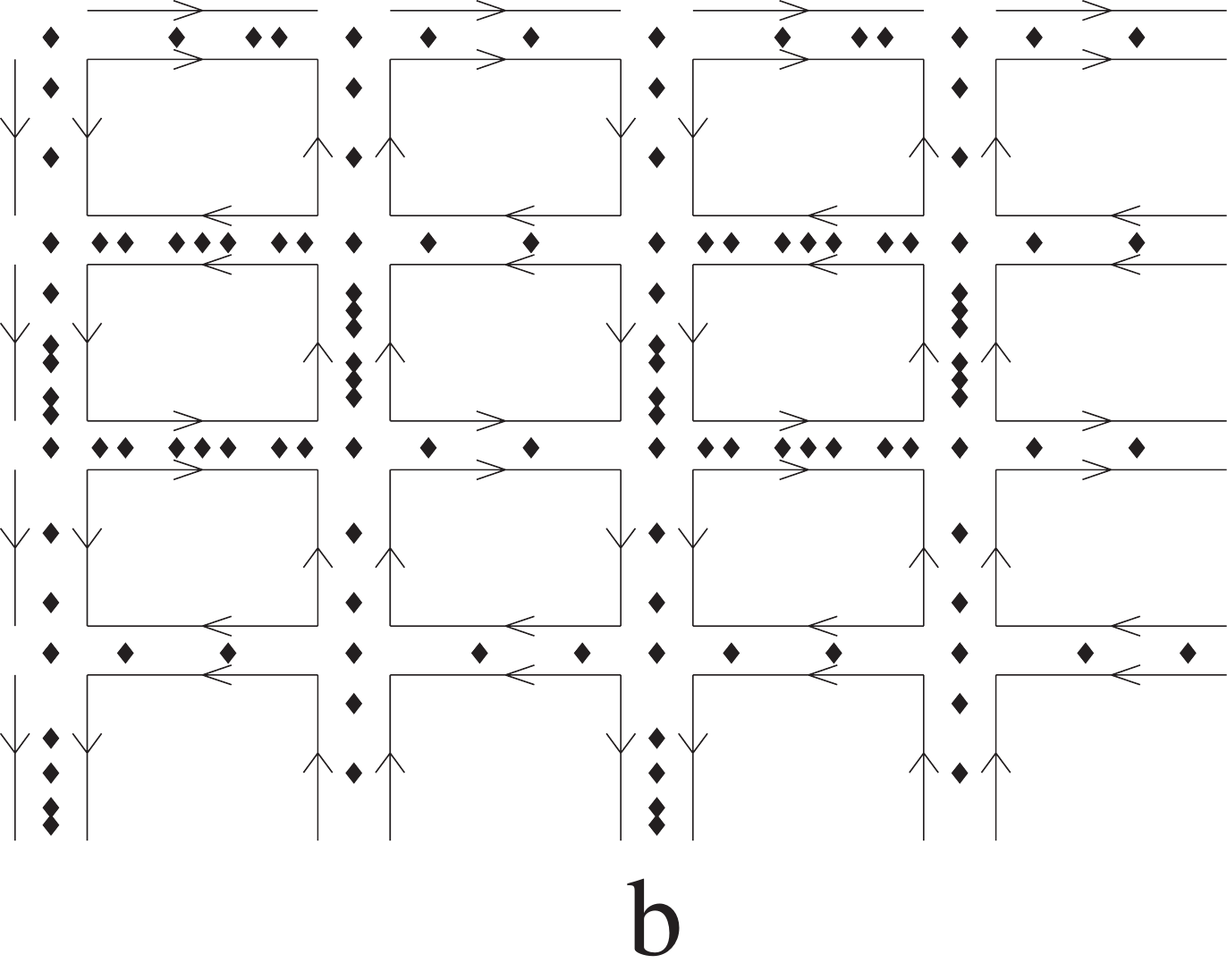} &
        \includegraphics[width=4cm]{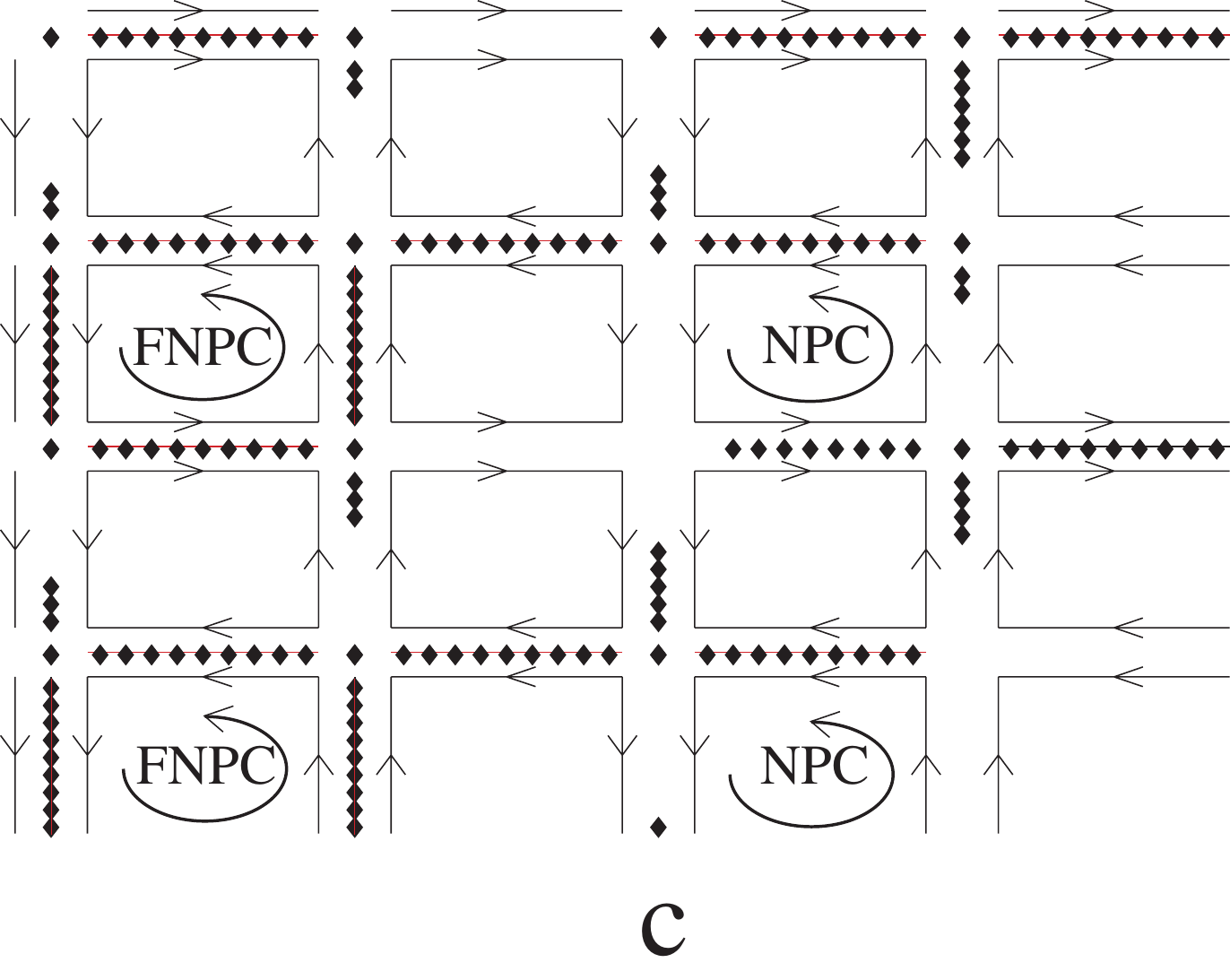}
      \end{tabular}
    \end{center}
    \caption{The periodic regime in the cases: (a) free phase
      ($d=64/304\approx 0.21$), (b) saturation phase ($d=120/304\approx 0.39$),
      and (c) freeze phase ($d=188/304\approx 0.61$). NPC: non-priority road circuit,
      FNPC: full non-priority road circuit.}
    \label{phases}
\end{figure}

In Figure~\ref{phases}, we show the initial and the asymptotic
positions of the vehicles in a $4\times 4$ regular city on a
torus, for the three phases (free, saturation, and freeze phases)
appearing in the symmetric case (where all the roads have the same
size, so $r=1/2$, and where the turning percentages are
all equal to $1/2$). The traffic behavior during the free and the
saturation phases are almost the same as in the case of two
junctions. However, the recession and the freeze phases can appear
for smaller densities comparing to the case of two junctions. Note
that in the case of two junctions only one non-priority road
circuit exists, whose length is also the sum of the sizes of all
the non-priority roads of the system. In the case of a regular
city or a system with more than two junctions, the blocking state
can be reached before all the non-priority roads fill up, that is,
when the global density $d$ satisfies $d<r$.

%--------------------------
\section{Traffic Control}
%--------------------------

To control the traffic in a regular city, we use traffic lights.
We suppose the existence of a traffic light at
the exit of each road, which is also an entry to a junction. The
traffic lights have only green and red colors with the standard
traffic meaning. To construct a regular city controlled
by traffic lights, we follow the same approach as in the case of a
regular city managed with the priority rule. That is, we develop a
model of one junction controlled by a traffic light
(see \citep{PHD} where we have presented many Petri net models with
traffic lights), then we generalize it to a regular city on a
torus (see also \citep{PHD} where the construction of the
controlled regular city is explained.) From the numerical simulation of these models we
see that the fundamental diagram makes sense, that is that the car
initial position does not hava notable influence on the average
flow. The average flow depends only on the car density and on the
light strategy. In this section we discuss the different light
strategies used and their corresponding fundamental diagrams for a
city on a torus.

%-----------------------------------
\subsection{Open loop control}
%-----------------------------------

In applying an open-loop control, one can consider a periodic
traffic plan over time for each junction or a periodic plan for
all the network. Here we consider the latter case. First, we fix
the period called the \emph{cycle}. Then we fix the duration of
time assigned to the green color and the duration of time assigned
to the red color. The timings are chosen with consideration of the
traffic on each road, but they are chosen only once and used
thereafter.

We test this approach on regular cities where all the roads have
the same size, and all the turning proportions are equal to $1/2$.
We fix the cycle to 4 units of time. Indeed, since the car
flow cannot exceed $1/4$ at the junctions (see
subsection~\ref{subsec2-1}), $4$ units of time is the shortest
possible value of the cycle for our model. This cycle time is not
realistic but gives the largest flow. The durations of green and red
colors are fixed to 2 units of time. 

The fundamental diagram obtained with this
policy is shown in Figure~\ref{policies}, where we compare it to
the other policies described below. For the right priority policy
the freeze phase starts at the density~$1/2$.  Therefore the
open-loop control extends the saturation phase and reduces the
freeze phase, the free phase being almost unchanged. 

%--------------------------------------
\subsection{Local feedback control}
%--------------------------------------

The local feedback control depends on the
state of the traffic on the roads.  By \emph{local} feedback we
mean that at every instant and at every junction the control
depends only on the traffic on the roads entering the junction and
not on the other roads.  Let us clarify it. We denote by :
\begin{itemize}
\item $R_1$ and $R_2$ the two roads entering the junction, \item
$n_1$ and $n_2$ their sizes, \item $z_1^k$ and $z_2^k$  the number
of vehicles on the corresponding road at time $k$, \item $b_1^k$
and  $b_2^k$ given by
              $\begin{cases}
                1 & \text{if at time $k$ there is a vehicle intending to enter}\\
                  & \text{the junction, on the corresponding road,} \\
                0 & \text{otherwise.}
              \end{cases}$
\end{itemize}
The feedback control $u^k$ at time $k$ is given by:
$$u^k=\begin{cases}
        \text{green for } R_1 \text{ and red for } R_2 & \text{if } n_2b_1^k+z_1^k\geq n_1b_2^k+z_2^k,\\
        \text{green for } R_2 \text{ and red for } R_1 & \hbox{otherwise.}
       \end{cases}$$

This feedback gives the green color to the road with a
vehicle intending to enter the junction when there is only one
such a road, and to the more crowded road (relatively to its size)
in the other cases. This policy is applied to all the junctions.

The fundamental diagram obtained is given in
Figure~\ref{policies}. We can see that the local feedback control is
clearly better than the open-loop control since the freeze phase range is reduced to
zero without any worsening of the other phases. 

\begin{figure}[htbp]
  \begin{center}
    \includegraphics[width=6cm]{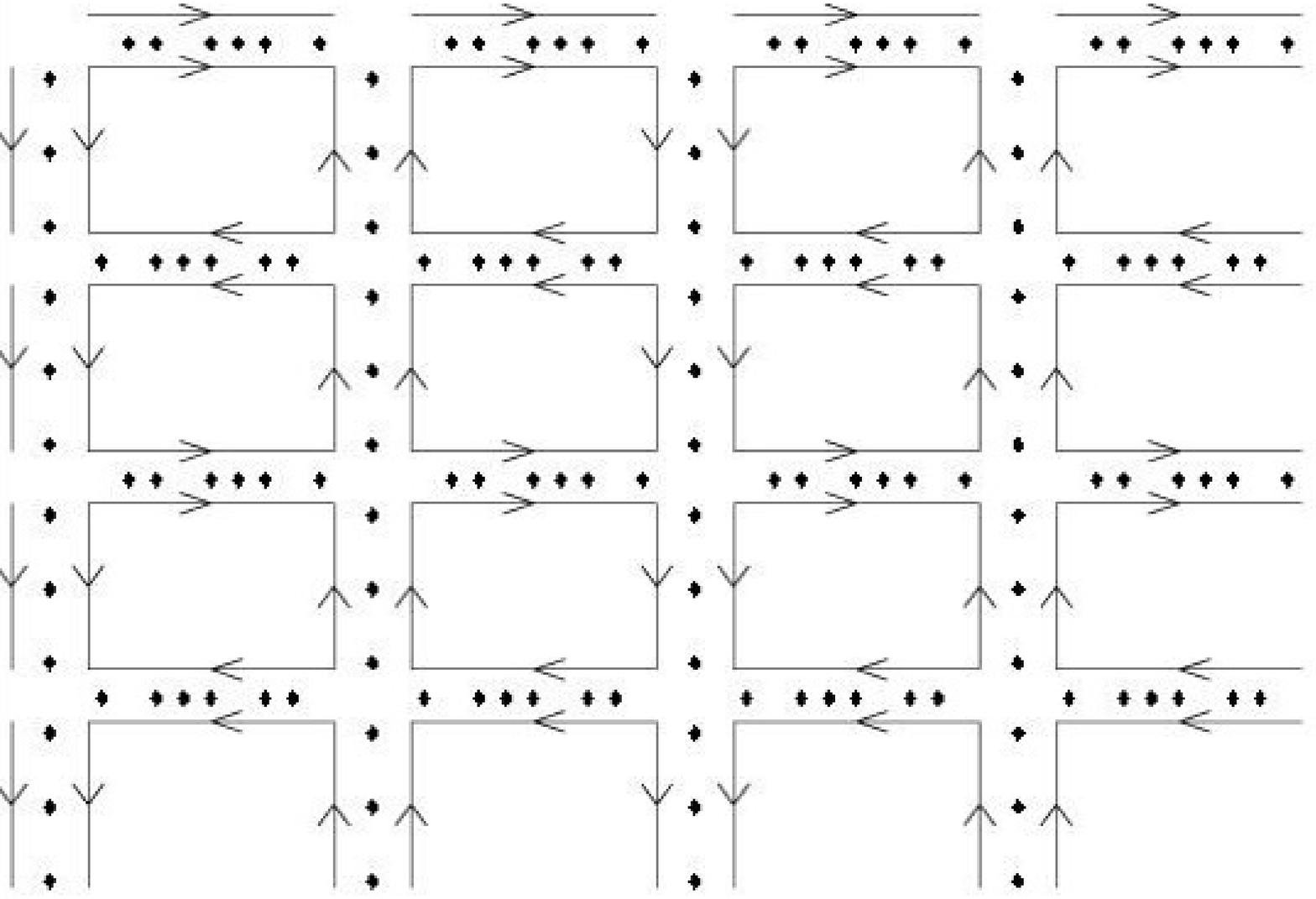}\hspace{0.7cm}
    \includegraphics[width=6cm]{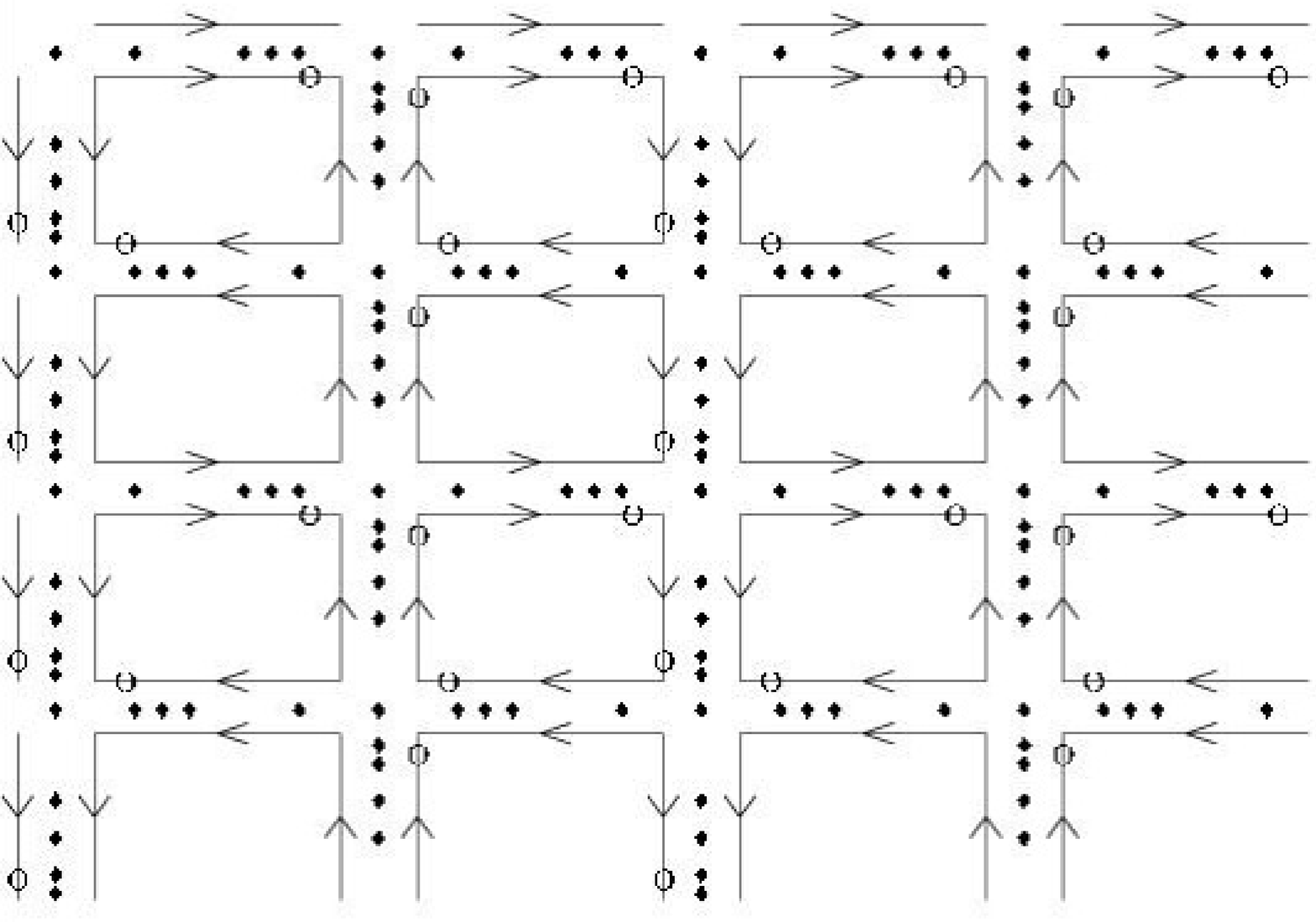}
    \caption{Stationary regimes. On the left side: open loop light control.
        On the right side: local feedback light control.}
    \label{great}
  \end{center}
\end{figure}

In Figure~\ref{great}, we compare open-loop and feedback traffic
controls in terms of the distribution of vehicles at the periodic
regime. We
see that the local feedback induces a stationary regime where the
vehicles are more uniformly distributed (on the roads) comparing
to the open loop control case.

The objective in this example is to show \emph{the existence}
of initial car configurations, for which it is better (in
terms of car distribution on the city) to manage the traffic by
a local feedback control than by an open-loop control. 

%------------------------------------------
\subsection{Global feedback control}
%------------------------------------------

In this section, we follow the point of view developed in TUC
(Traffic Urban Control, see \citep{PAP02}) and use a global
feedback based on a linear quadratic stabilization around a
nominal trajectory of the traffic in a city. In this modeling, the
streets are seen as car inventories flowing from one road to
another. We suppose the existence of a traffic light that we can
control at each ingoing street of the junctions.

At time $k$ for road $i$ we denote $x_i^k$ as the number of
vehicles, $\bar{x}_i$ as the nominal number of vehicles wanted,
$u_i^k$ as the flow outgoing from the road $i$ during the green
phase, and $\bar{u}$ as the nominal flow wanted in this road. We
solve the linear quadratic control problem~:
\begin{align}
    \min_{u\in\CU} & \sum_{k=0}^{+\infty}(x^k-\bar{x})'Q(x^k-\bar{x})+(u^k-\bar{u})'R(u^k-\bar{u})
       \nonumber\\
     & (x^{k+1}-\bar{x})=(x^k-\bar{x})+B(u^k-\bar{u}). \label{dynn}
\end{align}
where $B$ is a matrix describing the interconnections of the
streets (each junction has a proportion to enter the outgoing
street), and $Q$ and $R$ are weight diagonal matrices that we have
to choose empirically to obtain a regulator working in a
satisfactory way. The result of this optimization is a global
feedback used to determine the timing of the traffic light in the
microscopic modeling of the city. The fundamental diagram obtained
with this control is given in Figure~\ref{policies}.

%----------------------------------------------
\subsection{Fundamental diagram comparison of traffic light policies}
%----------------------------------------------

\begin{figure}[htbp]
  \begin{center}
    \includegraphics[width=6cm]{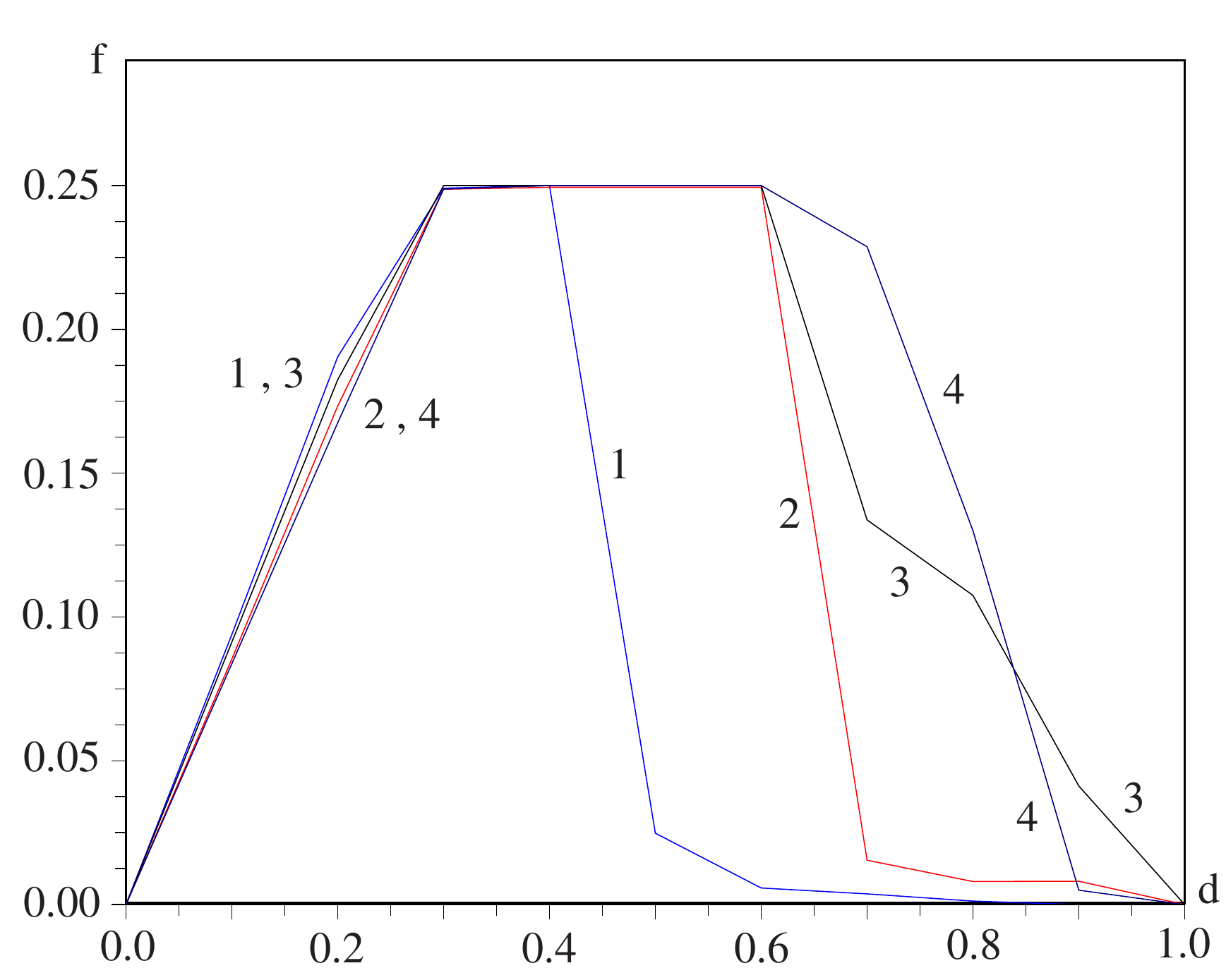}
    \caption{Comparison of traffic control policies on a regular
    city (set on a torus). \textbf{1. }Priority rule. \textbf{2. }Open loop signal light control.
      \textbf{3. }Local feedback signal light control. \textbf{4. }Global feedback signal light control.}
  \label{policies}
  \end{center}
\end{figure}

Figure~\ref{policies} shows the fundamental diagrams obtained with
the four junction policies: -- priority to the right rule (diagram
1), -- open-loop control (diagram 2), -- local feedback control
(diagram 3), and -- global feedback control (diagram 4). For low
densities, the flows are almost the same and are equal to the
density. However, during this phase, the open-loop control policy
(diagram 2) and the global feedback control policy (diagram 4)
slow down slightly the flow with respect to the right policy
control (diagram 1) or to the local feedback control policy
(diagram 3). In the case of open-loop control policy, this can be
explained by the fact that at each junction, a vehicle can be
stopped by a traffic light on a road even if there is no vehicle
intending to move into the junction from the other road. This fact
is also observed in the case of global feedback control. It is due
to the implementation constraint of these feedback control
policies which impose a minimum green and red time.

The critical density at which the freezing appears in the case of
the priority to the
right policy is 1/2. This value represents
also the ratio between the sum of the sizes of the non-priority
roads and the sum of the sizes of the priority roads for a
symmetric city. (All the roads have the same size, thus the ratio
is equal to 1/2.) All the other policies improve the density at
which appear the recession and the freeze phases. The best one is
the global feedback.

The maximal flow, obtained in the saturation phase, corresponds to the saturation
of the junctions which are all the same in this regular city case.

As it can be easily predicted, the flow obtained for low
densities (densities less than $1/4$) and by using any of the
policies, is equal to the car density. Indeed, for low densities,
the vehicles reach asymptotic regimes where they can move freely,
independently of the traffic control policy applied at the
junctions. For densities $d$ in $[0,1/2)$, the fundamental
diagram is the same for all the control policies. 
We have understood from simulations that in the case
of priority to the right rule, we need to increase considerably
the sizes of the roads in order to obtain the value $1/4$ of the
flow for a density less than but very close to $1/2$. On the
diagram~1 of Figure~\ref{policies}, the maximal density for which
we obtained the value $1/4$ of the flow is around $0.4$. It is
possible to obtain a maximal density much closer to $1/2$ by
increasing the sizes of the roads, but long time would be required
to simulate the system.

The result obtained for the feedback control policies is
surprising. The global feedback always should be better than the
local feedback, but the implementation is different. In the global
feedback case, we impose a light cycle that we do not impose for
the local feedback case. This is the reason why the local feedback
seems better than the global one at the end of the recession
phase. Up to this consideration the global feedback is the best
strategy. The improvement with respect to the open-loop appears
only for high densities. All these remarks are quite natural but
we see that the light control\;--even with the open-loop
policy--\;significantly improves the diagram with respect to the
priority to the right managing without large deterioration at low
densities.

%-----------------------------------------
\subsection{Response time comparison}
%-----------------------------------------

\begin{figure}[htbp]
  \begin{center}
    \includegraphics[width=6cm]{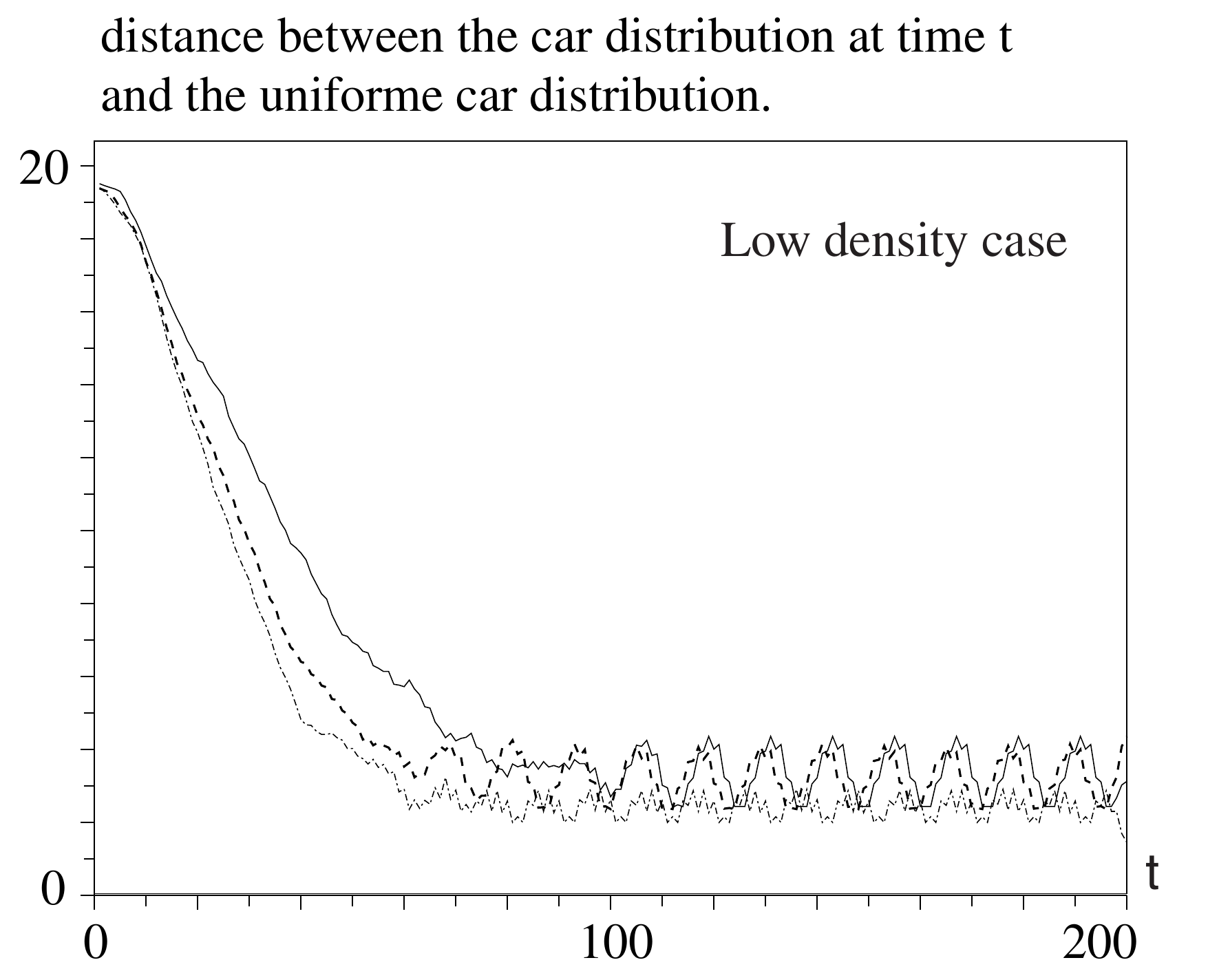}\hspace{1cm}
    \includegraphics[width=6cm]{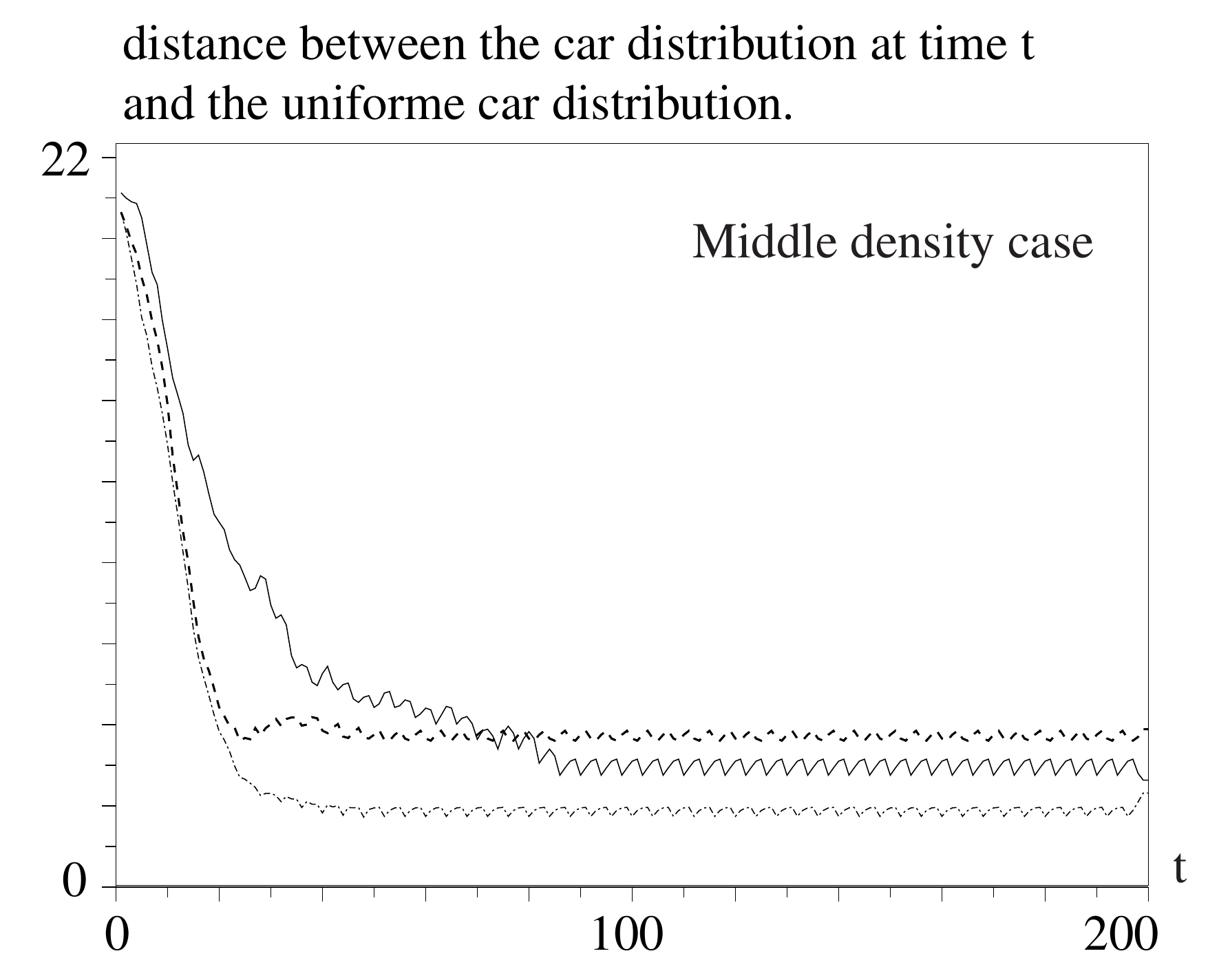}
    \caption{Distance of the vehicle distribution to the uniform distribution
      as function of time. On the left side: a low density
      case. On the right side: a middle density case.
      Continuous line: open loop control. Bold discontinuous line:
      local feedback control. Alternated discontinuous line:
      global feedback control.}
    \label{reponse1}
  \end{center}
\end{figure}

Figure~\ref{reponse1} gives the response times (i.e. time to recover
the stationary or periodic regime after a disturbance) obtained for
the open-loop, local and global feedback policies. In this figure,
we plot the distance between the present vehicle distribution and a
uniform car distribution on the streets as a time function. The city
has four horizontal and four vertical roads.

Independently from the density, the feedback control gives a better
time of response than the open loop control. Both feedback
controls are quite similar in terms of time response, but the
global feedback asymptotic regime is closer to the uniform
distribution.

\begin{remark}
  A regular city can be closed in different ways. We think that the way
  the network is closed will not have impact on the qualitative results
  obtained here. We would have the same traffic phases.
  In the case where junctions are managed with the priority rule, non priority road circuits
  would play the same role on the beginning of the freeze phase, independently
  of the network topology.
\end{remark}

%---------------------------------------------
\section{Junction Design Improvement}
%---------------------------------------------

As explained
above in section~\ref{sec2}, the flow is limited by the maximal
junction flow which is equal to $1/4$. It is possible to improve
the traffic by adding one car place in the junctions. It is easy
to check that with this improvement, the bound will be $1/2$
instead of $1/4$.

Let us give the analytical result obtained for the elementary
system of one circular road taking
shape of the numeral eight (8), with one junction. The only
change on the model is on the constraints~(\ref{const}):
\begin{equation}
  \begin{cases}
    0\leq a_i\leq 1,\;\; 1\leq i\leq n+m\;,\\
    0\leq a_n+a_{n+m}\leq 1.
  \end{cases}
\text{ becomes  }
  \begin{cases}
    0\leq a_i\leq 1,\;\; 1\leq i\leq n+m\;,\\
    0\leq a_n+a_{n+m}\leq 2.
  \end{cases}
\end{equation}

In this new case, we are able to obtain
results analogous to Theorem~\ref{eigen} and Corollary~\ref{fond-cor2}.

\begin{theorem}\label{eigen2}
  There exists an additive eigenvalue $\lambda$ satisfying:
  $$0=\max\left\{\min \left\{d-(1+\rho)\lambda,\;
      r-d-\left(2r-1+\rho\right)\lambda\right\},\;-\lambda \right\}.$$
\end{theorem}

\begin{corollary}\label{fond-cor22}
  For large values of $n$ and $m$ such that $n>m-2$ (which is the case $r\geq 1/2$),
  a non negative eigenvalue $\lambda$ is given by:
  $$\lambda=\max\left\{\min\left\{d,\;\frac{r-d}{2r-1}\right\},0\right\}.$$
\end{corollary}

\begin{figure}[htbp]
  \begin{center}
    \includegraphics[width=4cm,height=3.5cm]{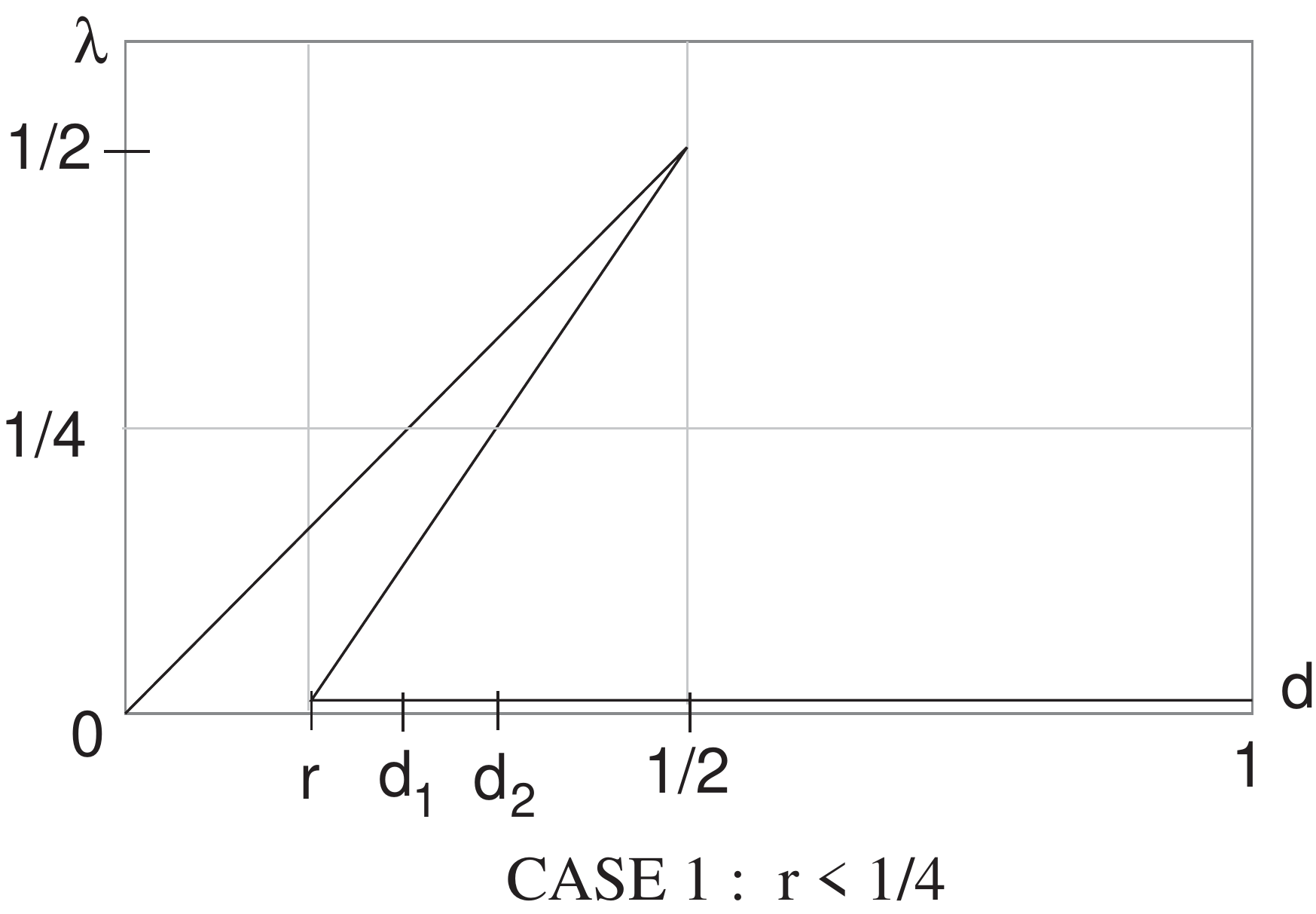}\hspace{5mm}
    \includegraphics[width=4cm,height=3.5cm]{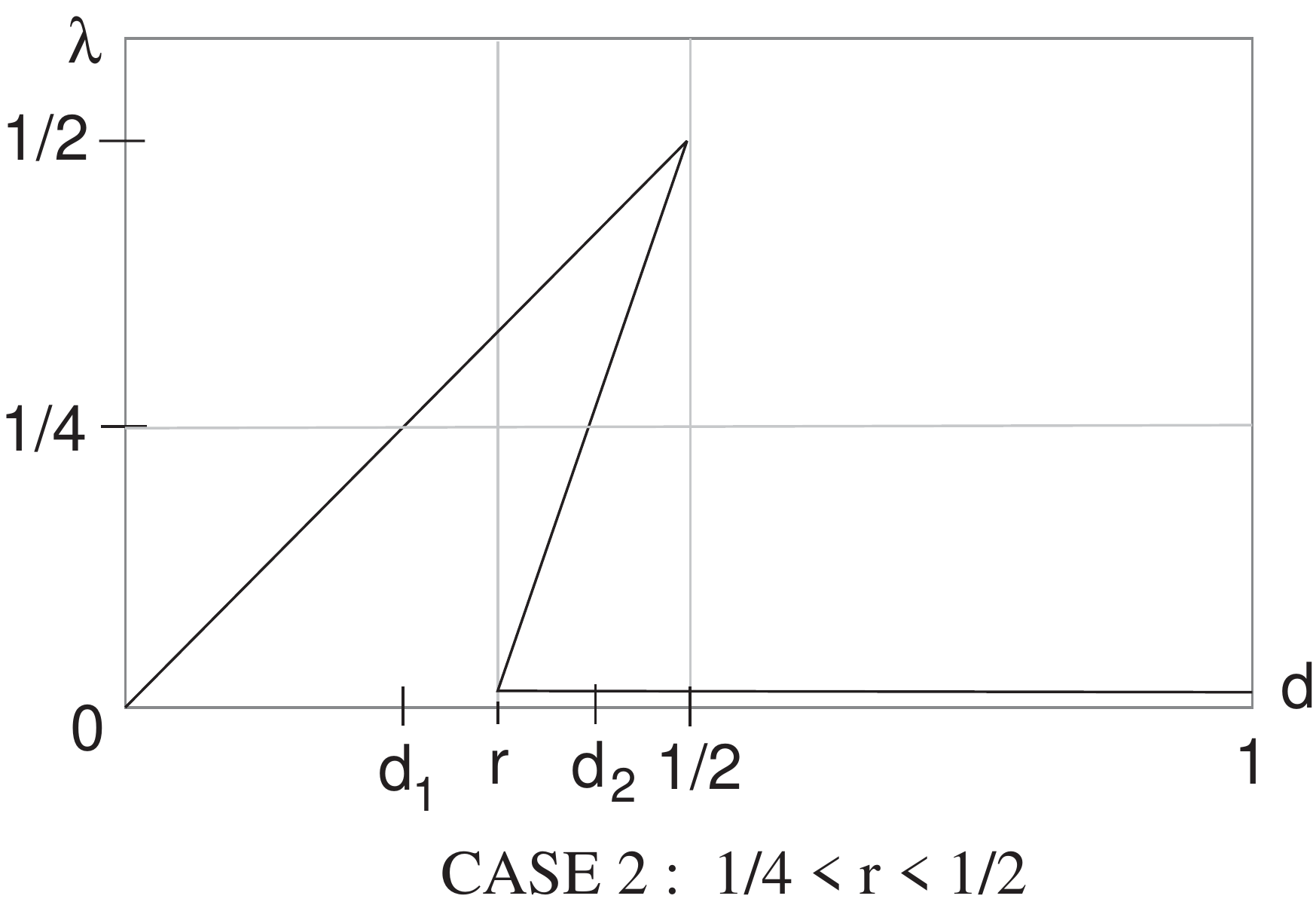}\hspace{5mm}
    \includegraphics[width=4cm,height=3.5cm]{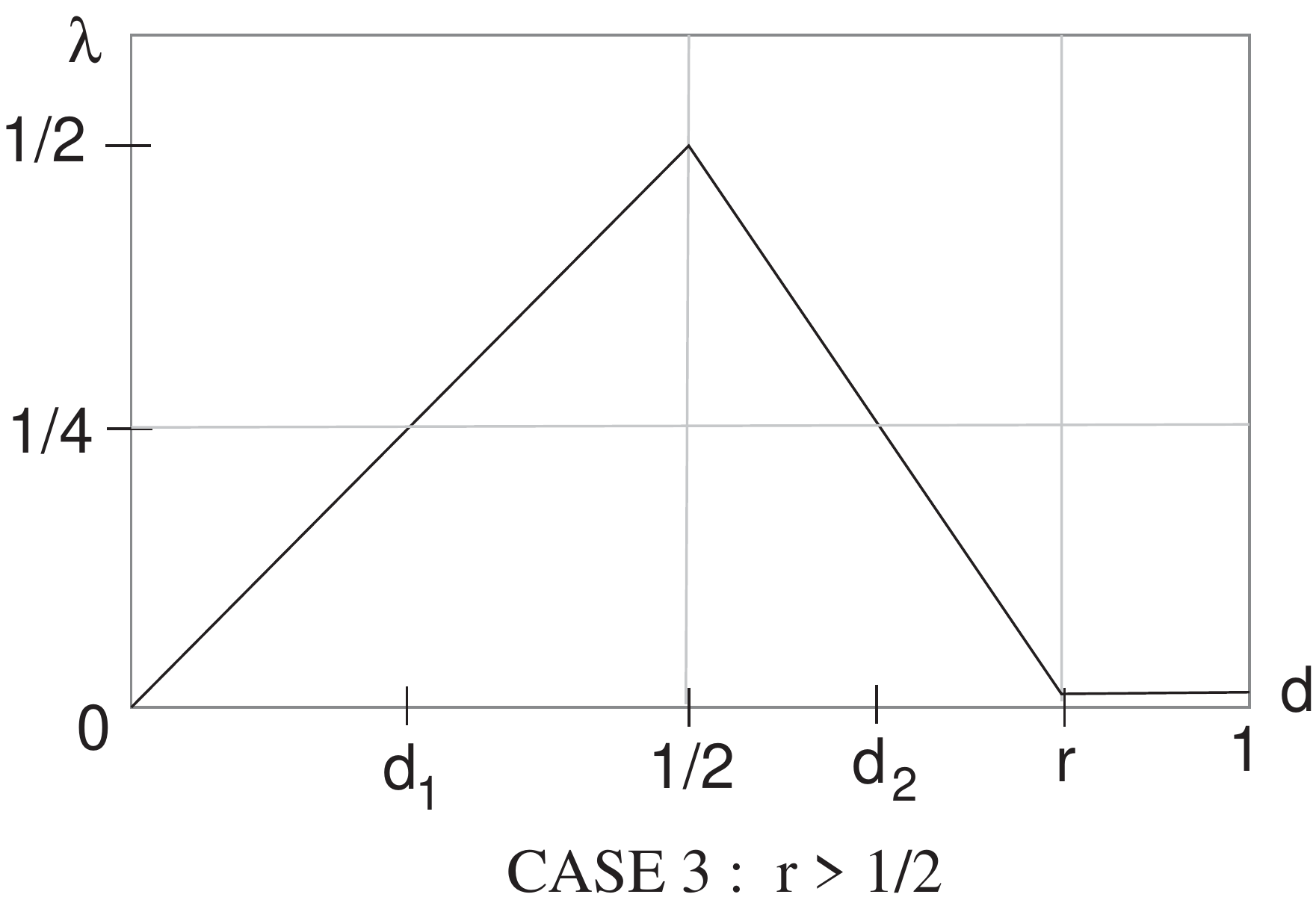}
    \caption{The curve of $\lambda$ given in Theorem~\ref{eigen2} depending on $d$, in the case of a large junction.}
    \label{diaggg}
  \end{center}
\end{figure}

\begin{figure}[htbp]
  \begin{center}
    \includegraphics[width=4cm,height=3.5cm]{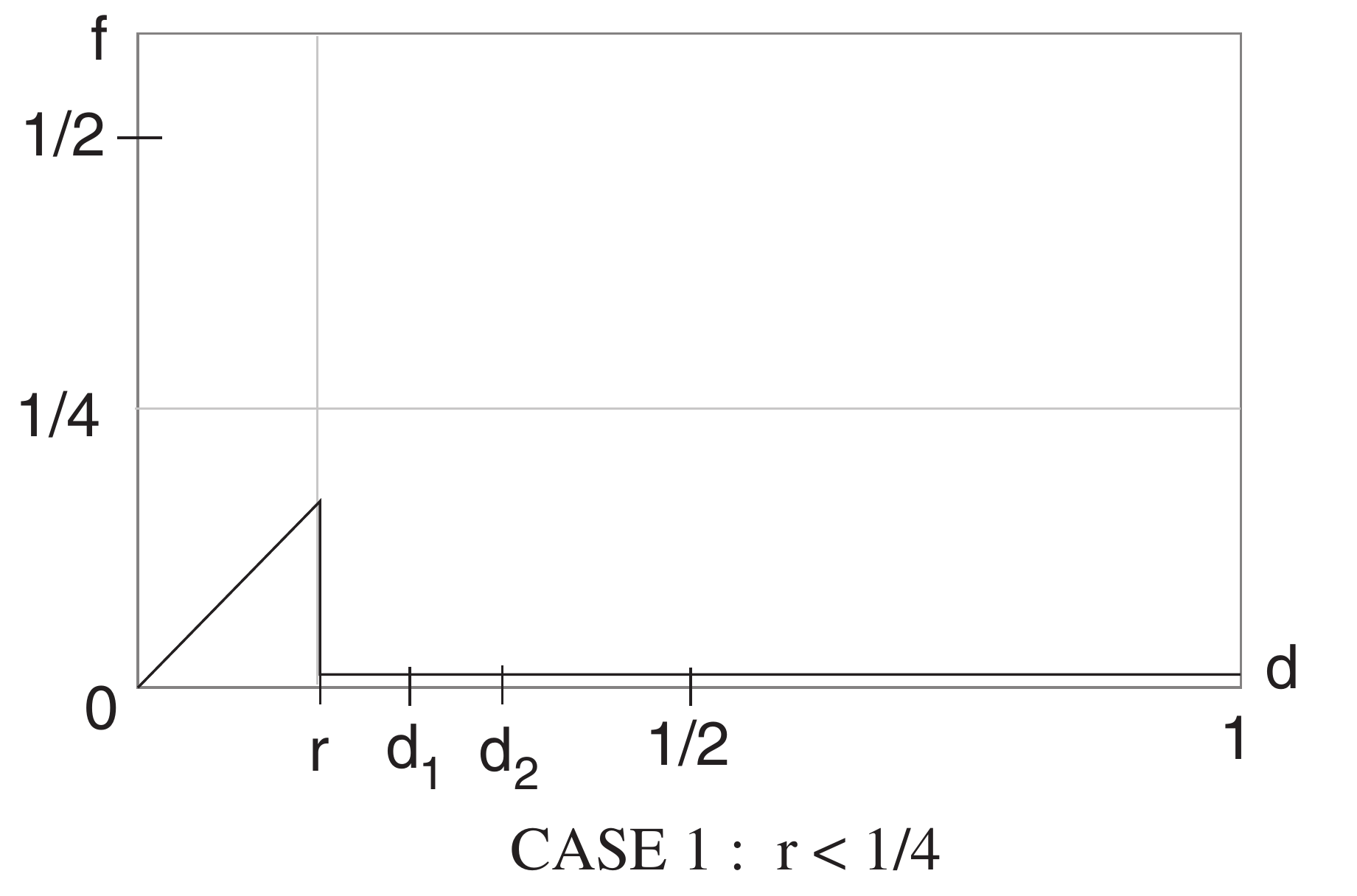}\hspace{5mm}
    \includegraphics[width=4cm,height=3.5cm]{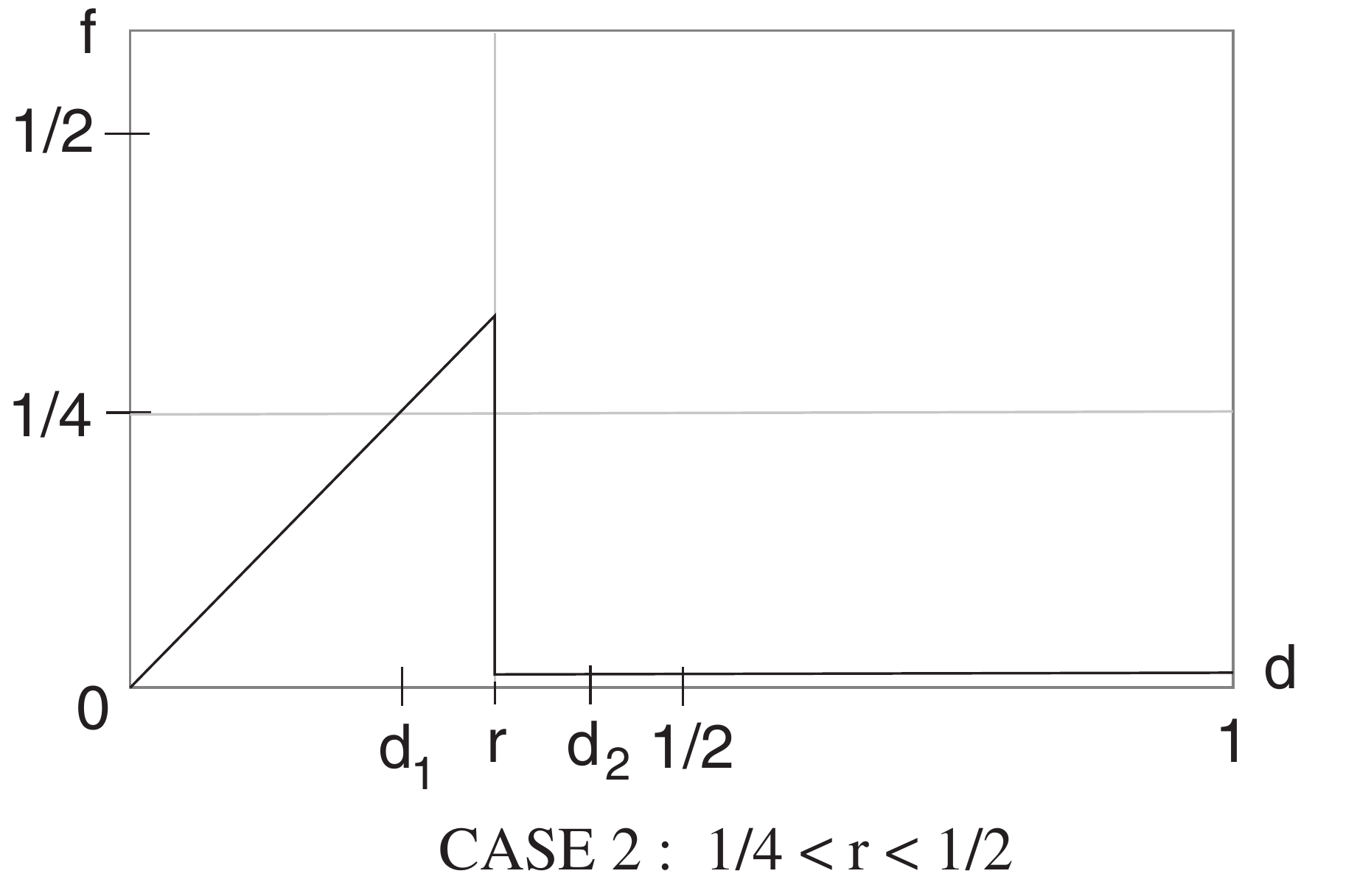}\hspace{5mm}
    \includegraphics[width=4cm,height=3.5cm]{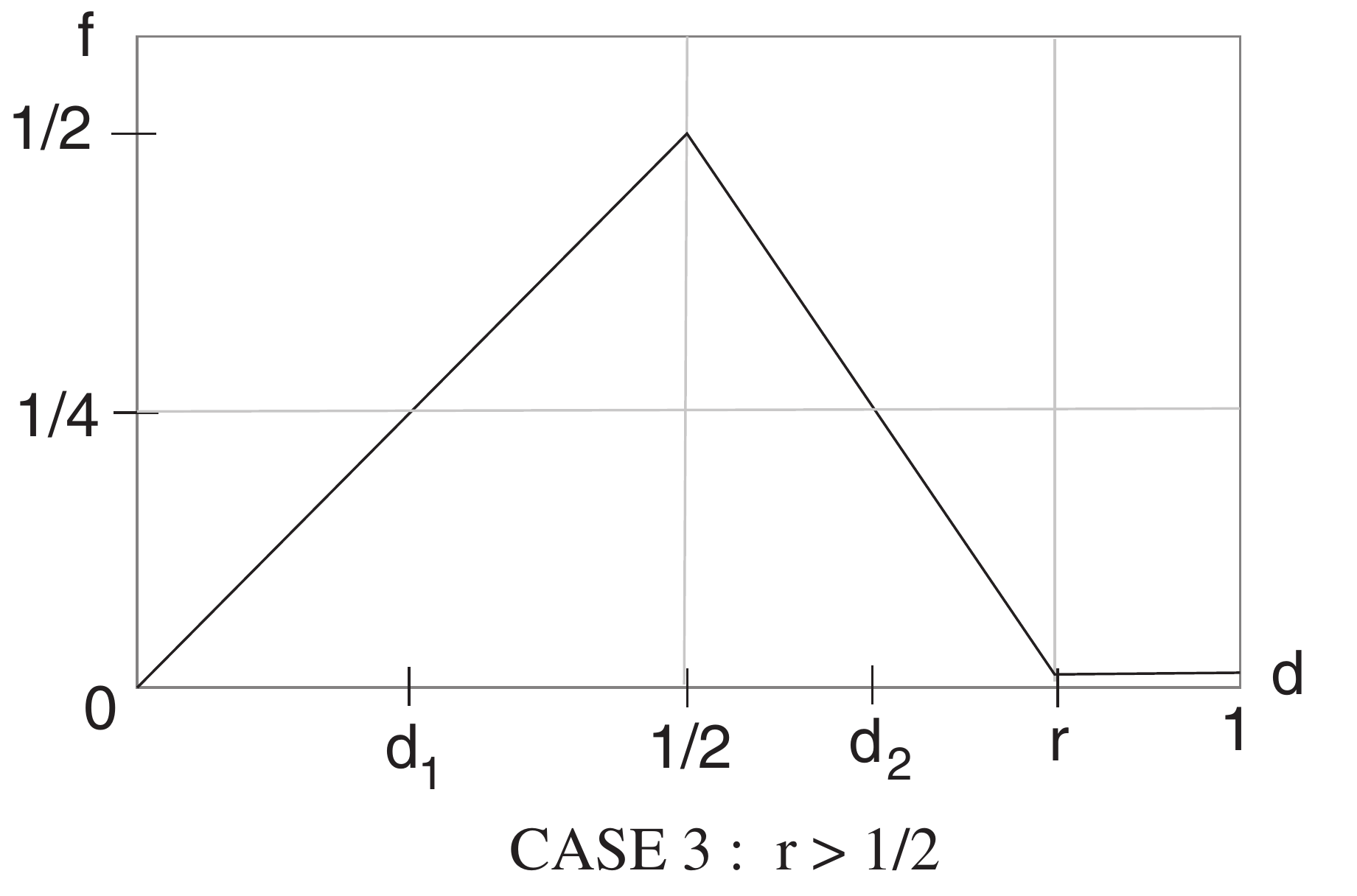}
    \caption{Summary of the fundamental diagrams obtained numerically in the case of a large junction.}
    \label{diag2-num2}
  \end{center}
\end{figure}

From this result, we obtain a good approximation of the
fundamental diagram~:
\begin{equation}
  f=\max\left\{\min\left\{d\;,\frac{r-d}{\max\{2r-1,\;0\}}\right\},0\right\},
\end{equation}
which is shown  in Figure~\ref{diag2-num2}.

%--------------------------
\section{Conclusion}
%--------------------------

The main contribution of this paper is the discussion of the
phases that can appear in the traffic on a road network. We
discuss in detail the four traffic phases (free, saturation,
recession and freeze phases) first on elementary systems and then
on larger road systems. We identify a key parameter that is a
ratio between the sum of the sizes of the non-priority roads and
the sum of the sizes of the priority roads. We discuss how the
freeze phase appears when a circuit of non-priority road jam is
created. On the fundamental diagram, we study the improvement
obtained thanks to open loop, local and global feedback traffic
light policies. The importance of the junction capacity
optimization is revealed clearly by this diagram.

The analytic results are obtained only in the case of two roads
with one junction. It would be interesting to analyze the simplest
case of two circuits of non-priority roads of different sizes in
greater detail to find out if there exists a new phenomenon.

More realistic models should confirm the very natural qualitative
results obtained here, but it would be worthwhile to see them on a
large set of more sophisticated simulations.

The fundamental diagram should have sense also in the case of open
networks for zones of a system where the car density stay almost
fixed. The fundamental diagram must be considered, for traffic, as
the analogue of the ideal gas law. It should be useful to build 2D
fluid macroscopic traffic models based on this law but this
possibility has not been explored for the moment.

\bibliographystyle{elsarticle-harv}
\bibliography{biblio}

\end{document}